\documentclass[
	framed_theorems,
	framed_proofs,
	abstract=true,
	print,
	english,
]{OCPreprint}

\usepackage{pgfplots}
\pgfplotsset{compat=1.16}
\pgfplotsset{plot coordinates/math parser=false}
\newlength\figureheight
\newlength\figurewidth

\usepackage{subcaption}

\usepackage{enumitem}
\setlist[enumerate]{label=(\alph*)}

\let\cite\parencite

\allowdisplaybreaks

\usepackage{marginnote}

\DeclareNameAlias{sortname}{family-given}
\DeclareNameAlias{default}{family-given}

\newcommand\Bad{Q}
\newcommand\Pen{P}

\newcommand\vol{\operatorname{vol}}

\newcommand{\byu}{(\beta,y,u)}

\newcommand{\bybub}{(\beta,y_{\beta},u_{\beta})}

\newcommand{\byuhat}{(\hat{\beta},\hat{y},\hat{u})}
\newcommand{\byubar}{(\bar{\beta},\bar{y},\bar{u})}
\newcommand{\byutilde}{(\tilde{\beta},\tilde{y},\tilde{u})}

\newcommand{\byugamma}{(\beta_{\gamma},y_{\gamma},u_{\gamma})}
\newcommand{\byugammae}{(\beta_{\gamma_{1}},y_{\gamma_{1}},u_{\gamma_{1}})}
\newcommand{\byugammaz}{(\beta_{\gamma_{2}},y_{\gamma_{2}},u_{\gamma_{2}})}

\newcommand{\byugammaT}{(\beta_{\gamma_T},y_{\gamma_T},u_{\gamma_T})}

\newcommand{\byugammabarT}{(\beta_{\bar{\gamma}_T},y_{\bar{\gamma}_T},u_{\bar{\gamma}_T})}
\newcommand{\byugammaubarT}{(\beta_{\text{\b{$\gamma$}}{}_T},y_{\text{\b{$\gamma$}}{}_T},u_{\text{\b{$\gamma$}}{}_T})}

\newcommand{\byuk}{(\beta_{k},y_{k},u_{k})}

\newcommand{\byuT}{(\beta_{T},y_{T},u_{T})}

\usepackage{float}
\floatstyle{ruled}
\newfloat{algorithm}{tbp}{loa}
\providecommand{\algorithmname}{Algorithm}
\floatname{algorithm}{\protect\algorithmname}

\date{\today}
\publishers{}

\title[Finding global solutions of IOCs]{%
	Finding global solutions of some inverse optimal control problems 
	using penalization and semismooth Newton methods%
	\footnote{
		This research was supported by the German Research Foundation (DFG) under grant number WA 3636/4-2
		within the priority program ``Non-smooth and Complementarity-based Distributed Parameter
		Systems: Simulation and Hierarchical Optimization'' (SPP 1962).%
	}%
}

\begin{document}
\author[Friedemann, Harder, Wachsmuth]{%
	Markus Friedemann\footnote{%
		Brandenburgische Technische Universität Cottbus-Senftenberg,
		Institute of Mathematics,
		03046 Cott\-bus, Germany,
		\url{https://www.b-tu.de/fg-optimale-steuerung},
		\email{markus.friedemann@b-tu.de}
	},
	Felix Harder\footnote{%
		Brandenburgische Technische Universität Cottbus-Senftenberg,
		Institute of Mathematics,
		03046 Cott\-bus, Germany,
		\url{https://www.b-tu.de/fg-optimale-steuerung},
		\email{felix.harder.0@gmail.com}
	},
	Gerd Wachsmuth\footnote{%
		Brandenburgische Technische Universität Cottbus-Senftenberg,
		Institute of Mathematics,
		03046 Cott\-bus, Germany,
		\url{https://www.b-tu.de/fg-optimale-steuerung},
		\email{gerd.wachsmuth@b-tu.de},
		\orcid{0000-0002-3098-1503}%
}}

\maketitle

\begin{abstract}
	We present a method 
	to solve a special class 
	of parameter identification problems
	for an elliptic optimal control problem
	to global optimality. 
	The bilevel problem is reformulated via
	the optimal-value function of the lower-level problem.
	The reformulated problem 
	is nonconvex and standard regularity conditions 
	like Robinson's CQ are violated. 
	Via a relaxation of the constraints,
	the problem can be decomposed
	into a family of convex problems
	and this is the basis for
	a solution algorithm.
	The convergence properties are analyzed.
	It is shown that a penalty method can
	be employed to solve this family of problems
	while maintaining convergence speed.
	For an example problem, the use of the
	identity as penalty function allows
	for the solution by a semismooth Newton method.
	Numerical results are presented.
	Difficulties and limitations 
	of our approach to solve a nonconvex problem to
	global optimality are discussed.
\end{abstract}


\begin{msc}	
	\mscLink{49M20}, 
	\mscLink{49M15},
	\mscLink{49N45},
	\mscLink{90C26}
\end{msc}

\section{Introduction}
\label{sec:introduction}

In this paper we study an inverse problem in which we aim to identify finitely
many parameters of an optimal control problem with a linear
partial differential equation.
This results in an infinite-dimensional bilevel optimal control problem.
The concept of bilevel optimization is discussed in
\cite{Bard1998,%
Dempe2002,%
DempeKalashnikovPerezValdesKalashnykova2015,%
ShimizuIshizukaBard1997},
while
\cite{HinzePinnauUlbrichUlbrich2009,%
LewisVrabieSyrmos2012,%
Troeltzsch2009,%
Troutman1996}
present a comprehensive introduction to optimal control. 
Bilevel optimal control problems are also studied in
\cite{FischLenzHolzapfelSachs2012,%
Hatz2014,%
KalashnikovBenitaMehlitz2015,%
KnauerBueskens2010},
for example.
To be more precise, we consider the parametric optimization problem
\begin{equation*}
	\tag{LL$(\beta)$}
	\label{eq:LL}
	\begin{aligned}
		\min_{y\in Y,\,u\in U} 
		\quad &f(\beta, y, u ) \\
		\text{s.t.} \quad & A y - B u = 0, \\
		\quad & u \in \Uad,
	\end{aligned}
\end{equation*}
where $\beta\in \Bad\subset\R^n$ is a parameter,
and the sets $\Bad$, $\Uad$, the linear operators $A$, $B$,
the spaces $U$, $Y$,
and the function $f$ are such that \cref{assumptions} is satisfied.
Here $u\in\Uad$ is the control,
$y\in Y$ is the state, and $Ay=Bu$ describes an elliptic PDE\@. 
\Cref{assumptions} guarantees that the solution of \eqref{eq:LL} is unique
for each $\beta\in\Bad$, see \cref{lem:quadratic_growth_lower_level}.

The problem \eqref{eq:LL} is also called the lower-level problem.
The upper-level problem
under investigation is
\begin{equation*}
	\tag{UL}
	\label{eq:UL}
	\begin{aligned}
		\min_{\beta \in \R^n} \quad & F(\beta, y, u ) \\
		\text{s.t} \quad & \beta \in \Bad, \\
		\quad & (y,u) = \Psi(\beta),
	\end{aligned}
\end{equation*}
where $\Psi(\beta)$ describes the unique solution of \eqref{eq:LL}.
Our main motivation for studying \eqref{eq:UL} is the purpose
of identifying an unknown parameter $\beta$
from some (possibly perturbed) measurements of $\Psi(\beta)$,
see also \cref{sec:para_id}.

Together, the problems \eqref{eq:LL} and \eqref{eq:UL} constitute
the bilevel optimization problem. Necessary optimality conditions
of bilevel optimal control problems, i.e.\ hierarchical optimization problems
with two decision layers, where at least one decision maker has to solve an optimal control
problem, are derived in 
\cite{BenitaMehlitz2016,
Mehlitz2017,
Harder2021,
MehlitzWachsmuth2015:1,
Ye1995,
Ye1997}. 
Recently, solution theory for
inverse optimal control problems of partial differential equation was developed
in 
\cite{HarderWachsmuth2018:1,
HollerKunischBarnard2018}.
We also note that optimal control
problems with variational inequality constraints such as optimal control of the 
obstacle problem 
(see \cite{HarderWachsmuth2017:2})
can be viewed as a bilevel optimal control problem.
Regarding the numerical solution of the presented problem type,
there mainly exist 
(to the best of our knowledge)
methods for inverse optimal control problems
with ordinary differential equations,
see
\cite{AlbrechtUlbrich2017,
AlbrechtLeiboldUlbrich2012,
Hatz2014,
HatzSchloederBock2012}. 
The corresponding algorithms
tend to replace the lower-level problem with their optimality conditions. 
A different approach was introduced in 
\cite{DempeHarderMehlitzWachsmuth2018:1},
where the authors solved
a special class of inverse problems of partial differential equations by
exploiting the optimal-value function of the parametric optimal control problem. 
The optimal-value function $\varphi \colon\Bad\to\R$ of \eqref{eq:LL} is defined by
\begin{equation}
	\label{eq:OVF}
	\varphi(\beta)
	:=
	\inf\set[\big]{ f(\beta, y, u) \given (y,u) \in Y \times \Uad, A y = B u }
	=
	f(\beta, \Psi(\beta))
	.
\end{equation}
The idea of using the optimal-value function in bilevel optimization problems can
be traced back to \cite{Outrata1990}.
With the help of the optimal-value function,
the hierarchical problem \eqref{eq:UL}
can be transformed
into the single-level problem
\begin{equation*}
	\tag{OVR}
	\label{eq:OVR}
	\begin{aligned}
		\min_{\beta,y,u} \quad & F(\beta, y, u ) \\
		\text{s.t.} \quad & \beta \in \Bad, \\
		& f(\beta, y, u) \le \varphi(\beta), \\
		& A y - B u = 0, \\
		\quad & u \in \Uad.
	\end{aligned}
\end{equation*}
We call this optimization problem the optimal-value reformulation of \eqref{eq:UL}.
This resulting nonconvex surrogate problem does not satisfy standard constraint
qualifications such as Robinson's CQ\@.
However, in \cite[Theorem 5.12]{DempeHarderMehlitzWachsmuth2018:1}
the authors were able to prove
prove necessary optimality
conditions of Clarke-stationary type via a relaxation approach.
Furthermore, \cite[Algorithm 1]{DempeHarderMehlitzWachsmuth2018:1} introduces a
solution algorithm using a piecewise affine
approximation $\xi$ of the optimal-value function $\varphi$
with $\xi \ge \varphi$,
which leads to the relaxed optimization problem
\begin{equation*}
	\tag{OVR$(\xi)$}
	\label{eq:OVR_xi}
	\begin{aligned}
    \min_{\beta,y,u} \quad & F(\beta, y, u ) \\
    \text{s.t.} \quad & \beta \in \Bad, \\
		& f(\beta, y, u) \le \xi(\beta), \\
		& A y - B u = 0, \\
		\quad & u \in \Uad.
	\end{aligned}
\end{equation*}
If $f$ and $F$ are convex,
this problem can be split into finitely many convex subproblems for which a global
solution can be obtained.
The original problem can then be solved by iteratively
improving the approximation $\xi$ of the optimal-value function,
see \cite[Theorem 6.5]{DempeHarderMehlitzWachsmuth2018:1}.
In this paper we start with the same approach to derive a global solution
scheme. We slightly deviate in the construction of the piecewise affine
approximation by starting with a triangulation of the admissible set for the
upper-level control variable and subsequently enforce some regularity on
further divisions. 
In addition to proving convergence of the global solution
scheme in \cref{thm:conv_alg_1}, 
this will allow us to link convergence speed to the size of
the elements of the partition (see \cref{thm:convergence_speed}).
In order to solve \eqref{eq:OVR_xi}, we also consider the penalty problem
\begin{equation*}
	\tag{OVRP$(\xi)$}
	\label{eq:OVRP_xi}
	\begin{aligned}
    \min_{\beta,y,u} \quad & F(\beta, y, u ) + \gamma P( f(\beta, y, u) - \xi(\beta) ) \\
    \text{s.t.} \quad & \beta \in \Bad,  \\
		& A y - B u = 0, \\
	  \quad & u \in \Uad.
	\end{aligned}
\end{equation*}
Here, $P \colon \R \to \R$ is a penalty function and $\gamma > 0$.
Interestingly, we will see that it possible to choose the identity $P(x) = x$
as a penalty function.
This has several benefits.
On the one hand, we show in
\cref{lem:OVRP_D_Id_opt_gamma_is_multiplier_for_OVR_D}
that a finite penalty parameter can be chosen such that one obtains
the solution of \eqref{eq:OVR_xi}. 
On the other hand,
the choice of the identity
results in much simpler derivatives of the objective of
\eqref{eq:OVRP_xi}
and this enables us
to use a semismooth Newton method to solve the
subproblems efficiently, see \cref{sec:semismooth_newton}. 

Solving nonconvex problems to global optimality is an intricate issue,
and, hence, 
we expect
difficulties.
Indeed, our approach has some limitations concerning the obtained convergence
speed, see
\cref{convergence_speed_issues}.
Especially in a practical
setting convergence speed deteriorates with an increasing dimension of the upper-level
variable (curse of dimensionality).

Let us describe the structure of this paper.
In \cref{sec:preliminaries} we present the used notation as well as the main
governing assumption in addition to some preliminary theory related to optimal
control problems. We proceed by introducing a global solution algorithm
(\cref{alg:global_solution_of_IOC})
in \cref{sec:algo} and prove its convergence in \cref{thm:conv_alg_1}. 
Further we present some
convergence speed estimates in 
\cref{thm:convergence_speed} related to the size and regularity of the
elements in the partition. 
To ensure this property, we derive a simple method for refining the partition 
in arbitrary finite dimensions while keeping some regularity properties of the
elements, see \cref{lem:subdivision_by_hypercube}. 
On top of this foundation we introduce our penalty approach
(\cref{alg:global_solution_of_IOC_penalty}) in
\cref{sec:penalty}.
We show that there exists a choice of the penalty
parameter (see \cref{lem:OVRP_D_Id_opt_gamma_is_multiplier_for_OVR_D}), 
for which one can expect to find the solution to the subproblems from
\cref{alg:global_solution_of_IOC}. 
A method
for solving the penalty subproblems by means of a semismooth Newton method is
presented in \cref{sec:para_id}.
We show its superlinear convergence in
\cref{thm:fast_convergence}.
The corresponding
implementation of our algorithm for solving the 
inverse optimal control problem
and a numerical example is covered in
\cref{sec:numerical_experiments}.

\section{Preliminaries}
\label{sec:preliminaries}
\subsection{Notation}
The norm in a (real)
Banach space $X$
is denoted by
$\|\cdot\|_{X}$. 
Let $B_{X}^{\varepsilon}(x)$
denote the closed $\varepsilon$-ball 
centered at $x\in X$
with respect to $\|\cdot\|_{X}$.
Furthermore, $X\dualspace$ is the topological
dual of $X$ and 
$\langle\cdot,\cdot\rangle_{X} : X\dualspace \times X \to \R$
denotes the corresponding dual pairing.
For a set $A\subset X$
we denote by
$\operatorname{conv} A$,
$\operatorname{cone} A$, $\operatorname{cl} A$,
$\operatorname{int} A$ and $\partial A$
the convex hull, the conical hull, the closure,
interior and the boundary of $A$, respectively. 
For a Banach space $Y$,
the space of all bounded
linear operators
from $X$ to $Y$ is denoted by $L[X,Y]$
and for some operator $F\in L[X,Y]$
the adjoint is called
$F\adjoint \in L[Y\dualspace,X\dualspace]$.
For a convex set $C\subset X$ and a point $x\in C$ we denote
by
\begin{align*}
	\RR_C(x) &:= \cone(C-x),
	\\
	\NN_C(x) &:= \set{x\dualspace\in X\dualspace\given 
	\dual{x\dualspace}{y-x}_X\leq0,\;\forall y\in C}
\end{align*}
the radial cone and the normal cone to the set $C$ at the point $x \in C$,
respectively.
For $x \not\in C$, we set $\NN_C(x) := \emptyset$.

The set
$\R^{n}$ denotes the usual $n$-dimensional real
vector space, equipped with
the Euclidean norm 
$\|\cdot\|_{\R^{n}}$.
The sets $\R_{+},\R_{-}$ represent
the nonnegative and nonpositive
numbers respectively. 
For an arbitrary bounded and open set 
$\Omega\subset\R^{d}$,
the space of equivalence classes
of measurable,
$q$-integrable functions
is given by
$L^{p}(\Omega)$, $q \in [1,\infty)$.
Similarly, $L^\infty(\Omega)$
denote the space of essentially bounded
(equivalence classes of) measurable functions.
The space of functions
on the set $\Omega$
for which the $m$-th derivatives
exist in the Sobolev sense in $L^p(\Omega)$
is denoted by
$W^{m,p}(\Omega)$.
Furthermore, we use the notations
$H_0^1(\Omega) = \cl_{W^{1,2}(\Omega)}(C_c^\infty(\Omega))$
and
$H^{-1}(\Omega)
\coloneqq
H_{0}^{1}(\Omega)\dualspace$
for the Sobolev
space with first order derivatives
and homogeneous boundary conditions
and its dual space.

A mapping 
$J:X\rightarrow Y$
is called Fréchet 
differentiable at
$x\in X$ 
if there exists
an operator 
$J'(x)\in L[X,Y]$ 
such that
\begin{equation}
	\lim_{\|d\|_{X}\rightarrow 0}
    \frac{\|J(x+d)-J(x)-J'(x)d\|_{Y}}
         {\|d\|_{X}}
  = 0.
\end{equation}
In this case,
$J'(x)$
is called the Fréchet derivative
of $J$ at $x$.
If $X\ni x\mapsto J'(x)\in L[X,Y]$
is well defined and continuous 
in a neighborhood of $x$
then $J$ is said to be 
continuously Fréchet differentiable at $x$.

\subsection{Assumptions}
Throughout this work we utilize the following standing assumption.

\begin{assumption}[Standing assumption]\leavevmode
	\label{assumptions}
	\begin{enumerate}
		\item \label{asmp:function_spaces} The spaces $Y$ and $U$ are (real) Hilbert spaces.

		\item \label{asmp:S_beta}
			The set 
			$\Bad \subset \R^{n}$ 
			is a nonempty bounded polyhedron,
			i.e., a nonempty and bounded intersection of finitely many closed halfspaces.
			We assume that $\Bad$ possesses a nonempty interior.

    \item \label{asmp:U_ad} The set $\Uad\subset U$ is nonempty, closed and convex.
		\item \label{asmp:operators} The operator $A\in L[Y,Y\dualspace]$ is an isomorphism
			and $B\in L[U,Y\dualspace]$ is a linear bounded operator.
			We denote by $S := A^{-1} B \in L[U, Y]$
			the control-to-state map.

		\item \label{asmp:F_f_convex} 
			The functionals 
			$F\colon \Bad  \times Y
			\times U\rightarrow \R$ 
			and 
			$f\colon \Bad  \times Y
			\times U\rightarrow \R$
			are assumed to be bounded
			from below, convex and continuously Fréchet differentiable.
		\item \label{asmp:F_f_Lipschitz}
			The upper-level objective functional $F$ and the partial derivative
      $f'_\beta$ are assumed to be
			Lipschitz continuous on bounded sets,
			whereas $f'_u$ and $f'_y$ are Lipschitz continuous w.r.t.\ $\beta$ on bounded sets,
			i.e.,
			for every $M \ge 0$ there exists a constant $L_M \ge 0$
			such that
			\begin{equation*}
				\begin{aligned}
					\mspace{-48mu}
					\norm{
						f'_\beta(\beta_1, y_1, u_1)
						-
          f'_\beta(\beta_2, y_2, u_2)
					}_{\R^n}
					& \le
					L_M \, \parens{
						\norm{\beta_1 - \beta_2}_{\R^n}
            +        
						\norm{y_1 - y_2}_Y
						+
						\norm{u_1 - u_2}_U
					},\\
					\abs{
						F(\beta, y_1, u_1)
						-
						F(\beta, y_2, u_2)
					}
					&\le
					L_M \, \parens{
						\norm{y_1 - y_2}_Y
						+
						\norm{u_1 - u_2}_U
					},\\
					\mspace{-48mu}
	        \norm{
						f'_u(\beta_1, S(u), u)
						-
						f'_u(\beta_2, S(u), u)
					}_{U\dualspace}
					&\le
					L_M \, \norm{
						\beta_1 - \beta_2
					}_{\R^n}
 					\\
					\mspace{-48mu}
          \norm{f'_y(\beta_1,S(u),u) 
          - f'_y(\beta_2,S(u),u)}_{Y\dualspace}
           &\le L_M\norm{\beta_1-\beta_2}_{\R^n}
				\end{aligned}
			\end{equation*}
			hold for all $\beta,\beta_1,\beta_2 \in \Bad$,
      $y_1, y_2 \in B^M_Y(0)$ and $u, u_1, u_2 \in \Uad \cap B^M_U(0)$.

		\item \label{asmp:f_strongly_convex}
			The reduced lower-level objective
			$u \mapsto f(\beta, S(u), u)$
			is assumed
			to be strongly convex
			with respect to the control
			with constant $\mu>0$ independent of $\beta \in \Bad$,
			i.e.,
			\begin{equation*}
				f(\beta,S(u_{2}),u_{2}) 
				\geq
				f(\beta,S(u_{1}),u_{1})
				+
				\dual{f_{y}'(\cdot)}{S(u_{2} - u_{1})}
				+
				\dual{
				f_{u}'(\cdot)}{u_{2} - u_{1}
				}
				+
				\frac{\mu}
				{2}
				\norm{
					u_{2} - u_{1}
				}^{2}_{U}
			\end{equation*}
			holds for all
			$\beta \in \Bad$ and
			$u_{1},u_{2}\in\Uad$.
			Here,
			$f_y'(\cdot)$ and $f_u'(\cdot)$
			denote the partial derivatives of $f$
			w.r.t.\ $y$ and $u$
			at the point
			$(\beta,S(u_{1}),u_{1})$.
	\end{enumerate}
\end{assumption}

\subsection{Preliminary results}
\label{ssec:preliminary_results}
Let the optimization problem
\begin{equation*}
	\label{eq:OP}
	\tag{OP}
	\begin{aligned}
		\min_{x\in X}  \quad & J(x)\nonumber\\
		\text{s.t.} \quad & g(x)\in C
	\end{aligned}
\end{equation*}
be given, with continuously Fréchet differentiable mappings $J:X\rightarrow \R$, $g:X\rightarrow Y$
between Banach spaces $X$, $Y$ and $C\subset Y$ being nonempty, closed and convex.
A feasible point $x\in X$ of \eqref{eq:OP} satisfies the Karush-Kuhn-Tucker (KKT) conditions if
\begin{equation}
	\exists\lambda
	\in\mathcal{N}_{C}(g(x)):
	\qquad
	J'(x) + g'(x)\adjoint \lambda = 0
	.
\end{equation}
If $x$ is a local solution of \eqref{eq:OP} which satisfies Robinson's constraint qualification
\begin{equation}
	g'(x)X
	-
	\mathcal{R}_{C}(g(x))
	=
	Y
	,
\end{equation}
then the KKT conditions hold, see \cite{ZoweKurcyusz1979} and \cite[Theorem~3.9]{BonnansShapiro2000}. 
Due to \cref{assumptions}, the lower-level problem
fits into the setting of \eqref{eq:OP}.
The KKT system for the lower level
for a parameter
$\tilde\beta$ 
in a solution
$(\tilde y,\tilde u)$
then reads
\begin{equation}
	\begin{aligned}
		0&= f'_{y}\byutilde + A\adjoint\tilde p, \\
		0&= f'_{u}\byutilde - B\adjoint \tilde p+\tilde\nu, \\
		0&= A\tilde y - B\tilde u , \\
		\tilde\nu &\in \NN_{\Uad}(\tilde u),
	\end{aligned}\label{eq:OVRP_KKT}
\end{equation}
where $\tilde p\in Y$ (we identify $Y\bidualspace$ with $Y$), $\tilde \nu\in U\dualspace$ are multipliers.
Note that Robinson's CQ is satisfied due to the surjectivity of $A$.
Thus, for a minimizer of the lower-level problem
there exist multipliers such that the KKT system \eqref{eq:OVRP_KKT} is satisfied.

We can now prove that the assumption 
of strong convexity for the lower level 
implies a quadratic growth condition
in the solution.
\begin{lemma}
	\label{lem:quadratic_growth_lower_level}
	For every $\beta \in \Bad$,
	the lower-level problem \eqref{eq:LL}
	has a unique solution $(y_\beta, u_\beta)$.
	Moreover, the quadratic growth condition
	\begin{equation}
		\label{eq:quad_growth}
		f(\beta,S(u),u)
		\geq
		f(\beta, y_\beta, u_\beta)
		+
		\frac{\mu}
		{2}
		\norm{
			u - u_\beta
		}^{2}_{U}
		\qquad
		\forall u\in\Uad
	\end{equation}
	is satisfied with the parameter $\mu > 0$
	from \cref{assumptions}\ref{asmp:f_strongly_convex}.
\end{lemma}
\begin{proof}
  Existence of a solution follows from the direct method
	of calculus of variations.
	Note that the boundedness of the minimizing sequence
	follows from the strong convexity.

	Let $(y_\beta, u_\beta)$
	denote a solution of \eqref{eq:LL}.
	Utilizing the strong convexity in the solution 
	$\bybub$ 
	yields
	\begin{equation*}
				f(\beta,S(u),u) 
				\geq
				f\bybub
					+
					\dual{
					f_{u}'(\cdot)}{u - u_\beta
					}
					+
					\dual{f'_{y}(\cdot)}{S(u - u_\beta)}
					+
					\frac{\mu}
							 {2}
					\norm{
					u - u_\beta
					}^{2}_{U}
	\end{equation*}
	for all $u \in U$,
	where $f_u'(\cdot)$
	and $f_y'(\cdot)$
	denote the partial derivatives of $f$ in $\bybub$.
	By using the KKT conditions 
	with multipliers $p$, $\nu$ we obtain
	\begin{equation*}
		\begin{aligned}
			\dual{f'_{u}(\cdot)}{u - u_\beta}
			+
			\dual{f'_{y}(\cdot)}{S(u - u_\beta)} 
			&=
			\dual{f'_{u}(\cdot) + S\adjoint f'_{y}(\cdot)}{u-u_\beta}\\
			&=
			\dual{f'_{u}(\cdot) - S\adjoint A\adjoint p}{u-u_\beta}\\
			&=
			\dual{f'_{u}(\cdot) - B\adjoint p}{u-u_\beta} \\
			&=
			\dual{-\nu}{u-u_\beta}
      \ge 0
			\qquad\forall u \in \Uad.
		\end{aligned}
	\end{equation*}
	The last inequality holds since $\nu \in \NN_{\Uad}(u_\beta)$
	and $u \in \Uad$.
	Hence, one gets the quadratic growth condition
	\eqref{eq:quad_growth}.
	This also yields uniqueness of the solution.
\end{proof}
Next, we introduce the
solution operator for 
\eqref{eq:LL}.
\begin{definition}
	\label{def:Psi}
	We denote by $\Psi : \Bad \rightarrow Y \times U$
	the solution mapping of the lower-level problem
	which maps $\beta \in Q$
  to the corresponding unique solution $(y_\beta,u_\beta)$
	given in \cref{lem:quadratic_growth_lower_level}.
	We further denote by $\psi^y(\beta) \in Y$
	and $\psi^u(\beta) \in U$
	the components of $\Psi(\beta)$. As an abbreviated
  notation we introduce $y_\beta\coloneqq \psi^y(\beta)$ and
  $u_\beta\coloneqq \psi^u(\beta)$.
\end{definition}
We will now prove that the function
$\Psi$ 
is globally Lipschitz continuous.
Local Lipschitz continuity follows already by
\cite[Lemma~3.1.6]{Harder2021}. 
However, by 
\cref{assumptions}\ref{asmp:F_f_Lipschitz}
we have a stronger assumption
on the derivative of $f$.
Thus, we can adopt the arguments from
\cite[Lemma~3.1.6]{Harder2021} 
to obtain global Lipschitz continuity.

\begin{lemma}
	\label{lem:continuous_dependency_on_gamma}
	Let $X, V$ be Banach spaces,
	and let
	$C\subset X, \hat\Bad\subset V$ be
	nonempty, closed and convex sets.
	Further,
	let
	$J:X\times V\to\R$ and $\mu > 0$ be given
	such that
	for all $p \in \hat\Bad$,
	the function 
	$J(\cdot,p)$
	is strongly convex with parameter 
	$\mu$ 
	on the feasible set $C$ and Fréchet differentiable.
	Then, the solution operator 
	$\psi:\hat\Bad\to X$ 
	for the parametrized optimization problem
	\begin{equation}
		\begin{aligned}
			\min_{x}  \quad & J(x,p)\nonumber\\
			\text{s.t.} \quad & x\in C
		\end{aligned}
	\end{equation}
  exists and
	we have the estimate
  \begin{equation*}
    \norm{\psi(p_2) - \psi(p_1)}_X
    \le
    \mu^{-1}\norm{J'_x(\psi(p_2),p_1)
        - J'_x(\psi(p_2),p_2)}_{X\dualspace}
				\qquad\forall p_1, p_2 \in \hat\Bad.
  \end{equation*}
\end{lemma}
\begin{proof}
The existence of 
$\psi$ 
follows by standard arguments for convex optimization problems
with strongly convex objectives.

We now consider fixed elements 
$p_{1},p_{2}\in \hat\Bad$ 
and their corresponding unique minimizers 
$\psi(p_{i}) = x_{i}\in C$, 
$i\in \set{1,2}$.
The associated optimality conditions are
\begin{equation}
	\dual{J_{x}'(x_{i},p_{i})}
			 {\hat{x}-x_{i}}
	\geq
	0
	\quad
	\forall \hat{x}\in C.
\end{equation}
If we now add these inequalities
with the special choices 
$\hat x = x_{3-i}$,
we obtain the estimate
\begin{align*}
	0
	&\le
	\dual{J_{x}'(x_{1},p_{1}) - J_{x}'(x_{2},p_{2})}
	{x_2-x_{1}}
	\\
	&\le
	\dual{
		J_{x}'(x_{1},p_{1})
		-
		J_{x}'(x_{2},p_{1})
		+
		J_{x}'(x_{2},p_{1})
		-
		J_{x}'(x_{2},p_{2})
	}
	{x_2-x_{1}}
	\\
	&\le
	-\mu\norm{ x_2 - x_1 }_X^2
	+
	\norm{J_{x}'(x_{2},p_{1}) - J_{x}'(x_{2},p_{2})
	}_{X\dualspace} \norm{x_2 - x_1}_X.
\end{align*}
In the last step,
we have used the strong convexity of $J(\cdot, p_1)$. Dividing the last
inequality by
$\mu\norm{x_2 - x_1}_X$ yields the claim.
\end{proof}
\begin{corollary}
	\label{cor:Psi_Lipschitz}
	The function
  $\Psi$ from \cref{def:Psi} is Lipschitz continuous on $\Bad$.
	Moreover,
	there exists a constant $M_\Psi \ge 0$
	such that
	\begin{equation*}
		\norm{\beta}_{\R^n},
		\norm{\psi^y(\beta)}_{Y},
		\norm{\psi^u(\beta)}_U
		\le
		M_\Psi
		\qquad\forall \beta \in \Bad.
	\end{equation*}
\end{corollary}
\begin{proof}
	We start by proving the boundedness.
	From \cref{lem:quadratic_growth_lower_level}, we get
\begin{equation*}
  f\bybub
	+ \frac{\mu}{2} \norm{\hat u - u_\beta}_U^2
	\leq f(\beta,S(\hat u),\hat u)
	\qquad\forall \beta \in \Bad
\end{equation*}
for a fixed $\hat u\in \Uad$.
Further, $f(\cdot,S(\hat u),\hat u) : \R^n \to \R$ is continuous,
thus it is bounded on the compact set $\Bad$.
Hence, one has
\begin{equation*}
  f\bybub 
	+ \frac{\mu}{2} \norm{\hat u - u_\beta}_U^2
  \le
	C
	\qquad\forall \beta \in \Bad
\end{equation*}
for some constant $C \in \R$.
Together with the assumption that $f$ is bounded from below 
(see \cref{assumptions}\ref{asmp:F_f_convex})
we get an upper bound for
$\norm{\psi^u(\beta)}_U = \norm{u_\beta}_U$. This also allows us to bound
$\norm{\psi^y(\beta)}_Y=\norm{S(\psi^u(\beta))}_Y\le\norm{S}\norm{\psi^u(\beta)}_U$, since $S$ is a linear
bounded operator by assumption.
Since $\Bad$ is bounded, $\beta \in \Bad$ is bounded as well.
We choose $M_\Psi$ to be the largest of the previously discussed
bounds for $\norm{\beta}_{\R^n}, \norm{\psi^y(\beta)}_Y$ and $\norm{\psi^u(\beta)}_U$.

In order to prove the Lipschitzness of $\Psi$,
we want to apply \cref{lem:continuous_dependency_on_gamma} 
to the state-reduced lower-level problem,
i.e., with the setting
\begin{equation*}
	x = u, \quad
	C = \Uad, \quad
	p = \beta, \quad
	\hat\Bad = \Bad, \quad
	J(x,p) = J(u,\beta) := f(\beta, S(u), u).
\end{equation*}
\cref{assumptions} yields that the assumptions of \cref{lem:continuous_dependency_on_gamma}
are satisfied.
From the chain rule,
we get
\begin{equation*}
	J_x(u, \beta) = f'_u(\beta, S(u), u) + S\adjoint f'_y(\beta, S(u), u).
\end{equation*}
Now, \cref{lem:continuous_dependency_on_gamma} yields
\begin{equation*}
  \begin{aligned}
  \norm{\psi^u(\beta_1) - \psi^u(\beta_2)}
  &\le
  \mu^{-1}\big(\norm{f'_u(\beta_1,\psi^y(\beta_1),\psi^u(\beta_1) 
      -
      f'_u(\beta_2,\psi^y(\beta_1),\psi^u(\beta_1)}_{U\dualspace}\\
      &\qquad+
      \norm{S\adjoint}\norm{f'_y(\beta_1,\psi^y(\beta_1),\psi^u(\beta_1) 
      -
    f'_y(\beta_2,\psi^y(\beta_1),\psi^u(\beta_1)}_{Y\dualspace}\big).
  \end{aligned}
\end{equation*}
By owing to
\cref{assumptions}\ref{asmp:F_f_Lipschitz}
with $M = M_\Psi$,
this yields the desired Lipschitz continuity of $\psi^u$.
Consequently, the Lipschitz continuity of $\psi^y$
follows due to the continuity of $S$.
\end{proof}
We can use this property 
to prove the existence of solutions for 
\eqref{eq:OVR}.

\begin{theorem}
	There exists a solution for
	\eqref{eq:OVR}.
\end{theorem}
\begin{proof}
	The lower-level problem admits to a unique solution. 
	Therefore the solution operator 
	$\Psi$
	of the lower-level optimization problem 
	can be used to reduce 
	\eqref{eq:UL}
	to an optimization problem in $\R^{n}$:
	\begin{equation*}
		\begin{aligned}
			\min_{\beta}  \quad & F(\beta,\psi^y(\beta),\psi^u(\beta))
			\\
			\text{s.t.} \quad & \beta\in \Bad .
		\end{aligned}
	\end{equation*} 
	By
	\cref{assumptions}\ref{asmp:F_f_convex} 
	$F$ 
	is continuous. 
	Thus with the Lipschitz continuity of 
	$\Psi$ 
	it
	follows that
	$\beta \mapsto F(\beta,\psi^y(\beta),\psi^u(\beta))$
	is continuous.
	Moreover,
	$\Bad \subset \R^{n}$ 
	is compact by 
	\cref{assumptions}\ref{asmp:S_beta}.
	The existence of a solution follows
	from the celebrated Weierstraß theorem.
\end{proof}
We finally mention that
more general results on the existence of solutions
for bilevel optimal control problems
are given in
\cite{MehlitzWachsmuth2019:1}. 
In particular,
our result is covered by the second part of \cite[Theorem 16.3.5]{MehlitzWachsmuth2019:1}.

In order to use interpolation error estimates,
we prove regularity of the optimal-value function $\varphi$.
\begin{corollary}
	\label{col:varphi_deriv_lipschitz}
	The optimal-value function is Fréchet differentiable
	on the interior  of $\Bad$ and the derivative is Lipschitz continuous.
	In particular, we have $\varphi \in W^{2,\infty}(\Bad)$.
\end{corollary}
\begin{proof}
	The differentiability of $\varphi$ can be shown as in
	\cite[Theorem~3.2.6]{Harder2021}.
	This also yields
	the expression
	$\varphi'(\beta) = f'_{\beta}(\beta,\psi^{y}(\beta),\psi^{u}(\beta))$,
	for the derivative.
	By combining this
	with the Lipschitz continuity of
	$\Psi$ (see \cref{cor:Psi_Lipschitz}) and
	\cref{assumptions}\ref{asmp:F_f_Lipschitz},
	we get the Lipschitz continuity of $\varphi'$ on the interior of $\Bad$.
	This yields $\varphi' \in W^{1,\infty}(\Bad)$,
	see \cite[Exercise~1.x.14]{BrennerScott2008},
	and, consequently,
	$\varphi \in W^{2,\infty}(\Bad)$.
\end{proof}

\section{Algorithm}
\label{sec:algo}
In this section,
we present an algorithm to solve 
\eqref{eq:OVR}
under the given 
\cref{assumptions}.
The algorithm is similar to
\cite[Algorithm~1]{DempeHarderMehlitzWachsmuth2018:1},
with the main difference 
being the choice of the function 
$\xi$
which approximates the value function 
$\varphi$.
In that reference, the functions 
$\xi_k$ 
were defined via
\begin{equation*}
	\xi_k(x)
	:=
	\min
	\set*{ \sum_{i=1}^m\mu_i\varphi(x^i)
		\given
		0 \le \mu, \;
		\sum_{i = 1}^m \mu_i = 1, \;
		\sum_{i=1}^m\mu_i x^i=x
	},
\end{equation*}
where 
$X_k = \set{x^1, \ldots, x^m} 
\subset 
\R^n$ 
is a finite set.
The sets 
$X_k$ 
are assumed to be increasing w.r.t.\ 
$k$
and in order to achieve 
a uniform Lipschitz bound of 
$\xi_k$ 
on 
$\Bad $,
one has to require
$\Bad  \subset \interior \conv X_1$,
see 
\cite[Lemma~6.1, Example~6.1]{DempeHarderMehlitzWachsmuth2018:1}.
The reason for this extra assumption is
that it is not possible 
to a priori control the shape of the simplices
on which $\xi_k$ is affine.

We use a different method 
to obtain a bounded aspect ratio 
of all the simplices.
We choose a subdivision 
$\TT_k$ 
of 
$\Bad $ 
(recall that $Q$ is a bounded polyhedron)
into simplices.
On each simplex 
$T \in \TT_k$,
we define
$\xi_T \colon T \to \R$
as the affine interpolant of 
$\varphi$ 
in the vertices of 
$T$.
The function 
$\xi_{\TT_k}$ 
is obtained by combining 
$\xi_T$ 
for all 
$T \in \TT_k$, 
see  
\eqref{eq:def_xi_TT_k} 
below.
The advantage of this approach 
is that the approximation quality of 
$\xi_k$
can be controlled 
by the quality of the subdivision,
which is measured by the aspect ratio
\begin{equation*}
	\rho(T)
	:=
	\frac{\diam(B_T)}{\diam(T)}
	\qquad\forall T \in \TT_k
	,
\end{equation*}
where $B_T$ is the largest ball contained in $T \in \TT_k$,
see \cite[Def.~(4.2.16) and Eq.~(4.4.16)]{BrennerScott2008}.

We mention that our approach does not require continuity of $\xi_{\TT_k}$.
Therefore, we do not need 
any special assumptions on the subdivision,
in particular, we allow for hanging nodes.
In fact, it is enough to require
\begin{equation*}
	\bigcup_{T \in \TT_k} T = \Bad .
\end{equation*}
Therefore, if we have
two elements 
$T,S \in \TT_k$
with 
$T \cap S \ne \emptyset$,
the values of 
$\xi_T$ and $\xi_S$
may not agree on $T \cap S$.
For the definition of 
$\xi_{\TT_k} 
\colon 
\Bad  \to \R$, 
we choose
\begin{equation}
	\label{eq:def_xi_TT_k}
	\xi_{\TT_k}(\beta)
	:=
	\max_{T\in\TT_{k}} \xi_{T}(\beta)
	.
\end{equation}
This definition of 
$\xi_{\TT_{k}}$ 
ensures upper semicontinuity.

The main idea in \cref{alg:global_solution_of_IOC}
is to solve \eqref{eq:OVR_xi}
with $\xi = \xi_{\TT_k}$
and to successively refine
a simplex on which a solution is found.
\begin{algorithm}[t]
	\begin{enumerate}[label=\textbf{(S\arabic*)}]
		\item
			Let 
			$\TT_1$ 
			be a subdivision of 
			$\Bad $
			and select parameters
			$q, \rho \in (0,1)$
			with
			$\rho \le \min_{T \in \TT_1}\rho(T)$.
			Further, set
			$k:=1$.
		\item
			\label{alg_1_S2}
			For each
			$T\in\TT_k\setminus\TT_{k-1}$
			compute a global solution 
			$(\beta_{T},y_{T},u_{T})$ 
			of the convex optimization problem
			\begin{equation}
				\label{eq:OVR_affine_relaxed}
				\tag{OVR$(\xi,T)$}
				\begin{split}
          \min_{\beta,y,u} \quad  & F(\beta,y,u)\\
          \text{s.t.} \quad & \beta\,\in\,T,\\
          & 0\,\geq\, f(\beta,y,u)-\xi_{T}(\beta),\\
          & 0\, =\, A y - B u,\\
					\quad       & u\,\in\,\Uad.
				\end{split}
			\end{equation}
			Select 
			$\bar{T}_{k} \in \argmin_{T\in \TT_{k}} \set{F(\beta_{T},y_{T},u_{T})}$ 
			and define
			$\byuk := (\beta_{\bar{T}_{k}},y_{\bar{T}_{k}},u_{\bar{T}_{k}})$.

		\item
			\label{alg_1_S3}
			Compute 
			$\varphi(\beta_k)$. 
			If 
			$f\byuk = \varphi(\beta_k)$, then $\byuk$ is a global solution of 
			\eqref{eq:OVR} 
			(and, thus, of 
			\eqref{eq:UL}) 
			and the algorithm terminates. 
			Otherwise, we construct 
			$\TT_{k+1}$ 
			from 
			$\TT_k$ 
			by a refinement of 
			$\bar T_k$ 
			such that 
			$\vol(T) \leq q\cdot \vol(\bar{T}_{k})$
			and $\rho(T) \ge \rho$
			for all 
			$T\in\TT_{k+1}\setminus\TT_{k}$ .
			Set $k := k+1$ and go to \ref{alg_1_S2}.
	\end{enumerate}
	\caption{Computation of global solutions to 
	\eqref{eq:UL}\label{alg:global_solution_of_IOC}}
\end{algorithm}
In order for \cref{alg:global_solution_of_IOC} to be well-defined,
we need to guarantee the existence of global minimizers of
\eqref{eq:OVR_affine_relaxed}.
This can be shown by the direct method of calculus of variations.
The boundedness of $\beta$ follows from $\beta \in T$
and the boundedness of $(y,u)$ follows from $f(\beta, y, u) \le \xi_T(\beta)$,
cf.\ \cref{assumptions}\ref{asmp:f_strongly_convex}.

Under very mild assumptions
we can show the convergence towards global minimizers.
\begin{theorem}
	\label{thm:conv_alg_1}
	\Cref{alg:global_solution_of_IOC} 
	either stops at a global solution of 
	\eqref{eq:OVR}
	or the computed sequence 
	$\byuk$ 
	contains a subsequence 
	converging strongly in $\R^n \times Y \times U$ to a global solution of 
	\eqref{eq:OVR}. 
	If \eqref{eq:OVR} has a unique global solution
	$\byubar$,
	then the entire sequence $\byuk$ converges strongly to $\byubar$.
\end{theorem}
\begin{proof}
	The value function $\varphi$ 
	is convex and therefore 
	$\xi_{\TT_k}(\beta)
		\geq 
		\varphi(\beta)$. 
	Thus, the feasible set of
	\hyperref[eq:OVR_xi]{(OVR$(\xi_{\TT_k})$)} 
	contains the feasible set of
	\eqref{eq:OVR}.
	If the solution $\byuk$
	of
  \hyperref[eq:OVR_xi]{(OVR$(\xi_{\TT_k})$)}
	is feasible for 
	\eqref{eq:OVR},
	it is globally optimal for 
	\eqref{eq:OVR}. 
	Hence, the stopping criteria of the algorithm
	ensures that $\byuk$ is globally optimal for \eqref{eq:OVR}.
	It remains to discuss the case where \cref{alg:global_solution_of_IOC}
	does not terminate.
	We denote by
	$\byubar$ 
	a global solution of
	\eqref{eq:OVR}.
	Then
\begin{equation}
F\byuk \leq F\byubar\label{eq:UL_objective_value_estimate_for_proof_conv_alg_1}
\end{equation}
by the same argument. 
The feasible set
$\Bad $ 
is compact by 
\cref{assumptions}\ref{asmp:S_beta}.
This implies the existence of $N \in \R$ with $\varphi(\beta) \le N$
for all $\beta \in \Bad $.
Therefore, the estimate
\begin{equation*}
	N
	\ge
	\xi_{\TT_k}(\beta_k)
	\ge
	f\byuk
	\ge
	f(\beta_k, y_{\beta_k}, u_{\beta_k})
	+
	\frac\mu 2 \norm{ u_{\beta_k} - u_k }_U^2
\end{equation*}
(where we used \eqref{eq:quad_growth} in the last step)
together with the boundedness of $u_{\beta_k}$
shows the boundedness of $u_k$ in $U$.
The boundedness of $y_k$ in $Y$
follows from the properties of the linear operators $A$ and $B$.
Therefore the sequence 
$\byuk$ 
is bounded
by a constant $M \ge 0$
and contains 
a weakly convergent subsequence (without relabeling)
$\byuk
	\rightharpoonup 
	\byuhat$
	in $\R^n \times Y \times U$.
	In particular,
one has strong the convergence
$\beta_k
\rightarrow \bar{\beta}$,
since $\R^n$ is finite dimensional.

In order to estimate
the distance between $\varphi$ and its interpolant $\xi_{\TT_k}$,
we use the interpolation error estimate
\cite[Theorem 4.4.20]{BrennerScott2008}
(the required condition \cite[(4.4.16)]{BrennerScott2008}
is satisfied due to \ref{alg_1_S3} in \cref{alg:global_solution_of_IOC}).
We apply this result (for polynomial degree one with $m = 2$, $s=0$, $p=\infty$)
on each simplex $T \in \TT_k$
and obtain
\begin{equation}
	\norm{\xi_T - \varphi}_{L^\infty(T)}
	\leq
	C_{\rho} \, \diam (T)^{2} \norm{\varphi}_{W^{2,\infty}(T)}
	\qquad
	\forall T \in \TT_k,
	\label{eq:interpolation_error}
\end{equation}
where $C_\rho>0$ is a constant that depends on the regularization parameter $\rho$.
\Cref{col:varphi_deriv_lipschitz}
provides the upper bound
$\norm{\varphi}_{W^{2,\infty}(T)} \leq \norm{\varphi}_{W^{2,\infty}(\Bad )}
=: C_{\varphi}$.
We want to apply \eqref{eq:interpolation_error} for $\bar T_k\in\TT_k$,
where $\bar T_k$ is chosen as in the algorithm,
and also intend to show $\diam(\bar T_k)\to0$.
We will use the relation between diameter
and volume given by the aspect ratio 
of the simplices and argue by contradiction. 
We assume that 
$v\coloneqq\limsup_{k\rightarrow\infty}\vol(\bar{T}_{k})>0$.
Thus the set 
$\bar \TT_0 := \set{ \bar T_k \given k \in \N, \vol(\bar T_k) \ge v}$
is infinite. Now there has to be at least one simplex  
$T_{0} \in \TT_1$
that contains infinitely many simplices 
from $\bar \TT_0$,
i.e.,
the set
$\bar{\TT}_1 
\coloneqq 
\set{T \in \bar\TT_0 \given T \subsetneq T_{0}}$
is infinite.
These simplices are refined 
at least once and thus we have 
$\vol(T) \leq q \vol (T_0)$ for all $T \in \bar\TT_1$.
Again, one simplex in $\bar\TT_1$ has to contain
infinitely many of the simplices from $\bar\TT_1$
and we can
repeat the above argument.
This leads to a contradiction as
the volume of the simplices
is bounded from above by $q^{-l} \vol (T_0)$
and this contradicts the lower bound $v > 0$.
Hence, we have shown $\vol(\bar T_k) \to 0$.
Using the bound on the aspect ratio,
this implies $\diam(\bar T_k) \to 0$.
Indeed,
\begin{equation*}
	\diam(\bar{T}_{k}) 
	\le
	\frac{\diam(B_{\bar T_{k}})}{\rho}
	= \frac{2}{\rho}\parens*{\frac{\Gamma(\frac{n}{2}+1)\vol(B_{\bar T_{k}})}{\pi^{\frac{n}{2}}}}^{1/n}
	\leq
	\frac{2}{\rho}\parens*{\frac{\Gamma(\frac{n}{2}+1)\vol(\bar T_{k})}{\pi^{\frac{n}{2}}}}^{1/n}
	\to 0
	.
\end{equation*}
Now we are in position to apply \eqref{eq:interpolation_error} on $\bar T_k$.
This yields
\begin{equation}
	\begin{aligned}
		\varphi(\hat{\beta}) 
		\leq 
		f\byuhat 	&\leq \liminf_{k\rightarrow\infty}
											f\byuk
							\leq \limsup_{k\rightarrow\infty} 
											f\byuk
							\leq \limsup_{k\rightarrow \infty} 
											\xi_{\bar{T}_{k}}(\beta_{k})\\
							&\leq \limsup_{k\rightarrow \infty} 
											\left(
											\varphi(\beta_{k}) 
											+ C_\rho C_\varphi \,
												\diam(\bar{T}_{k})^{2}
											\right)
							= \varphi(\hat{\beta}).
	\end{aligned}
	\label{eq:limit_varphi_constraint}
\end{equation}
Note that we have used
the sequential weak lower semicontinuity of $f$
which follows from convexity and continuity
in \cref{assumptions}\ref{asmp:F_f_convex}.
Thus, \eqref{eq:limit_varphi_constraint} yields
feasibility of 
$\byuhat$
for
\eqref{eq:OVR}.
Similarly, $F$ is sequentially weakly lower semicontinuous.
Therefore, we can pass to the limit $k \to \infty$
in \eqref{eq:UL_objective_value_estimate_for_proof_conv_alg_1}
and obtain
\begin{equation} 
	F\byuhat
	\leq 
	\liminf_{k\rightarrow\infty}
		F\byuk
	\leq 
		F\byubar
		.
\end{equation}
This shows
that
$\byuhat$ 
is a global solution for 
\eqref{eq:OVR}.

Next, we prove the strong convergence of $y_k$ and $u_k$.
Strong convergence of the control 
$u_k$
can be obtained by exploiting
the quadratic growth condition from \cref{lem:quadratic_growth_lower_level}:
Note that 
$y_{k} = S(u_{k})$
by feasibility of 
$\byuk$ 
for
\hyperref[eq:OVR_affine_relaxed]{\textup{(OVR($\xi,\bar T_k$))}}.
Thus,
\cref{lem:quadratic_growth_lower_level}
and the
Lipschitz continuity of $f'_\beta(\hat \beta,\cdot,\cdot)$
from
\cref{assumptions}\ref{asmp:F_f_Lipschitz}
yield
\begin{align}
	\nonumber
	f\byuk 
	&\geq f(\hat{\beta},y_{k},u_{k})+\dual{f'_{\beta}(\hat{\beta},y_{k},u_{k})}{\beta_{k}-\hat{\beta}}
	\\
	\nonumber
	&\ge f(\hat{\beta},y_{k},u_{k}) 
	-\norm{f'_{\beta}(\hat{\beta},y_{k},u_{k})}_{\R^n}\norm{\beta_{k}-\hat\beta}_{\R^n}\\
	\nonumber
	&\ge f(\hat{\beta},y_{k},u_{k}) 
	\\\nonumber&\qquad{}-	\left(
			\norm{f'_{\beta}\byuhat}_{\R^n}
				+ L_M \norm{ y_{k}-\hat{y} }_Y
				+ L_M \norm{ u_{k}-\hat{u} }_U
		  \right)
			\norm{\beta_{k}-\hat{\beta}}_{\R^n}\\
			\label{eq:from_beta_k_to_hat_beta}
	&\geq 
	f\byuhat 
	+ \frac{\mu}{2}\norm{u_{k}-\hat{u}}^{2}_{U}
	-C\norm{\beta_{k}-\hat{\beta}}_{\R^n}
	.
\end{align}
Since
\eqref{eq:limit_varphi_constraint}
implies
$f\byuk \to f\byuhat$
and since $\beta_k \to \hat\beta$,
this inequality
yields
the strong convergence $u_k \to \hat u$ in $U$.
The continuity 
of the solution operator $S$ now
implies strong convergence of the states.

If the solution to 
\eqref{eq:OVR} 
is unique,
the convergence of the entire sequence
follows from a usual subsequence-subsequence argument.
\end{proof}

An important ingredient of \cref{alg:global_solution_of_IOC}
is the refinement of the simplices
in \ref{alg_1_S3}
such that the properties
involving the constants $q$ and $\rho$
are obtained.
In the two-dimensional case 
$\Bad\subset\R^{2}$ 
this can be done 
by splitting the triangle 
$\bar T_k$
into 
$4$
similar triangles by using the midpoints
of the edges. 
However, already in three dimensions 
this is not straightforward 
since
a general tetrahedron
cannot be divided into similar tetrahedrons.
In particular, a regular tetrahedron cannot
be split into smaller regular tetrahedra.
One, however, can use hypercubes
to construct a method 
of refinement 
that maintains a bounded aspect ratio.

\begin{lemma}
	\label{lem:subdivision_by_hypercube}
	For every (finite) subdivision $\TT_1$,
	there exist constants
	$q, \rho \in (0,1)$
	such that
	the refinement in
	\ref{alg_1_S3}
	of \cref{alg:global_solution_of_IOC}
	is always possible.
\end{lemma}

\begin{proof}
Let $S_n$ denote the permutations of $\set{1,2,\ldots, n}$.
We consider the hypercube 
$[0,1]^{n}$ 
and a permutation 
$\pi \in S_n$.
Then
$T_{\pi} 
\coloneqq 
\set{x\in\R^{n} \given 0\leq x_{\pi(1)} 
								 \leq \dots 
								 \leq x_{\pi(n)} 
								 \leq 1}$ 
describes a simplex. 
For each point $x$ in the hypercube 
there exists at least one permutation $\pi$
for which the definition of 
$T_{\pi}$ 
is consistent with the 
``$\leq$''-ordering 
of the components of $x$, i.e., $x \in T_\pi$.
Therefore 
$\bigcup_{\pi\in S_{n}}T_{\pi}
=
[0,1]^{n}$.
If we consider a point 
$x\in[0,1]^{n}$ with $x_{i} 
\neq 
x_{j}$ 
for all
$i\neq j$, 
then there exists only one permutation
$\pi$ 
such that
$x\in T_{\pi}$
since
the components of $x$
have a uniquely determined order. 
Furthermore, those points 
are dense in 
$[0,1]^n$ 
and this implies 
that two simplices constructed 
with two different permutations 
cannot have a 
$n$-dimensional intersection.
Moreover, different simplices $T_\pi$
can be matched by a permutation of the coordinates
and this implies
that the volume
of each $T_\pi$ is equal to $1/n!$
and the aspect ratio $\rho(T_\pi)$
is independent of $\pi$.

The hypercube can be split into 
$2^{n}$
smaller cubes.
By dividing these smaller cubes
again into simplices,
we arrive at
\begin{equation}
	T_{\pi}^{t} 
	\coloneqq 
	\set{x\in\R^{n} \given 0	\leq x_{\pi(1)}-t_{\pi(1)}
										\leq \dots 
										\leq x_{\pi(n)}-t_{\pi(n)}
										\leq 0.5
	}
	,
\end{equation} 
where we consider all possible 
$t \in \set{0, 0.5}^n$ and $\pi \in S_n$.
We observe that these simplices
are the translated and scaled versions of $T_\pi$.
In particular, we have
$T_\pi^t = \frac12 T_\pi + t$
and
this implies
$\vol(T_\pi^t) = 2^{-n} \vol(T_\pi) = 2^{-n} / n!$
and
$\rho(T_\pi^t) = \rho(T_\pi)$.

We argue that for all $\pi \in S_n$ and $t \in \set{0, 0.5}^n$,
there exists $\hat\pi \in S_n$ with $T_\pi^t \subset T_{\hat\pi}$.
Indeed, for $x \in T_\pi^t$,
the coordinates $x_i$ with $t_i = 0$
are smaller (or equal) than the coordinates
$x_j$ with $t_j = 0.5$.
Further, we have $x_{\pi(i_1)} \le x_{\pi(i_2)}$
if $t_{\pi(i_1)} = t_{\pi(i_2)}$ 
and $i_1 \le i_2$.
Thus, we can construct $\hat\pi$
by first taking the indices $\pi(i)$ with $t_{\pi(i)} = 0$
and afterwards the indices $\pi(j)$ with $t_{\pi(j)} = 0.5$.
Due to $\vol(T_\pi^t)=2^{-n}\vol(T_\pi)$ this implies
that every $T_\pi$ can be divided into $2^n$ smaller simplices
$T^{t^{(i)}}_{\pi^{(i)}}$ with $i = 1,\ldots, 2^n$.
Again, these smaller simplices have the same aspect ratio as $T_\pi$.

Repeating this subdivision
proves the assertion in the case
that $\TT_1 \subset \set{ T_\pi \given \pi \in S_n}$
with the constants
$q = 2^{-n}$, $\rho = \rho(T_\pi)$ for some fixed $\pi \in S_n$.

In the general case,
we map each simplex $T \in \TT_1$
to $T_\pi$ for some fixed $\pi \in S_n$
by an (invertible) affine transformation
$a:T\to T_\pi$.
The first part of the proof showed
that $T_\pi$ can be divided repeatedly
into smaller simplices.
In each subdivision step,
the volume is scaled down by $2^{-n}$
whereas the aspect ratio is constant.
By applying the inverse transformation $a^{-1}$,
we get a subdivision of $T$.
The ratio of volumes is invariant w.r.t.\ the affine transformation $a^{-1}$,
thus we can take $q = 2^{-n}$.
It remains to study the effect of the affine transformation $a^{-1}$
on the aspect ratio.
Every simplex that is the result of repeated refinement
of $T$ has the form $a^{-1}(\hat T)$, where $\hat T\subset T_\pi$
is a simplex which has the same aspect ratio as $T_\pi$.
We denote the largest balls in $\hat T$ and $a^{-1}(\hat T)$
by $B_{\hat T}$ and $B_{a^{-1}(\hat T)}$.
The ellipsoid
$a^{-1}(B_{\hat T})$
is contained in $a^{-1}(\hat T)$
and it contains a ball of diameter
$\theta_{\min} \diam(B_{\hat T})$,
where $\theta_{\min}$ is the smallest spectral value of the matrix $a'(0)^{-1}$.
Thus,
\begin{equation*}
	\diam(B_{a^{-1}(\hat T)})
	\ge
	\theta_{\min}
	\diam(B_{\hat T}).
\end{equation*}
Similarly, we get
\begin{equation*}
	\diam( a^{-1}(\hat T) )
	\le
	\theta_{\max}
	\diam( \hat T ),
\end{equation*}
where $\theta_{\max}$ is the largest singular value of the matrix $a'(0)^{-1}$.
This yields the estimate
\begin{equation*}
	\rho( a^{-1}(\hat T) )
	=
	\frac{\diam( B_{a^{-1}(\hat T)})}{\diam( a^{-1}(\hat T) )}
	\ge
	\frac{\theta_{\min}}{\theta_{\max}}
	\cdot
	\frac{\diam(B_{\hat T})}{\diam( \hat T )}
	=
	\frac{\theta_{\min}}{\theta_{\max}}
	\rho(\hat T)
	=
	\frac{\theta_{\min}}{\theta_{\max}}
	\rho(T_\pi)
	=:
	\rho_T
	>
	0.
\end{equation*}
Thus, the aspect ratio of every simplex that is the result
of repeated refinement of $T$ can be bounded from below by $\rho_T$.
Since $\TT_1$ is finite,
we can choose
$\rho = \min\set{\rho_T \given T \in \TT_1} > 0$.
\end{proof}

\begin{remark}
	The refinement technique of \cref{lem:subdivision_by_hypercube}
	always generates hanging nodes.
	The presented method is consistent 
	with splitting a triangle into 
	$4$ 
	similar parts using 
	the midpoints of the edges.
	In higher dimensions 
	there might exist more advanced methods. 
	Since \cref{alg:global_solution_of_IOC}
	only requires a bound on the aspect ratio,
	we can use the simple strategy from \cref{lem:subdivision_by_hypercube}.
\end{remark}

After we have proven 
the convergence
of \cref{alg:global_solution_of_IOC},
we want to get an estimate 
on the convergence speed. 
We establish a preliminary result 
on the error in the upper-level objective
induced by the approximation $\xi_T$ of $\varphi$.

\begin{lemma}
\label{lem:error_on_one_element}
Let $\TT$ be a subdivision of $\Bad$.
For $T\in\TT$ and any feasible point $\byu$ of 
\eqref{eq:OVR_affine_relaxed}
we have
\begin{equation}
	\abs{F\byu - F(\beta,y_{\beta},u_{\beta})} 
	\leq L_M(1+\|S\|)
	\sqrt{\frac{2C_{\rho}C_{\varphi}}
		{\mu}
	}
	\diam(T)
	,
\end{equation}
where $(y_\beta, u_\beta)$ is the solution of the lower-level problem associated with the parameter $\beta$,
see \cref{def:Psi}.
Here, $C_\rho$ is as in \eqref{eq:interpolation_error}
and $C_\varphi := \norm{\varphi}_{W^{2,\infty}(\Bad)}$.
The constant $M$ does not depend directly on $T$ but only on $\rho(T)$.
\end{lemma}

\begin{proof}
	We use the quadratic growth condition from \cref{lem:quadratic_growth_lower_level}
	to obtain
\begin{equation*}
	\xi_T(\beta) \geq f\byu
										\geq f(\beta,y_{\beta},u_{\beta}) 
													+ \frac{\mu}
																 {2}
														\|u-u_{\beta}\|_{U}^{2}
						 = \varphi(\beta) 
										+ \frac{\mu}
													 {2}
											\|u-u_{\beta}\|_{U}^{2}
											.
\end{equation*}
Next, we
apply the interpolation estimate 
\eqref{eq:interpolation_error} 
to get
\begin{equation}
	\norm{u-u_{\beta}}^{2}_{U} 
	\leq 
	\frac{2C_\rho C_\varphi \diam(T)^{2}}
			 {\mu}
			.
			\label{eq:norm_bound_control}
\end{equation}
In order to apply the Lipschitz assumption from \cref{assumptions},
we define
$M := M_\Psi + \max\{1,\norm{S}\} \sqrt{2 C_\rho C_\varphi / \mu} \diam(Q)$,
where $M_\Psi$ is given in \cref{cor:Psi_Lipschitz}.
Due to \eqref{eq:norm_bound_control},
all quantities are bounded by $M$.
Thus,
	\begin{align*}
		\abs{F\byu - F(\beta,y_{\beta},u_{\beta})} 
			&\leq L_M (\norm{S u - S u_\beta}_{Y} + \norm{u - u_\beta}_{U})\\
			&\leq L_M (1+\|S\|)\norm{u - u_\beta}_{U}\\
			&\leq L_M (1+\|S\|)
						\sqrt{\frac{2C_{\rho}C_{\varphi}}
											 {\mu}
						}
						\diam(T)
						.
						\qedhere
	\end{align*}
\end{proof}
\begin{theorem}
\label{thm:convergence_speed}
Let $\TT$ be a subdivision of $\Bad$ and suppose that the upper-level objective functional
satisfies a quadratic growth condition 
for a solution $\byubar$ 
of 
\eqref{eq:OVR}
in the sense that
\begin{equation}
	\label{eq:quad_growth_UL}
	F(\beta, y_\beta, u_\beta)
	\geq 
	F\byubar
		+ G
		\norm{\beta - \bar\beta}_{\R^n}^{2}
			\qquad\forall \beta \in \Bad
\end{equation}
holds for some constant $G > 0$.
Let
$T\in\TT$ 
be an element
satisfying the condition
\begin{equation}
	\label{eq:small_T}
	\diam(T) 
	<
		\frac{G}{L_{M}
					(1+\|S\|)
					\sqrt{\frac{2C_{\rho}C_{\varphi}}
										 {\mu}}
							  }
	\dist(T, \bar\beta)^{2} 
	.
\end{equation}
Then,
for any feasible point
$\byu$ of the relaxed problem \eqref{eq:OVR_affine_relaxed}
we have
\begin{equation*}
	F\byu > F\byubar.
\end{equation*}
The constants appearing in \eqref{eq:small_T}
have the same meaning as in \cref{lem:error_on_one_element}.
\end{theorem}
\begin{proof}
	Let $T \in \TT$ satisfy \eqref{eq:small_T}
	and let
$\byu$ 
be feasible to 
\eqref{eq:OVR_affine_relaxed}.
By using the quadratic growth condition \eqref{eq:quad_growth_UL} and \cref{lem:error_on_one_element}
we obtain
\begin{equation}
	\label{eq:alg_conv_speed_calculation}
	\begin{aligned}
		F\byu - F\byubar 
		&=
		F\bybub 
		- F\byubar
		+
		F\byu 
		- F\bybub
		\\
		&\geq G
		\norm{\beta - \bar\beta}_{\R^n}^{2}
		-L_M (1+\norm{S})
		\sqrt{\frac{2C_{\rho}C_{\varphi}}
			{\mu}
		}
		\diam(T)
		\\
		&>
		G \norm{\beta - \bar\beta}_{\R^n}^{2}
		- G \dist(T, \bar\beta)^2
		\ge
		0.
	\end{aligned}
\end{equation}
This shows the claim.
\end{proof}

\begin{remark}
  \label{convergence_speed_issues}
	We give some interpretation of
	\cref{thm:convergence_speed}.
	Let $\byubar$ be
  a solution to \eqref{eq:OVR} satisfying the growth condition
  \eqref{eq:quad_growth_UL}.
	Let $T \in \TT$ satisfy \eqref{eq:small_T}
	and let $\byu$ be a feasible point of
  \hyperref[eq:OVR_affine_relaxed]{\textup{(OVR($\xi, T$))}}.
	Further, let $\bar T \in \TT$ be a simplex with $\bar\beta \in \bar T$.
	Then, a solution $(\beta_{\bar T}, y_{\bar T}, u_{\bar T})$ of
	\hyperref[eq:OVR_affine_relaxed]{\textup{(OVR($\xi, \bar T$))}}
	satisfies
	\begin{equation*}
		F\byu > F\byubar
		\ge
		F(\beta_{\bar T}, y_{\bar T}, u_{\bar T})
		.
	\end{equation*}
  Hence,
	\cref{alg:global_solution_of_IOC}
	will never refine the simplex $T$
	and, consequently,
	this simplex will be ignored in the
	subsequent iterations
	of the algorithm.

	\cref{thm:convergence_speed}
	also has a quantitative implication.
	We consider a subdivision of $Q$ into simplices
	of diameter $h$.
	According to \eqref{eq:small_T},
	the minimizer $\bar \beta$
	cannot occur in simplices $T$
	with $h < C \dist(T, \bar\beta)^2$,
	with some constant $C > 0$.
	That is, we only have to consider simplices
	with $\dist(T, \bar\beta) \le \sqrt{h/C}$.
	The number of simplices satisfying this condition
	is roughly of the order $h^{n/2 - n} = h^{-n/2}$.

  If we are able to improve \eqref{eq:small_T}
	to $\diam(T) < C \dist(T, \bar \beta)^\alpha$ for some $\alpha \in [1,2)$,
	see the discussion below,
	this number of simplices improves to $h^{-n (1-1/\alpha)}$.
	In particular, in the case $\alpha = 1$,
	we expect a constant number of simplices.
\end{remark}
\begin{remark}
	\label{rem:better_small_T}
There are two possibilities to improve condition
\eqref{eq:small_T}. First, if one has a stronger growth condition for the upper-level objective functional,
i.e.,
\begin{equation}
	\label{eq:better_then_quad_growth_UL}
	F(\beta, y_\beta, u_\beta)
	\geq 
	F\byubar
		+ G
		\norm{\beta - \bar\beta}_{\R^n}^{\alpha}
			\qquad\forall \beta \in \Bad
\end{equation}
for some $\alpha \in [1, 2)$,
then we can use
$\dist(T,\bar\beta)^\alpha$
instead of
$\dist(T,\bar\beta)^2$
in \eqref{eq:small_T},
cf.\ \eqref{eq:alg_conv_speed_calculation}.
In particular, $\alpha = 1$
might be possible if $\bar\beta$ is located on the boundary of $\Bad$
or if the reduced objective is non-smooth at $\bar\beta$.

Second, we
can improve \cref{thm:convergence_speed} if
$F'\byubar = 0$.
For simplicity, we discuss the case that
$F$ is quadratic, i.e.,
\begin{equation}
\begin{aligned}
	F\byu
	&=
	F\bybub + F'\bybub (\byu - \bybub)\\
	&\qquad+ \frac12 F''\bybub [\byu - \bybub]^2.
\end{aligned}
	\label{eq:F_quadratic}
\end{equation}
In particular, the second derivative is constant.
Together with the Lipschitz continuity of $F'$ and 
$\Psi$ (see \cref{cor:Psi_Lipschitz}),
we readily obtain
\begin{equation*}
	\norm{ F'( \beta, y_\beta, u_\beta) }_{\R^n \times Y\dualspace \times U\dualspace}
	=
	\norm{ F'( \beta, y_\beta, u_\beta) -  F'\byubar}_{\R^n \times Y\dualspace \times U\dualspace}
	\le
	C \norm{\beta - \bar \beta}_{\R^n}
	.
\end{equation*}
Using this estimate and \eqref{eq:norm_bound_control} in \eqref{eq:F_quadratic},
we find
\begin{align*}
	\abs{F\byu - F\bybub}
	&
	\le
	C \norm{\beta - \bar \beta}_{\R^n}\diam(T)
	+ C \diam(T)^2
	\\&
	\le
	C \dist(T,\bar\beta)\diam(T) + C \diam(T)^2 .
\end{align*}
By using this estimate in \eqref{eq:alg_conv_speed_calculation},
we see that \eqref{eq:small_T} can be replaced by
$\diam(T) < c \dist(T, \bar\beta)$
for some $c > 0$.
Note that $F'\byubar = 0$ is highly restrictive. 
However, the positive influence on the convergence speed
can already be expected if the first derivative of 
$F$ is close to zero in the solution. 
The approach can be applied to
non-quadratic objective functionals 
$F$ by replacing \eqref{eq:F_quadratic} by a Taylor expansion
and requiring that $\norm{F''}$ is bounded on bounded subsets.
\end{remark}

\cref{alg:global_solution_of_IOC}
can still be sped up substantially without additional restrictions.
In \ref{alg_1_S3},
we have to evaluate $\varphi(\beta_k)$,
and for this purpose we calculate the lower-level solutions
$(y_{\beta_k}, u_{\beta_k})$.
Therefore
$(\beta_k, y_{\beta_k}, u_{\beta_k})$
is a feasible point of \eqref{eq:OVR}
and, thus,
$F(\beta_k, y_{\beta_k}, u_{\beta_k})$
is an upper bound for the
minimal objective value of \eqref{eq:OVR}.
On the other hand,
the computed values $F(\beta_T, y_T, u_T)$
for $T \in \TT$
are lower bounds for the possible objective value of \eqref{eq:OVR}
restricted to $T$.
Hence, all elements $T \in \TT$
with
$F(\beta_T, y_T, u_T) > F(\beta_k, y_{\beta_k}, u_{\beta_k})$
cannot contain a solution of \eqref{eq:OVR}
and can be ignored in later iterations.
Furthermore, the simplices can be sorted by $F(\beta_T, y_T, u_T)$
and multiple simplices may be refined in each iteration.
This results in a larger number of auxiliary problems
which have to be solved in the next iteration
(recall that \eqref{eq:OVR_affine_relaxed} has to be solved on refined elements only).
These problems are independent
of each other and can be solved in parallel. 

Finally,
we demonstrate that in most cases,
the value-function constraint
in
\eqref{eq:OVR_affine_relaxed}
will be satisfied
with equality.
To study the issue we introduce the problem
\begin{equation}
	 \label{eq:UP_beta}
	\begin{aligned}
		\min_{\beta,y,u} \quad & F(\beta, y, u ) \\
		\text{s.t.} \quad & A y - B u = 0, \\
		\quad & \beta\in\Bad,\quad u\in\Uad.
	\end{aligned}
\end{equation}
This problem is a relaxation of \eqref{eq:OVR},
since we neglected the optimality of $(y,u)$ for the lower level.
We expect that this problem has a smaller optimal value than \eqref{eq:OVR}.
\begin{lemma}
	\label{lem:f_equals_xi}
	Suppose that the infimal value of \eqref{eq:UP_beta}
	is smaller than the infimal value of \eqref{eq:OVR}. 
  Let $\byuk$ be defined as in 
  \cref{alg:global_solution_of_IOC}\ref{alg_1_S2}.
	Then,
	the constraint
	$f\byuk \leq \xi_{\bar T_k}(\beta_k)$
	is satisfied with equality
	for all $k$ large enough and for which $\xi_{\TT_k}$
	is continuous at $\beta_k$.
\end{lemma}        
\begin{proof}
Let 
$\byutilde$ 
be a global solution for 
\eqref{eq:UP_beta}. 
Note that global solutions
$\byubar$ 
to \eqref{eq:OVR}
are not globally optimal for
\eqref{eq:UP_beta}. 
The construction of the sequence 
$\byuk$ according to
\cref{alg:global_solution_of_IOC}
yields a monotonically increasing sequence $F\byuk$.
By 
\cref{thm:conv_alg_1} one gets
$F\byuk \rightarrow F\byuhat = F\byubar$.
Due to $F\byutilde < F\byubar$,
we have
$F\byutilde < F\byuk$ for sufficiently large $k$.

We argue by contradiction and assume that
$f\byuk < \xi_{\bar T_k}(\beta_k)$ for some large $k$ for which
$\xi_{\TT_k}$ is continuous at $\beta_k$.
We consider a convex combination
$(1-s) \byuk + s \byutilde$,
$s \in (0,1)$,
and check that it is a feasible point of
\hyperref[eq:OVR_xi]{\textup{(OVR($\xi_{\TT_k}$))}}
for $s$ small enough.
The constraint $Ay = Bu$ 
is linear and the admissible sets
$\Bad$ and $\Uad$ are convex. 
Moreover,
since $f$ is continuous (see \cref{assumptions}) and
since $\xi_{\TT_k}$
is continuous by assumption, we have
\begin{equation*}
  f((1-s) \byuk + s \byutilde)
  <
  \xi_{\TT_k}( (1-s) \beta_k + s \tilde\beta)
  \qquad
  \forall s\in(0,\varepsilon]
  .
\end{equation*}
for some $\varepsilon > 0$.
Now the convexity of the upper-level objective functional $F$ 
(see \cref{assumptions}\ref{asmp:F_f_convex})
implies
\begin{equation*}
	F((1-s) \byuk + s \byutilde)
	\le
	(1 - s) F\byuk
	+
	s F\byutilde
	<
	F\byuk
\end{equation*}
for all $s \in (0,\varepsilon]$.
This contradicts the optimality of $\byuk$ from 
\cref{alg:global_solution_of_IOC}\ref{alg_1_S2}.
\end{proof}
Note that the piecewise linear function $\xi_{\TT_k}$ is continuous
if the triangulation $\TT_k$ does not possess hanging nodes.
Otherwise, it might be discontinuous at all facets containing hanging nodes.

\section{Penalty approach}
\label{sec:penalty}
The subproblems 
\eqref{eq:OVR_affine_relaxed}  
presented in 
\cref{alg:global_solution_of_IOC}
are already subject to convex constraints,
however, the nonlinear inequality constraint 
$f\byu \leq \xi(\beta)$ 
still may introduce difficulties 
when implementing the solution algorithm.
In particular,
this constraint is of a rather unusual form
in an optimal control context, see \cref{sec:para_id}.
Using a penalty method
for this complicated constraint
the treatment of the subproblems 
\eqref{eq:OVR_affine_relaxed}
can be simplified
since this inequality constraint is incorporated
into the objective functional.
Any additional error that is
introduced by the penalty approach
has to be compared to the 
error induced by the relaxation of the problem
with the affine interpolation of the
optimal-value function. 

By replacing the subproblems in \cref{alg:global_solution_of_IOC}
with a penalty approach,
we arrive at \cref{alg:global_solution_of_IOC_penalty}
for which we now provide
some further comments.
In a classical penalty method
the penalty parameter
depends only on the iteration counter $k$.
In \cref{alg:global_solution_of_IOC_penalty},
we allow an additional dependence on the simplex $T$.
Indeed, if $\gamma_{k,T}$ is independent of $k$,
it is sufficient to solve the auxiliary problems
\eqref{eq:OVR_affine_penalty}
only on the new cells $T \in \TT_{k+1}  \setminus \TT_k$.
Otherwise, we would need to solve these problems
on \emph{all} cells in \emph{each} iteration.
The stopping criterion in \ref{item:A2S3}
is justified in the first part of the proof of the upcoming 
\cref{lem:penalty_convergence}.

\begin{algorithm}[t]
	\begin{enumerate}[label=\textbf{(S\arabic*)}]
		\item\label{item:A2S1}
			Let 
		$\TT_1$ 
		be a subdivision of 
		$\Bad $
		and select parameters
		$q, \rho \in (0,1)$
		and a non-decreasing
		function $P \colon \R \to \R$
		with $P(0) = 0$.
		Further, set
		$k:=1$.
			
	\item\label{item:A2S2}
			For every simplex $T \in \TT_k$,
			choose $\gamma_{k,T} > 0$
			and
			compute a global solution 
			$(\beta_{k,T},y_{k,T},u_{k,T})$ 
			of the optimization problem
			\begin{equation}
			\label{eq:OVR_affine_penalty}
			\tag{OVRP$(T, \gamma_{k,T})$}
				\begin{split}
          \min_{\beta,y,u} \quad  & F(\beta,y,u)
						                      + \gamma_{k,T}\Pen (f(\beta, y, u) - \xi_T(\beta) )\\
        \text{s.t.} \quad & \beta\, \in\, T,\\
        & 0\, =\, A y - B u,\\
        \quad  & u\, \in\, \Uad.
				\end{split}
			\end{equation}
			Select
			\begin{equation*}
				\bar{T}_{k}
					\in \argmin_{T\in \TT_{k}}
					\set*{F(\beta_{k,T},y_{k,T},u_{k,T}) 
						+ \gamma_{k,T}P\big(f(\beta_{k,T},y_{k,T},u_{k,T}) 
					-\xi_{T}(\beta_{k,T})\big)}
			\end{equation*}
			and set 
			$\byuk 
			= (\beta_{k, \bar{T}_{k}},y_{k, \bar{T}_{k}},u_{k, \bar{T}_{k}})$.

		\item\label{item:A2S3}
			Compute $\varphi(\beta_k)$. 
			If 
			$f\byuk
			=
			\varphi(\beta_k)$, 
			then $\byuk$ is a global solution of
			\eqref{eq:OVR} 
			(and, thus, of 
			\eqref{eq:UL})
			and the algorithm terminates. 
			Otherwise, we construct 
		$\TT_{k+1}$ 
		from 
		$\TT_k$ 
		by a refinement of 
		$\bar T_k$ 
		such that 
		$\vol(T) 
			\leq 
			q\cdot \vol(\bar{T}_{k})$
			and $\rho(T) \ge \rho$
		for all 
		$T\in\TT_{k+1}\setminus\TT_{k}$ .
		Set 
		$k := k+1$ 
		and go to 
		\ref{item:A2S2}.

	\end{enumerate}
	\caption{Computation of global solutions to 
					 \eqref{eq:UL}
					with penalty approach
					\label{alg:global_solution_of_IOC_penalty}
					}
\end{algorithm}

\begin{lemma}
	\label{lem:penalty_existence}
	Let the penalty function
	$\Pen \colon \R \to \R$ 
	be non-constant, non-decreasing and convex.
	Then, for every simplex $T \subset \Bad$ and $\gamma_{k,T} > 0$,
	the problem
	\eqref{eq:OVR_affine_penalty}
	possesses a solution.
\end{lemma}
\begin{proof}
	From the monotonicity and convexity of $P$, we get $\Pen(s) \to \infty$ for $s \to \infty$.
	For a minimizing sequence $\byuk$,
	the boundedness of $\beta_k$ follows from $\beta_k \in T$.
  Since $F$ is bounded from below by
  \cref{assumptions}\ref{asmp:F_f_convex} and since $\gamma_{k,T} > 0$,
	the expression
	$\Pen (f\byuk - \xi_T(\beta_k) )$ is bounded from above.
	Due to the properties of $\Pen$,
	the sequence $f\byuk$ is bounded from above.
	Thus, the boundedness of $(y_k,u_k)$ follows from 
	\cref{assumptions}\ref{asmp:f_strongly_convex}.
	Now, the remaining part of the proof is clear
	since the objective is continuous and convex,
	hence, weakly sequentially lower semicontinuous.
\end{proof}
\subsection{Standard penalization}
\label{sec:standard_penalization}
We first prove the convergence of \cref{alg:global_solution_of_IOC_penalty}
for a typical penalty function $\Pen$.

\begin{theorem}
\label{lem:penalty_convergence}
	Let the penalty function
	$\Pen \colon \R \to \R$ 
	be monotone and convex, such that
	$\Pen(s)
	=
	0$ for all $s\leq 0$
	and
	$P(s) > 0$ for all $s > 0$.
	If 
	$\gamma_{k, \bar T_k}\rightarrow \infty$, 
		\cref{alg:global_solution_of_IOC_penalty}
	either stops at a global solution of 
	\eqref{eq:OVR}
	or the computed sequence 
	$\byuk$ 
	contains a subsequence 
	converging strongly in $\R^m \times Y \times U$ to a global solution of 
	\eqref{eq:OVR}. 
	If \eqref{eq:OVR} has a unique global solution
	$\byubar$,
	then the entire sequence
	$\byuk$ 
	converges strongly to
	$\byubar$.
\end{theorem}

\begin{proof}
	A global solution $\byubar$ to
	\eqref{eq:OVR}
	is feasible for 
	\eqref{eq:OVR_affine_penalty}
	if 
	$\bar{\beta}\in T$. 
	By definition of 
	$\byuk$
  and the assumed properties for the penalty function 
  $\Pen$
	one obtains the estimate
	\begin{equation}
    \begin{aligned}
		F\byuk
		&\le
		F\byuk
		+ \gamma_{k,\bar T_k}\Pen\big(f\byuk - \xi_{\bar T_k}(\beta_{k})\big)\\
    &\le
    F(\beta_{k,T},y_{k,T},u_{k,T})
    +  \gamma_{k,T}\Pen\big(
                                f(\beta_{k,T},y_{k,T},u_{k,T}) 
                              - \xi_{T}(\beta_{k,T})
                               \big)\\
    &\le 
		F\byubar
		.
		\label{eq:UL_objective_value_estimate_for_proof_conv_alg_2}
  \end{aligned}
	\end{equation}
	If 
	\cref{alg:global_solution_of_IOC_penalty}	
	terminates in \ref{item:A2S3}, then the condition
	$f\byuk = \varphi(\beta_k)$
	implies feasibility of 
	$\byuk$ 
	for
	\eqref{eq:OVR}
	while
	\eqref{eq:UL_objective_value_estimate_for_proof_conv_alg_2}
	ensures global optimality.

	It remains to check the case that
	\cref{alg:global_solution_of_IOC_penalty}
	does not terminate.
	From \eqref{eq:UL_objective_value_estimate_for_proof_conv_alg_2}
	and
	\cref{assumptions}\ref{asmp:F_f_convex}
	we get a constant $C \ge 0$
	such that
	\begin{equation}
		\Pen\big(f\byuk - \xi_{\bar T_k}(\beta_{k})\big) 
		\leq 
		\frac{F\byubar - F\byuk}
		{\gamma_{k,\bar T_k}}
		\le
		\frac{C}
		{\gamma_{k,\bar T_k}}
		\to 0
		.
		\label{eq:P_estimate}
	\end{equation}
	Using that $P$ is non-decreasing
	and that $\xi_{\bar T_k}(\beta_k)$ is bounded from below
	(since $\varphi$ is bounded from below on $Q$),
	we get that
	$f\byuk$ is bounded from above.
	From \cref{lem:quadratic_growth_lower_level} we get
	\begin{equation*}
		f\byuk
		\ge
		f( \beta_k, y_{\beta_k}, u_{\beta_k} )
		+
		\frac\mu2 \norm{ u_k - u_{\beta_k}}_U^2
		.
	\end{equation*}
	Since $f$ is bounded from below
	and since $u_{\beta_k}$ is bounded by \cref{cor:Psi_Lipschitz},
	we obtain the boundedness of $u_k$ in $U$.
	The boundedness of the solution operator $S$
	then implies boundedness of the state
	$y_k = S u_k$ in $Y$.
	Thus, the sequence 
	$\byuk$ 
	is bounded and contains 
	a weakly convergent subsequence (without relabeling),
	$\byuk
	\rightharpoonup 
	\byuhat$. 
	The parameter
	$\beta_{k}$ 
	converges strongly because 
	$\beta\in  \Bad\subset\R^{n}$ 
	is finite dimensional.
	It remains to check optimality of the weak limit $\byuhat$
	and the strong convergence.

	From \eqref{eq:P_estimate} we obtain
	$\limsup_{k\to\infty} f\byuk-\xi_{\bar T_k}(\beta_k) \leq0$.
	Arguing as in \cref{thm:conv_alg_1},
	we obtain $\diam(\bar T_k)\rightarrow 0$.
	Together with the interpolation error estimate
	\eqref{eq:interpolation_error} we get
	\begin{equation*}
		\begin{aligned}
			0 &\leq 
			\liminf_{k\to\infty} \paren[\big]{f\byuk-\varphi(\beta_k)}
			\leq
			\limsup_{k\to\infty} \paren[\big]{f\byuk-\varphi(\beta_k)}
			\\ &\leq 
			\limsup_{k\to\infty} \paren[\big]{f\byuk-\xi_{\bar T_k}(\beta_k)
			+ C_{\rho}C_{\varphi} \diam(\bar{T}_{k})^{2}}
			=0.
		\end{aligned}
	\end{equation*}
	In particular, we have $f\byuk-\varphi(\beta_k)\to0$.
	This implies
	\begin{equation*}
		0 \leq f\byuhat - \varphi(\hat\beta)
		\leq
		\lim_{k\to\infty} \paren[\big]{f\byuk-\varphi(\beta_k)}
		=0.
	\end{equation*}
	Therefore, $\byuhat$ is feasible for \eqref{eq:OVR}
	and $f\byuk\to f\byuhat$ holds.
	Then we can argue as in \eqref{eq:from_beta_k_to_hat_beta}
	and obtain
	strong convergence for the control $u$.
	Since the solution operator $S$ is continuous,
	this proves the strong convergence of the subsequence
	$\byuk\rightarrow\byuhat$. 
	Finally, due to
	\begin{equation*}
		F\byuhat = \lim_{k\to\infty} F\byuk \leq F\byubar
	\end{equation*}
	we know that $\byuhat$ is a global minimizer of \eqref{eq:OVR}.

	Analogous to 
	\cref{thm:conv_alg_1},
	the usual subsequence-subsequence argument 
	can be used to obtain strong convergence of
	the entire sequence if the solution to
	\eqref{eq:OVR} is unique.
\end{proof}
\begin{remark}
We observe from
\cref{lem:penalty_convergence}
that it is sufficient to have the penalty parameter $\gamma_{k,T}$
being solely dependent on the simplex $T$.
A possibility is the choice
$\gamma_{k,T} = \upsilon(\diam(T))$ with a function $\upsilon$ satisfying $\upsilon(t) \to \infty$ for $t \to 0$.
A direct benefit is that the solution of the subproblem on a fixed simplex is now independent of
the iteration and only needs to be carried out once,
as in \cref{alg:global_solution_of_IOC}.
\end{remark}

\subsection{Direct penalization}
\label{sec:direct_penalization}

The problem \hyperref[eq:OVR_affine_penalty]{(OVRP($T,\gamma_{k,T}$))} can be further simplified 
if instead of a penalty function 
as described in 
\cref{lem:penalty_convergence} 
a direct penalization 
$\Pen = \Id$ 
is considered. 
It is clear that
we cannot use
$\Pen = \Id$ 
as a penalty function
for a general optimization problem.
The reason is that this function would
reward overachieving the penalized constraint.
Our constraint
$f\byu - \xi_T(\beta)$ 
cannot be arbitrarily negative
and this renders the usage of
$\Pen = \Id$ 
possible.
This choice, however, has implications 
on the choice of the penalty parameter.
The difference between
the lower-level objective functional
and the interpolation of the optimal-value function
$f\byu - \xi_T(\beta)$ can be negative.
Thus, arbitrarily increasing $\gamma$ does not work. 
The penalty parameter $\gamma$
needs to be set specifically for each simplex.

\begin{corollary}
\label{lem:OVRP_Id_solution_scheme}
We consider
\cref{alg:global_solution_of_IOC_penalty} 
with 
$\Pen = \Id$ 
and we assume that the penalty parameters
satisfy
\begin{equation*}
	\gamma_{\bar T_k} \to \infty,
	\quad
	\gamma_{\bar T_k} \diam(\bar T_k)^2 \to 0
\end{equation*}
as $k \to \infty$.
Then,
$\byuk$ 
contains a strongly convergent subsequence 
and all accumulation points are globally optimal for
\eqref{eq:OVR}. 
If  
\eqref{eq:OVR} 
admits to a unique global minimizer $\byubar$
then the entire sequence
$\byuk$
converges strongly towards this minimizer.
\end{corollary}
\begin{proof}
The argumentation follows the lines of the proof of \cref{lem:penalty_convergence}.
Therefore, we just comment on the differences.
The interpolation error estimate
\eqref{eq:interpolation_error}
allows for a lower bound for the violation of the constraint, i.e.\
\begin{equation}
	f\byuk
		- \xi_{\bar T_k}(\beta_k) 
	\geq 
	\varphi(\beta_k) 
	- \xi_{\bar T_k}(\beta_k) 
	\geq 
		- C_{\rho}C_{\varphi}
			\diam(\bar T_k)^{2}.
			\label{eq:f_xi_from_below}
\end{equation}
When using $\Pen=\Id$,
an upper bound follows as in \eqref{eq:P_estimate}
and we have
\begin{equation}
		f(\beta_k,y_k,u_k) - \xi_{\bar T_k}(\beta_{k}) 
			\leq \frac{F\byubar - F\byuk}
								 {\gamma_{\bar{T}_{k}}}
			\leq \frac{C_{F}}
								 {\gamma_{\bar{T}_{k}}}
								 .
	\label{eq:limit_varphi_constraint_penalty_P_Id}
\end{equation}
We can now argue as in \cref{lem:penalty_convergence}
and obtain $\byuk\to \byuhat$ along a subsequence, where
$\byuhat$ is a feasible point of \eqref{eq:OVR}.
In order to achieve optimality of $\byuhat$,
we combine \eqref{eq:f_xi_from_below} and 
\eqref{eq:limit_varphi_constraint_penalty_P_Id}
and obtain
\begin{equation*}
	F\byuk
	\le
	F\byubar
	+
	C_\rho C_\varphi \gamma_{\bar T_k} \diam(\bar T_k)^2
	\to
	F\byubar + 0 ,
\end{equation*}
which implies $F\byuhat=\lim_{k\to\infty} F\byuk \leq F\byubar$.
The remaining part of the proof
follows the proof of \cref{lem:penalty_convergence}.
\end{proof}

The next lemma addresses the continuous dependence
of the solution on the penalty parameter.

\begin{lemma}
\label{lem:OVRP_D_Id_continuous_dependence}
	We suppose that
$F(\cdot, S(u), u)$
is strongly convex (w.r.t. $\beta$) with constant $\mu_\beta > 0$,
independent of the control $u$.
	Then, \hyperref[eq:OVR_affine_penalty]{\textup{(OVRP($T,\gamma$))}} has a unique
	solution $\byugamma$ for all $\gamma > 0$.
Further, let
  $0<\gamma_a\leq\gamma_T<\infty$ and
	$\gamma_a \le\hat\gamma$.
	Then,
  \begin{equation*}
		\norm{\beta_{\gamma_T} - \beta_{\hat\gamma}}_{\R^n}
		+\norm{u_{\gamma_T} - u_{\hat\gamma}}_{U} \leq
      C_{\mu_\beta,\gamma_a,\gamma_T}\abs{\gamma_T-\hat\gamma}.
  \end{equation*}
\end{lemma}
\begin{proof}
The existence of a solution to
\hyperref[eq:OVR_affine_penalty]{\textup{(OVRP($T,\gamma$))}}  
follows from
\cref{lem:penalty_existence}. For
$\gamma_a\leq\gamma$,
the strong convexity of $f$
implies that the reduced objective of
\hyperref[eq:OVR_affine_penalty]{\textup{(OVRP($T,\gamma$))}}  
is
strongly convex w.r.t.\ $u$ with constant $\gamma_a \mu$
on the feasible set.
This gives
uniqueness of the state $y_\gamma=S(u_\gamma)$ and of the control
  $u_\gamma$. With
	the additional assumption on $F$,
  we get the uniqueness of $\beta_\gamma$.

  Next, we want to apply \cref{lem:continuous_dependency_on_gamma}
  to the state reduced variant of
	\hyperref[eq:OVR_affine_penalty]{\textup{(OVRP($T,\gamma$))}},
	i.e., we apply the setting
	\begin{align*}
		x &= (\beta, u), \quad
		C = T \times \Uad, \quad
		p = \gamma, \quad
		\hat \Bad = [\gamma_a,\infty), \quad
		\\
		J(x,p) &= J( (\beta, u),\gamma) := F(\beta, S(u), u) + \gamma ( f(\beta, S(u), u) - \xi_T(\beta) ).
	\end{align*}
  \cref{assumptions} ensures that the assumptions
	of
  \cref{lem:continuous_dependency_on_gamma} are satisfied.
	Thus, \cref{lem:continuous_dependency_on_gamma} implies
	\begin{equation*}
		\norm{\beta_{\hat\gamma} - \beta_{\gamma_T}}_{\R^n} + \norm{u_{\hat\gamma}-u_{\gamma_T}}_{U}
		\le
		C_{\mu_\beta, \gamma_a}
		\norm{J'_x( (\beta_{\gamma_T}, u_{\gamma_T}),\hat \gamma) - J'_x( (\beta_{\gamma_T}, u_{\gamma_T}), \gamma_T)}_{\R^n \times U\dualspace}
		.
	\end{equation*}
	Now, the derivative $J'_x( (\beta, u), \gamma) $ contains the two components
	\begin{align*}
		&F_\beta'(\beta, S(u), u) + \gamma ( f_\beta'(\beta, S(u), u) - \xi_T'(\beta )),
		\\
		&F_u'(\beta, S(u), u) + S\adjoint F_y'(\beta, S(u), u) + \gamma ( f_u'(\beta, S(u), u) + S\adjoint f_y'(\beta, S(u), u)).
	\end{align*}
	Thus, the above estimate implies
	\begin{equation*}
		\norm{\beta_{\hat\gamma} - \beta_{\gamma_T}}_{\R^n} + \norm{u_{\hat\gamma}-u_{\gamma_T}}_U
		\le
		C_{\mu_\beta, \gamma_a}
		\abs{\hat\gamma - \gamma_T}
		\parens*{
			C_{1,\gamma_T}
			+
			C_{2,\gamma_T}
		}
		,
	\end{equation*}
	with
	\begin{align*}
		C_{1, \gamma_T}
		&=
		\norm{ f_\beta'(\beta_{\gamma_T}, S(u_{\gamma_T}), u_{\gamma_T}) - \xi_T'(\beta_{\gamma_T} ) }_{\R^n}
		,
		\\
		C_{2, \gamma_T}
		&=
		\norm{ f_u'(\beta_{\gamma_T}, S(u_{\gamma_T}), u_{\gamma_T}) + S\adjoint f_y'(\beta_{\gamma_T}, S(u_{\gamma_T}), u_{\gamma_T}) }_{U\dualspace}
		.
	\end{align*}
	This shows the claim.
\end{proof}
The problem \hyperref[eq:OVR_affine_penalty]{(OVRP($T,\gamma_T$))} is a relaxation of 
\eqref{eq:OVR_affine_relaxed} and consequently the objective functional attains 
a smaller minimal value and represents 
a lower bound to the minimal objective value of
\eqref{eq:OVR_affine_relaxed}.
Since this lower bound depends on the chosen penalty parameter $\gamma_T$,
we try to adjust this parameter to obtain the largest possible lower bound.
We will now show that it is reasonable
to aim for a choice of the penalty parameter such that
the equality
$f\byugammaT = \xi(\beta_{\gamma_T})$
holds for the solution $\byugammaT$ of
\hyperref[eq:OVR_affine_penalty]{(OVRP($T,\gamma_T$))}. In the expected case where no solution to 
\eqref{eq:UP_beta} is feasible for 
\eqref{eq:OVR},
this specific penalty parameter results in the largest possible minimal objective value for
\hyperref[eq:OVR_affine_penalty]{(OVRP($T,\gamma_T$))}. 
\begin{lemma}
\label{lem:OVRP_D_Id_best_possible_gamma}
Let the state reduced
functional $F$ be strongly
convex with respect to $\beta$ with constant $\mu_\beta$ independent of the
control $u$.
Let a simplex $T$ be given and, again, $\Pen = \Id$.
Further, we assume the existence of $\beta \in T$
with $\varphi(\beta) < \xi_T(\beta)$.
For $\gamma \ge 0$, we denote a solution to \hyperref[eq:OVR_affine_penalty]{\textup{(OVRP($T,\gamma$))}} by $\byugamma$.
\begin{enumerate}[label=(\alph*)]
	\item If $f\byutilde \leq \xi_T(\tilde{\beta})$ 
		for one global solution
		$\byutilde$ to
		\hyperref[eq:OVR_affine_penalty]{\textup{(OVRP($T,0$))}}
		then
		the choice
		$\gamma_T = 0$
		yields
		the largest minimal objective value 
		for 
		\hyperref[eq:OVR_affine_penalty]{\textup{(OVRP($T,\gamma_T$))}}.
		\label{lem:opt_gamma_case_0}
		
	\item If $f\byutilde > \xi_T(\tilde{\beta})$ 
		for all global solutions 
		$\byutilde$ to
		\hyperref[eq:OVR_affine_penalty]{\textup{(OVRP($T,0$))}}
		then
		there exists 
		$\gamma_T > 0$ 
		such that
		$f\byugammaT 
		= 
		\xi_T(\beta_{\gamma_T})$
		and this choice of $\gamma_T$
		results in
		the largest minimal objective value 
		for 
		\hyperref[eq:OVR_affine_penalty]{\textup{(OVRP($T,\gamma_T$))}}.
		\end{enumerate}
	\label{lem:opt_gamma_case_hat}
\end{lemma}
The existence of $\beta \in T$
with $\varphi(\beta) < \xi_T(\beta)$
is equivalent to $\varphi$ being not affine on $T$.
Thus, this assumption is not very restrictive.
\begin{proof}
	\begin{enumerate}[label=(\alph*)]
		\item
			For any $\gamma \ge 0$
			we have
			\begin{align*}
				F\byugamma
				+ \gamma (f\byugamma - \xi_T(\beta_{\gamma})) 
					&\le
				F\byutilde 
					+ \gamma (f\byutilde - \xi_T(\tilde{\beta})) 
					\\
				&\leq 
				F\byutilde
				.
			\end{align*}
			Hence,
			the infimal value of
			\hyperref[eq:OVR_affine_penalty]{(OVRP($T,\gamma_T$))} 
			is maximized for $\gamma_T = 0$.
	
		\item 
		  \label{lem:gamma_for_f_xi_equality_b}
			We prove the existence of 
			$\gamma_T > 0$
			with
			$f\byugammaT - \xi_T(\beta_{\gamma_T}) = 0$
			by the intermediate value theorem.
			Therefore, we have to provide
			penalty parameters
			$\text{\b{$\gamma$}}{}_T, \bar\gamma_T > 0$ with
			$f\byugammaubarT - \xi_T(\beta_{\text{\b{$\gamma$}}{}_T}) \geq 0$ and
			$f\byugammabarT - \xi_T(\beta_{\bar\gamma_T}) \leq 0$.
			The required continuous dependence w.r.t.\ $\gamma > 0$
			follows from \cref{lem:OVRP_D_Id_continuous_dependence}.

			We first construct $\bar\gamma_T$.
			By assumption $F$ is bounded from below by a constant $C \in \R$
			and there exists a $\beta\in T$, such that
			$\varphi(\beta)=f\bybub<\xi_T(\beta)$.
			Thus, we can choose $\bar\gamma_T > 0$ such that
			\begin{equation}
				F\bybub 
					+ \bar\gamma_T(f\bybub - \xi_T(\beta)) \le C
					.
				\label{eq:OVRP_D_upper_bound}
			\end{equation}
			It follows that
			$f\byugammabarT - \xi_T(\beta_{\bar\gamma_T}) \leq 0$.
			
			The existence of 
			$\text{\b{$\gamma$}}{}_T$
			is proven by contradiction.
			Assume that there is no $\gamma > 0$ with 
			$f\byugamma - \xi_T(\beta_{\gamma}) \ge 0$.
			For $\gamma \searrow 0$, the bound
			$f\byugamma < \xi_T(\beta_\gamma)$ and the quadratic growth condition from 
			\cref{lem:quadratic_growth_lower_level} implies 
			boundedness of the control $u_\gamma$ whereas 
			the continuity of the solution operator 
			yields boundedness of the state $y_\gamma = Su_\gamma$. 
			The parameter $\beta_\gamma \in T$ is bounded as well.
			Thus, one obtains the existence 
			of a weak accumulation point $\byubar$ for $\gamma \searrow 0$.
			It is clear that
			$\byubar$ is feasible for
			\hyperref[eq:OVR_affine_penalty]{\textup{(OVRP($T,0$))}}
			and we show that it is even a solution.
			By optimality, we get the inequality
			\begin{equation*}
				F\byugamma 
				+ \gamma(f\byugamma - \xi_{T}(\beta_\gamma))
				\le
				F\byutilde 
				+ \gamma(f\byutilde - \xi_{T}(\tilde\beta))
			\end{equation*}
			and
			\begin{equation*}
				\lim_{\gamma\searrow 0}
				\gamma(f\byugamma - \xi_{T}(\beta_\gamma))
				=
				\lim_{\gamma\searrow 0} \gamma(f\byutilde - \xi_{T}(\tilde\beta))
				=
				0
			\end{equation*}
			follows by boundedness of $f\byugamma$.
			Thus,
			\begin{equation*}
				F\byubar
				\le
				\liminf_{\gamma\searrow 0} F\byugamma 
				\le
				F\byutilde,
			\end{equation*}
			where we take the limes inferior along the weakly convergent subsequence.
			Thus, $\byubar$ is a solution to
			\hyperref[eq:OVR_affine_penalty]{\textup{(OVRP($T,0$))}}.
			Similarly, passing to the limit inferior in
			$f\byugamma - \xi_T(\beta_{\gamma}) < 0$
			yields
			$f\byubar - \xi_T(\bar\beta) \le 0$.
			This contradicts the assumption
			and yields the existence of
			$\text{\b{$\gamma$}}{}_T$.

			By the intermediate value theorem,
			we conclude the existence of $\gamma_T > 0$
			with
			$f\byugammaT - \xi_T(\beta_{\gamma_T}) = 0$.

			It remains to prove that
			this choice of $\gamma_T$
			results in
			the largest infimal objective value 
			for 
			\hyperref[eq:OVR_affine_penalty]{\textup{(OVRP($T,\gamma_T$))}}.
			It is clear that
			$f\byugamma - \xi_T(\beta_\gamma)$
			is non-increasing w.r.t.\ $\gamma$.
			Thus,
      it follows with \cref{lem:OVRP_D_Id_continuous_dependence} that
			\begin{equation*}
				\set[\big]{
					\gamma_T \in [0,\infty) \given
					f\byugammaT - \xi_T(\beta_{\gamma_T}) = 0
				}
				=
				[\gamma_a,\gamma_b ]
				\subset
				\R_{+}.
			\end{equation*}
			For $\gamma_b < \gamma_1 < \gamma_2$,
			we have
			$f\byugammae -\xi_T(\beta_{\gamma_{1}}) 
			<
			0$
			and, thus,
			the optimality of $\byugammaz$
			for
			\hyperref[eq:OVR_affine_penalty]{(OVRP($T,\gamma_{2}$))} 
			implies
			\begin{align*}
				&F\byugammae 
					+ \gamma_{1} (f\byugammae-\xi_T({\beta_{\gamma_{1}}})) \\
				&\quad > F\byugammae 
					+ \gamma_{2} (f\byugammae-\xi_T({\beta_{\gamma_{1}}})) \\
				&\quad \geq F\byugammaz
					+ \gamma_{2} (f\byugammaz-\xi_T({\beta_{\gamma_{2}}})).
			\end{align*}
			It follows that the objective value of 
			\hyperref[eq:OVR_affine_penalty]{(OVRP($T,\gamma$))}
			is monotonically decreasing for $\gamma > \gamma_b$ 
			and, similarly, one can show that it is monotonically increasing for $\gamma < \gamma_a$
			and constant on $[\gamma_a, \gamma_b]$.
			Thus, all
			$\gamma_T\in[\gamma_a,\gamma_b]$ maximize the minimal objective value of 
			\hyperref[eq:OVR_affine_penalty]{(OVRP($T,\gamma_{T}$))}.
	\end{enumerate}
\end{proof}
In general it is not possible to check which case of
\cref{lem:OVRP_D_Id_best_possible_gamma} applies. 
However, the proof suggests that after solving 
\hyperref[eq:OVR_affine_penalty]{(OVRP($T,\gamma$))}
the value $f\byugamma - \xi_T(\beta_\gamma)$ 
can be checked to infer 
whether the choice of the penalty parameter 
$\gamma$ was adequate, 
too small or too large. 
Furthermore, when splitting the simplices in
\cref{alg:global_solution_of_IOC_penalty},
the approximation $\xi_T$ of the optimal-value function $\varphi$ cannot increase in any point $\beta\in\Bad$.
Together with the feasibility of the solution to the refined problems 
for the problem 
on the original simplex $T$, this yields that the minimal 
objective value may only remain constant or increase if the same penalty parameter 
$\gamma_T$ is used for a subproblem.
We therefore suggest starting with $\gamma = 0$ and
then using a heuristic to find a $\gamma_T$. 
The refined problems can inherit the parameter 
$\gamma_T$ as a starting point instead of zero. 
This approach covers both cases of 
\cref{lem:OVRP_D_Id_best_possible_gamma}
without the need to calculate all solutions of 
\eqref{eq:UP_beta}. Once a $\gamma_T$ is found such that $f\byugamma - \xi_T(\beta_{\gamma_T}) > 0$ one can be sure that all subproblems are of case 
\cref{lem:OVRP_D_Id_best_possible_gamma}\ref{lem:gamma_for_f_xi_equality_b}, because $\xi$ is decreasing with further refinement of the simplices.
\begin{lemma}
\label{lem:OVRP_D_Id_opt_gamma_is_multiplier_for_OVR_D}
Let the state reduced
functional $F$ be strongly
convex with respect to $\beta$ with constant $\mu_\beta$ independent of the
control $u$.
Let a simplex $T$ be given and, again, $\Pen = \Id$.
Further, we assume the existence of $\beta \in T$
with $\varphi(\beta) < \xi_T(\beta)$.
For $\gamma \ge 0$, we denote a solution of \hyperref[eq:OVR_affine_penalty]{\textup{(OVRP($T,\gamma$))}} by $\byugamma$.
Let the penalty parameter $\gamma_T$ 
be chosen as described in 
\cref{lem:OVRP_D_Id_best_possible_gamma},
i.e.,
we have one of the following cases:
\begin{enumerate}[label=(\alph*)]
	\item
		$\gamma_T = 0$ and
		$f\byutilde \leq \xi_T(\tilde{\beta})$ 
		for one global solution
		$\byutilde$ of
		\hyperref[eq:OVR_affine_penalty]{\textup{(OVRP($T,0$))}},
	\item
		$\gamma_T > 0$ and $f\byugammaT = \xi_T(\beta_{\gamma_T})$.
\end{enumerate}
Then, the point $\byutilde$ or $\byugammaT$, respectively,
is a solution of
\eqref{eq:OVR_affine_relaxed}
and
$\gamma_T$
is a multiplier corresponding to the constraint 
$f(\beta,y,u)\leq\xi_T(\beta)$
in the optimality system for 
\eqref{eq:OVR_affine_relaxed}.
\end{lemma}

\begin{proof}
First, we consider the case $\gamma_T > 0$.
Note that $\byugammaT$ is feasible for
\eqref{eq:OVR_affine_relaxed}.
We denote by $\byuT$ a solution of                      
\eqref{eq:OVR_affine_relaxed}.
Then, the optimality of both points,
$f\byuT \le \xi(\beta_T)$
and
$f\byugammaT = \xi(\beta_{\gamma_T})$
yield
\begin{align*}
	F\byugammaT
	\ge
	F\byuT
	&\ge
	F\byuT
	+
	\gamma_T ( f\byuT - \xi(\beta_T))
	\\
	&\ge
	F\byugammaT
	+
	\gamma_T ( f\byugammaT - \xi(\beta_{\gamma_T}))
	\\&=
	F\byugammaT
	.
\end{align*}
This shows
$f\byuT = \xi(\beta_{T})$
and
$F\byugammaT = F\byuT$.
Hence, the triple $\byuT$ solves 
\hyperref[eq:OVR_affine_penalty]{(OVRP($T,\gamma_T$))}
and, by the uniqueness of the solution,
the solution is $\byuT = \byugammaT$.

Thus, $\byugammaT$ is globally optimal for 
\eqref{eq:OVR_affine_relaxed}.
The optimality system of 
\hyperref[eq:OVR_affine_penalty]{\textup{(OVRP($T,\gamma_T$))}}
can be interpreted as the KKT system of 
\eqref{eq:OVR_affine_relaxed}
and the parameter $\gamma_T$
in
\hyperref[eq:OVR_affine_penalty]{\textup{(OVRP($T,\gamma_T$))}}
becomes
a Lagrange multiplier
in the KKT system of
\eqref{eq:OVR_affine_relaxed}.
Note that Lagrange multipliers
for
\hyperref[eq:OVR_affine_penalty]{\textup{(OVRP($T,\gamma_T$))}}
exist
since
the
CQ by \cite{ZoweKurcyusz1979,Robinson1976:1} is satisfied.

Finally, we consider the case $\gamma_T = 0$.
Due to $f\byutilde < \xi_T(\tilde\beta)$,
the point $\byutilde$
is feasible for
\eqref{eq:OVR_affine_relaxed}.
Since
\hyperref[eq:OVR_affine_penalty]{\textup{(OVRP($T,0$))}}
is a relaxation of
\eqref{eq:OVR_affine_relaxed},
this shows that
$\byutilde$
is a solution of
\eqref{eq:OVR_affine_relaxed}.
The interpretation of $\gamma_T$ as a multiplier
is analogous to the case $\gamma_T > 0$.
\end{proof}
This lemma shows that the problem
\eqref{eq:OVR_affine_relaxed}
is equivalent (in some sense)
to
\hyperref[eq:OVR_affine_penalty]{\textup{(OVRP($T,\gamma_T$))}}
for the
``optimal'' value of $\gamma_T$, cf.\ \cref{lem:OVRP_D_Id_best_possible_gamma}.
In the application we have in mind,
the structure of
\hyperref[eq:OVR_affine_penalty]{\textup{(OVRP($T,\gamma_T$))}}
is much nicer,
since the ``complicated'' function $f$
appears in the objective
and not in the constraints.

\section{Parameter identification in an optimal control problem}
\label{sec:para_id}
In the previous section we discussed how a global minimizer for
\eqref{eq:OVR} can be found using
\cref{alg:global_solution_of_IOC_penalty}.
However, so far we did not introduce a solution scheme for the subproblems
\eqref{eq:OVR_affine_penalty}. 
In this section we will show that
one of the main advantages 
when introducing the direct penalization
(see \cref{sec:direct_penalization})
is that
the semismooth Newton method
is applicable.
This is demonstrated
by means of
a class of example problems.

\subsection{Problem formulation and properties}
We consider the bilevel optimization problem with the 
lower-level problem
\begin{equation*}
	\tag{LL$(\alpha)$}
	\label{eq:LL_alpha}
	\begin{aligned}
		\min_{y\in H_0^1(\Omega),\,u\in L^2(\Omega)}\quad &
    		\hat{f}(\alpha,y,u)
		\coloneqq
		\sum_{i=1}^{n}
		\frac{\alpha_{i}}{2}
		\|C_{i}y-y_{d,i}\|^{2}_{L^{2}(\Omega)}
		+\frac{\sigma_{l}}{2}
		\|u\|^{2}_{L^{2}(\Omega)}
		\\
    \text{s.t.} \quad & A y - B u = 0, \\
    \quad & u \in \Uad,
	\end{aligned}
\end{equation*}
and upper-level problem
\begin{equation*}
	\tag{UL}
	\label{eq:UL_alpha}
	\begin{aligned}
		\min_{\alpha\in\R^n,\,y\in H_0^1(\Omega),\,u\in L^2(\Omega)} \quad &
		\hat{F}(\alpha,y,u) 
		\coloneqq 
		\frac{1}{2}
		\|y-y_{m}\|_{L^{2}(\Omega)}^{2}
		+\frac{\sigma_{u}}{2}
		\|u-u_{m}\|_{L^{2}(\Omega)}^{2}
		+\frac{\sigma_{\alpha}}{2}
		\|\alpha\|_{\R^{n}}^{2}
		\\
    \text{s.t.} \quad & \alpha \in Q_\alpha,
		\\
		\quad & (y,u)\;\text{solves}\;\eqref{eq:LL_alpha}.
	\end{aligned}
\end{equation*}
As an underlying assumption let 
$\sigma_u,\sigma_l,\sigma_{\alpha} > 0$,
$y_m,y_{d,i},\in L^2(\Omega)$, where $\Omega\subset \R^l$ is an open and bounded set.
Moreover, let $Q_\alpha := [a_1,b_1]\times\dots\times [a_n,b_n]$
constitute a box constraint on $\alpha$,
where $a_i,b_i\in\R$ satisfies $0<a_i<b_i$ for all $i\in\set{1,\ldots,n}$.
We also require that the admissible set $\Uad$ has the structure
$\Uad=\{v \in L^2(\Omega) \mid 
u_a \leq v \leq u_b
\;
\text{a.e.\ in}\;\Omega\}$,
where 
$u_a,u_b \in L^2(\Omega)$ 
are functions
such that 
$\Uad$
is nonempty.
Further, let
$A:H_0^1(\Omega) \to H^{-1}(\Omega)$,
$B:L^2(\Omega) \to H^{-1}(\Omega)$,
$C_i:H_0^1(\Omega) \to L^2(\Omega)$
be bounded linear operators such that $A$ is bijective.

We also assume that $B$ can be extended to
an operator $B \in L[L^q(\Omega), H^{-1}(\Omega)]$
for some $q\in (1,2)$.
Additionally, we require $u_m,u_a,u_b\in L^{q'}(\Omega)$,
where $q'>2$ satisfies $1/q+1/q' = 1$.

We observe that the lower-level objective functional $\hat f$
is not convex with respect to all variables.
In particular, \cref{assumptions}\ref{asmp:F_f_convex} is not satisfied.
Additionally, the corresponding optimal-value function
is usually not convex either.
As \cref{alg:global_solution_of_IOC}
depends on convexity 
of the optimal-value function 
one has to first transform the problem 
in such a way that the new lower-level objective functional is convex.
For this purpose, we consider the simple substitution 
$\beta_i=1/\alpha_i$ .
We also define $\sigma_\beta := \sigma_\alpha$.
For the upper-level objective this substitution results in
\begin{equation*}
	F\byu
	\coloneqq
	\frac12\norm{y-y_m}_{L^2(\Omega)}^2
	+ \frac{\sigma_u}2 \norm{u-u_m}_{L^2(\Omega)}^2
	+ \frac{\sigma_\beta}{2} \sum_{i=1}^n \left(\frac1{\beta_i}\right)^2.
\end{equation*}
The constraint $\alpha\in Q_\alpha$ has to be transformed
to $\beta\in \Bad := [b_1^{-1}, a_1^{-1}] \times \cdots \times [b_n^{-1}, a_n^{-1}]$.
Observe that $Q$ is a compact subset of $(0,\infty)^n$
because $Q_\alpha$ is a compact subset of $(0,\infty)^n$.

One can check that $F$ is convex on 
$\Bad\times H_0^1(\Omega) \times L^2(\Omega)$
due to $\beta > 0$ for $\beta \in \Bad$.
The transformed lower-level objective is
\begin{equation}
	\label{eq:f_beta_def}
	f\byu 
	\coloneqq
	\sum_{i=1}^{n}
	\frac{1}{2\beta_{i}}
	\|C_{i}y-y_{d,i}\|^{2}_{L^{2}(\Omega)}
	+ \frac{\sigma_{l}}{2}
	\|u\|_{L^{2}(\Omega)}^{2}.
\end{equation}
We check that this $f$
is indeed convex on $Q\times H_0^1(\Omega)\times L^2(\Omega)$.
Here we use that 
for a Banach space $Y$,
the function 
$g:Y\to\R$, $y\mapsto \tfrac12\norm{y}^2_Y$ is convex and
for $\lambda > 0$ the so-called perspective of $g$ is given by
\begin{equation}
	Y \times ( 0,\infty) \ni (y,\lambda) \mapsto \lambda g(y / \lambda) 
	= \frac{1}{2 \lambda}\norm{y}^{2}_{Y}.
\end{equation}
It is known that the perspective of a convex function is convex
(e.g.\ one can simply generalize the proof of 
\cite[Lemma 2]{DacorognaMarechal2008} to Banach spaces).
Now convexity is preserved under composition with an affine function 
$y\mapsto Cy - y_d$. 
Thus, the function
$(\beta_i,y)\mapsto \frac1{2\beta_i}\norm{C_i y - y_{d,i}}_{L^2(\Omega)}^2$
is convex.
The convexity of $f$ follows.

With the above setting and observations, one can show that
the transformed problem satisfies \cref{assumptions}.

\subsection{Stationarity system for the direct penalization}
\label{sec:stationarity_system}
Classic choices
of the penalty function for 
\eqref{eq:OVR_affine_penalty},
e.g., $P = \max(0,\cdot)^2$,
will result in subproblems that are difficult to handle.
In particular,
the optimality system cannot be reformulated
as a simple projection formula.
We will see that
the direct penalization $P = \Id$ results 
in an easy to implement solution algorithm 
for 
\eqref{eq:OVR_affine_penalty}.
Computing the solution of 
\eqref{eq:OVR_affine_penalty}
requires the construction of 
$\xi$ 
and thereby the evaluation of 
$\varphi(\beta)$ 
at certain points. 
This equates to solving
single-level optimal control problems.

In order to state the stationarity conditions,
we first reformulate the condition $\beta\in T$. 
Recall that $T$ is a (non-degenerate) simplex.
Thus,
$T$
can be written as the intersection of $n+1$ half-spaces,
$T = \set{\beta\in \R^n \given K_T \beta\leq b_T}$, 
where $K_T\in\R^{(n+1)\times n}$ is a suitable matrix.
Clearly,
at most $n$ of these constraints may simultaneously hold 
with equality and that all those constraints that are satisfied 
with equality are linearly independent.
Thus, \eqref{eq:OVR_affine_penalty}
with $P = \Id$
takes the form
\begin{align*}
  \min_{\beta,y,u}  \quad  & F(\beta,y,u)
	+ \gamma_{k,T}\big(f(\beta, y, u)
	- \xi_T(\beta)\big)\\
  \text{s.t.} \quad  
	& K_T\beta -b_T \leq 0, \\
  & A y - B u = 0, \\
	\quad  & u  \in \Uad.
\end{align*}
The KKT system for 
\eqref{eq:OVR_affine_penalty}
with direct penalization
($P = \Id$)
is given by
\begin{subequations}
\label{eq:KKT_OVRP_T}
\begin{align}
	0 &= F'_\beta\byu + \gamma_{k,T} (f'_\beta\byu - a_T\beta) + K_T^\top z,
    \label{eq:KKT_OVRP_T_1}
	\\
	0 &=  F'_y\byu + \gamma_{k,T} f'_y\byu + A\adjoint p, \label{eq:KKT_OVRP_T_2}
  	\\
		0 &= F'_u\byu + \gamma_{k,T} f'_u\byu - B\adjoint p + \nu, \label{eq:KKT_OVRP_T_3}\\
  0 &= Ay - Bu, \label{eq:KKT_OVRP_T_4}\\
  z &\ge 0 \wedge K_T\beta - b_T\le 0\wedge z^\top(K_T\beta - b_T) = 0,
  \label{eq:KKT_OVRP_T_5}\\
  \nu &\in \NN_{\Uad}(u),\qquad u\in \Uad,\label{eq:KKT_OVRP_T_6}
\end{align}
\end{subequations}
where $p\in {H^{1}_{0}(\Omega)}$, $z\in \R^{n+1}$,
and $\nu\in L^2(\Omega)$ 
are the Lagrange multipliers.
The vector $a_{T}$ 
refers to the derivative of 
the affine function
$\xi_T$ on the simplex $T$. 
\begin{lemma}
	\label{lem:Robinson_regularity_penalty_problem}
	The feasible point $\byu$ is a local/global solution to
	\eqref{eq:OVR_affine_penalty}
	if and only if
	there exist multipliers 
	$p\in {H^{1}_{0}(\Omega)}$, $z\in \R^{n+1}$,
	and $\nu\in L^2(\Omega)$ 
	such that \eqref{eq:KKT_OVRP_T} holds.

	The solution and the corresponding multipliers 
	are 
	unique.
\end{lemma}
\begin{proof} 
	``$\Rightarrow$'': We check that the
	Robinson regularity condition for 
	the reformulated problem
	is satisfied.
	This condition reads
	\begin{equation*}
		\begin{bmatrix}
				A&-B&0\\
				0&0&K_T
		\end{bmatrix}
		\begin{pmatrix}
				H^{1}_{0}(\Omega)\\
				\mathcal{R}_{\Uad}(u)\\
				\R^{n}
		\end{pmatrix}
		-
		\begin{pmatrix}
				\{0\}\\
				\operatorname{cone}( \R_{-}^{n+1} - (K_T\beta - b_T))
		\end{pmatrix}
		=
		\begin{pmatrix}
				H^{-1}(\Omega)\\
				\R^{n+1}
		\end{pmatrix}
		.
	\end{equation*}
	The two lines of the equation
	are independent of each other.
	By assumption, $A$ is bijective,
	i.e.,
	$A(H^{1}_{0}(\Omega)) 
	= 
	H^{-1}(\Omega)$.
	For the second line we recall that
	the Robinson regularity condition
	is equivalent to the
	Mangasarian--Fromovitz condition
	for standard nonlinear optimization problems,
	see \cite[p.~71]{BonnansShapiro2000}.
	Thus, the second line is satisfied
	since we have assumed that the simplex $T$
	is non-degenerate,
	i.e., we even have the
	linear-independence constraint qualification
	for the system $K_T \beta \le b_T$.
	This shows the existence of multipliers,
	see
	\cite[Theorem~3.9]{BonnansShapiro2000}.

	``$\Leftarrow$'': This is clear since
	\eqref{eq:OVR_affine_penalty}
	is a convex problem.

	It remains to address the uniqueness.
	The uniqueness of the solution follows from the strict convexity
	of the objective.
	The second line of the KKT system gives uniqueness of the adjoint $p$,
	since $A$ is an isomorphism.
  Similarly one gets uniqueness of $\nu$ from the third line. 
  Regarding uniqueness of $z$ we observe that the matrix $K_T$ describing a non-degenerate
  simplex has rank $n$, even after removing
  an arbitrary line. 
	Additionally, there exists at least one inactive
  constraint, such that $z$ is equal zero in this component. After removing the
  corresponding component from $z$ and the respective column from $K^\top_T$ in
	the first line of \eqref{eq:KKT_OVRP_T},
  $z$ is obtained by inverting a square matrix of full rank. Thus, $z$ is
  unique.
\end{proof}
We introduce two auxiliary functions 
$h,\hat h:(0,\infty)^n\times H_0^1(\Omega)\to\R$
via
\begin{equation}
	\begin{aligned}
		\hat h(\beta,y) &\coloneqq
		\frac12\|y-y_m\|^2_{L^2(\Omega)} + 
		\gamma_{k,T}\left(
			\sum_{i=1}^{n}
			\frac{1}{2\beta_{i}}
			\|C_{i}y-y_{d,i}\|^{2}_{L^{2}(\Omega)}
			-\xi_T(\beta)
		\right),\\
		h(\beta,y) &\coloneqq  \hat h(\beta,y)
		+\frac{\sigma_\beta}
		{2}
		\sum_{i=1}^n 
		\left(
			\frac1{\beta_i}
		\right)^2.
	\end{aligned}
	\label{eq:h_definition}
\end{equation}
Note that $h$ represents the part of the objective function
of \eqref{eq:OVR_affine_penalty}
that does not depend on $u$.

Recall that 
$K_T\in\R^{(n+1)\times n}$,
$A:H_0^1(\Omega)\to H^{-1}(\Omega)$,
$B:L^2(\Omega)\to H^{-1}(\Omega)$
are bounded linear operators
and that 
$A$ 
is invertible. We define the function\\
$W:(0,\infty)^n
\times 
H_0^1(\Omega)
\times 
L^2(\Omega)
\times 
\R^{n+1}
\times H_0^1(\Omega)
\to
\R^{n}
\times
H^{-1}(\Omega)
\times
L^2(\Omega)
\times
\R^{n+1}
\times
H^{-1}(\Omega)
$
via
\begin{equation}
	W(\beta,y,u,z,p)
	\coloneqq
	\begin{pmatrix}
		h'_\beta(\beta,y) + K^\top_T z
		\\
		h'_y(\beta,y) + A\adjoint p
		\\
		u - \min\bigl(\max((B\adjoint p+\sigma_{u}u_{m}) / \hat\sigma,u_a),u_b\bigr)
		\\
		\max(K_T\beta-b_T,-z)
		\\
		Ay-Bu
	\end{pmatrix}
	.\label{W_definition}
\end{equation}
with $\hat\sigma \coloneqq \sigma_u + \gamma_{k,T}\sigma_l$.
Now we discuss the relation between the roots of $W$ and the optimality system. 
\begin{lemma}
	\label{lem:optimality_via_W_new}
	Let 
  $\beta\in T$, 
	$y\in H_0^1(\Omega)$, 
	$u\in L^2(\Omega)$
	be given.
	Then 
	$\byu $ 
	is the solution of 
	\eqref{eq:OVR_affine_penalty}
	if and only if there exist 
	$z\in\R^{n+1}$, 
	$p\in H_0^1(\Omega)$
	such that $W(\beta,y,u,z,p)=0$ with $h$ as defined in 
	\eqref{eq:h_definition}. 
\end{lemma}
\begin{proof}
	In view of \cref{lem:Robinson_regularity_penalty_problem},
	we have to check that \eqref{eq:KKT_OVRP_T}
	is equivalent to $W(\beta,y,u,z,p)=0$.

	It is clear that
	\eqref{eq:KKT_OVRP_T_1}, \eqref{eq:KKT_OVRP_T_2} and \eqref{eq:KKT_OVRP_T_4}
	are equivalent to
	lines 1, 2 and 5
	in \eqref{W_definition}.
	The complementarity conditions \eqref{eq:KKT_OVRP_T_5} on $z$ and $b_T - K_T \beta$
can be reformulated
via
\begin{align*}
	z&\geq 0,\quad b_T-K_T\beta\geq 0,\quad z^\top(b_T-K_T\beta)=0\\
	\Longleftrightarrow\quad 0 &=\min\paren{z,b_T-K_T\beta}
	\quad\Longleftrightarrow\quad 0 =\max\paren{-z,K_T\beta-b_T}.
\end{align*}

	A similar reformulation is standard for treating
	the gradient equation \eqref{eq:KKT_OVRP_T_3}
	in combination with the inclusion \eqref{eq:KKT_OVRP_T_6},
	see
	\cite[Theorem~2.28]{Troeltzsch2009}. 
	These two equations are equivalent to the projection formula
	\begin{equation*}
		u
		=
		\Proj_{\Uad}\parens[\big]{(B\adjoint p+\sigma_{u}u_{m}) / \hat\sigma}
		=
		\min\bigl(\max((B\adjoint p+\sigma_{u}u_{m}) / \hat\sigma,u_a),u_b\bigr),
	\end{equation*}
	i.e., line 3 in \eqref{W_definition}.
	Note that $\nu$ does not appear in \eqref{W_definition},
	but it is uniquely determined by \eqref{eq:KKT_OVRP_T_3}.
	This shows that the KKT system
	is equivalent to
	$W(\beta,y,u,z,p)=0$.
	This finishes the proof.
\end{proof}

\subsection{Semismooth Newton method for the subproblems}
\label{sec:semismooth_newton}
We have shown in 
\cref{lem:optimality_via_W_new}
that we can characterize the solution of the subproblem 
\eqref{eq:OVR_affine_penalty} with the nonlinear operator $W$.
An established way to solve problems with this structure
is
the semismooth Newton method,
cf.\ \cite{HintermuellerItoKunisch2002}. 
To this end, we verify 
the Newton differentiability of 
$W$ 
and the invertibility of the Newton matrix.
In order to
state the Newton derivative of $W$,
we need to define some index sets
and corresponding operators. 
We define
\begin{equation*}
	\begin{aligned}
		\AA_1(\beta,z)
		&:= \set{i\in\set{1,\ldots,n+1}\given (K_T\beta-b_T)_i \geq -z_i},
		\\
		\AA_2(\beta,z)
		&:= \set{i\in\set{1,\ldots,n+1}\given (K_T\beta-b_T)_i < -z_i},
		\\
		\AA_3(p)
		&:= \set{u_a \leq (B\adjoint p+\sigma_{u}u_{m} )/\hat\sigma \leq u_b} \subset\Omega.
	\end{aligned}
\end{equation*}
and for $i\in\set{1,2}$ we write 
$\chi_{\AA_i(\beta,z)} \in \R^{(n+1)\times (n+1)}$
for the diagonal matrix 
that whose 
$k$-th 
diagonal entry is
$1$ 
if 
$k\in\AA_i(\beta,z)$ 
and 
$0$
otherwise.
Similarly,
we write 
$\chi_{\AA_3(p)}:
L^2(\Omega)
\to 
L^2(\Omega)$
for the multiplication operator corresponding to the
characteristic function of $\AA_3(p)$
on the space 
$L^2(\Omega)$.

\begin{lemma}
\label{lem:W_semismooth}
The mapping 
$W$ 
is Newton differentiable
and
a Newton derivative
of $W$
at a point 
$(\beta,y,u,z,p)$ 
is given by
the block operator
\begin{equation*}
	W'(\beta,y,u,z,p)
	=
	\begin{bmatrix}
		h''_{\beta\beta}(\beta,y) & h''_{\beta y}(\beta,y) & 0 & K_T^\top & 0
		\\
		h''_{y\beta}(\beta,y) & h''_{yy}(\beta,y) & 0 & 0 & A\adjoint
		\\
		0 & 0 & \Id & 0 & -\hat\sigma^{-1}\chi_{\AA_3(p)}B\adjoint 
		\\
		\chi_{\AA_1(\beta,z)}K_T & 0 & 0 & -\chi_{\AA_2(\beta,z)} & 0
		\\
		0 & A & -B & 0 & 0
	\end{bmatrix}
	.
\end{equation*}
\end{lemma}
\begin{proof}
	To show Newton differentiability of $W$, one has to pay
  attention only to the third and fourth line as the others are Fréchet
  differentiable. For the fourth line one can use that in finite dimensions the
  composition of Newton differentiable functions is Newton differentiable
  cf.\ \cite[Proposition~2.9]{Ulbrich2011} 
	and combine this with the fact that $\max(\cdot,\cdot)$ 
	is Newton differentiable (see \cite[Proposition~2.26]{Ulbrich2011}).
	Furthermore, \cite[Theorem~3.49]{Ulbrich2011}
	can be used to show the Newton differentiability of the third line:
	If we use $m=3$, $\psi(s)=\min(\max(s_1,s_2),s_3)$, $r=r_i=2$,
	$G(p)= \paren{(B\adjoint p+\sigma_u u_m)/\hat\sigma , u_a,u_b}$
	in the setting of \cite[Section~3.3]{Ulbrich2011},
	then the required \cite[Assumption~3.32]{Ulbrich2011} is satisfied with $q_i=q'>2$,
	by the higher regularity
	$B\adjoint \in L[H_0^1(\Omega), L^{q'}(\Omega)]$.

	Consequently, the function 
	$H_0^1(\Omega)\ni p\mapsto
	\min\paren{\max\paren{(B\adjoint p+\sigma_u u_m)/\hat\sigma , u_a},u_b}$
	is Newton differentiable.

	Now a Newton derivative can be obtained using direct calculations and
	utilizing the index sets that are introduced above.
\end{proof}
The proof required a norm gap,
which was ensured by the higher regularity $B\adjoint \in L[H_0^1(\Omega),
L^{q'}(\Omega)]$ with $q' > 2$, which is intrinsic to our problem setting. 
This allowed us to prove the
Newton differentiability of $W$ in the spaces where $W$ is defined.
In particular when                 
adapting the Algorithm from \cite[Algorithm~3.10]{Ulbrich2011},
see \cref{alg:semismooth_Newton_OVRP},
this allows for the smoothing step
to be skipped.
This smoothing step is designed to treat the more general
case when Newton differentiability can only be shown by artificially introducing
a norm gap while the boundedness of the inverse of the derivative can only be
shown in the original setting 
(cf. \cite[Introduction to section 3]{Ulbrich2011}).
Note that \ref{item:A3S2} is well defined as long as
$\beta_i$ is positive, since the function $W$ is only defined for positive $\beta$.
This, however, does not influence the local convergence of \cref{alg:semismooth_Newton_OVRP}.
\begin{algorithm}[h]
	\begin{enumerate}[label=\textbf{(S\arabic*)}]
		\item\label{item:A3S0}
			Choose an initial point $(\beta_0,y_0,u_0,z_0,p_0) \in (0,\infty)^n
\times 
H_0^1(\Omega)
\times 
L^2(\Omega)
\times 
\R^{n+1}
\times H_0^1(\Omega)$ and set $i=0$
		\item\label{item:A3S1}
			If $W(\beta_i,y_i,u_i,z_i,p_i) = 0$, then STOP
		\item\label{item:A3S2}
			Compute $s_i$ from
			\begin{equation*}
				W'(\beta_i,y_i,u_i,z_i,p_i)s_i = - W(\beta_i,y_i,u_i,z_i,p_i)
			\end{equation*}
		\item\label{item:A3S3}
		Set $(\beta_{i+1},y_{i+1},u_{i+1},z_{i+1},p_{i+1}) = 
    (\beta_i,y_i,u_i,z_i,p_i) + s_i$, increment $i$ by one, 
		and go to step \ref{item:A3S1}
	\end{enumerate}
	\caption{Semismooth Newton method for \eqref{eq:OVR_affine_penalty}}
	\label{alg:semismooth_Newton_OVRP}
\end{algorithm}

To prove fast convergence of the semismooth Newton method, 
the uniform invertibility of the Newton derivative
$W'(\beta,y,u,z,p)$ is needed.
For this purpose, we convert the
Newton derivative $W'(\beta,y,u,z,p)$ 
into a self-adjoint operator,
since the latter type 
of operator is easier to handle.
For that purpose we fix a point
$(\beta,y,u,z,p)$.
We use the notation 
$I_1\in\R^{(n+1)\times l_1}$,
$I_2\in\R^{(n+1)\times (n+1 - l_1)}$,
$I_3: L^2(\AA_3(p))\to L^2(\Omega)$,
$I_4: L^2(\Omega\setminus \AA_3(p))\to L^2(\Omega)$,
to refer to the canonical embedding operators
that correspond to the index sets
$\AA_1(\beta,z)$, 
$\AA_2(\beta,z)$, 
$\AA_3(p)$, 
$\Omega\setminus\AA_3(p)$.
Here $l_1$ denotes the cardinality of $\AA_1(\beta, z)$.
We mention that 
$I_1^\top,
I_2^\top,
I_3\adjoint,
I_4\adjoint$
are the corresponding restriction operators 
and, consequently,
\begin{equation*}
	\begin{aligned}
		\chi_{\AA_1(\beta,z)}&= I_1I_1^\top,
		&
		\chi_{\AA_2(\beta,z)}&= I_2I_2^\top,
		&
		\chi_{\AA_3(p)}&= I_3I_3\adjoint,
		\\
		\Id_{\R^{n+1}}&=I_1I_1^\top + I_2I_2^\top,
		&
		\Id_{L^2(\Omega)}&=I_3I_3\adjoint + I_4I_4\adjoint.
	\end{aligned}
\end{equation*}
We define the linear operator
$
	\hat W'$
	from
	$
				\R^n
				\times 
				H_0^1(\Omega)
				\times 
				L^2(\AA_3(p))
				\times 
				\R^{l_1}
				\times
				H^1_0(\Omega)
				$
				to
				$
				\R^n
				\times 
				H^{-1}(\Omega)
				\times 
				L^2(\AA_3(p))
				\times 
				\R^{l_1}
				\times 
				H^{-1}(\Omega)
				$
via
\begin{equation*}
	\hat W' :=
	\begin{bmatrix}
		h''_{\beta\beta}(\beta,y) & h''_{\beta y}(\beta,y) & 0 
		& K_T^\top I_1 & 0
		\\
		h''_{y\beta}(\beta,y) & h''_{yy}(\beta,y) & 0 & 0 & A\adjoint
		\\
	0 & 0 & \hat\sigma \Id & 0 & -(B I_3)\adjoint 
		\\
		I_1^\top K_T & 0 & 0 & 0 & 0
		\\
		0 & A & -B I_3 & 0 & 0
	\end{bmatrix}
	.
\end{equation*}
It can be seen that 
$\hat W'$ 
is self-adjoint.
Note that the spaces 
on which 
$\hat W'$ 
operates depend on 
$\beta,z,p$.
The next 
\lcnamecref{lem:symmetrization}
gives us a relation between 
$\hat W'$
and 
$W'(\beta,y,u,z,p)$.

\begin{lemma}
	\label{lem:symmetrization}
	Let 
	$(\beta,y,u,z,p)
	\in 
	(0,\infty)^n
	\times 
	H_0^1(\Omega)
	\times 
	L^2(\Omega)
	\times 
	\R^{n+1}
	\times 
  H_0^1(\Omega)$ 
	be fixed.
	Furthermore, let two points
	$(\beta_1,y_1,u_1,z_1,p_1)
	\in 
	\R^n
	\times 
	H_0^1(\Omega)
	\times 
	L^2(\Omega) 
	\times 
	\R^{n+1}
	\times 
  H_0^1(\Omega)$ 
	and
	$(\beta_2,y_2,u_2,z_2,p_2)
	\in 
	\R^n
	\times 
	H^{-1}(\Omega)
	\times 
	L^2(\Omega) 
	\times 
  \R^{n+1}
	\times 
	H^{-1}(\Omega)$
	be given.
	Then
	
	\begin{equation}
		W'(\beta,y,u,z,p)
		\begin{pmatrix}
			\beta_1 \\ y_1 \\ u_1 \\ z_1 \\ p_1
		\end{pmatrix}
		=
		\begin{pmatrix}
			\beta_2 \\ y_2 \\ u_2 \\ z_2 \\ p_2
		\end{pmatrix}
		\label{eq:W'_byuzp}
	\end{equation}
	holds if and only if
	\begin{equation}
		\begin{aligned}
			\hat W'	\begin{pmatrix}
			\beta_1 \\ y_1 \\ I_3\adjoint u_1 \\ I_1^\top z_1 \\ p_1
		\end{pmatrix}
		&=
		\begin{pmatrix}
			\beta_2+ K_T^\top I_2I_2^\top z_2
			\\ y_2 \\ \hat\sigma I_3\adjoint u_2 
			\\ I_2^\top z_2 \\ p_2 + BI_4I_4\adjoint u_2
		\end{pmatrix},
		\qquad
		\begin{aligned}
			I_4\adjoint u_1 &= I_4\adjoint u_2,
			\\
			-I_2^\top z_1 &= I_2^\top z_2
		\end{aligned}
		\end{aligned}
		\label{eq:hatW'_byuzp}
	\end{equation}
	hold.
\end{lemma}
\begin{proof}
The proof can be carried out by direct calculation. 
We first assume 
\eqref{eq:W'_byuzp} to be valid.
Computing the application of $\hat W'$ yields
\begin{equation*}
		\hat W'
		\begin{pmatrix}
			\beta_1 \\ y_1 \\ I_3\adjoint u_1 \\ I_1^\top z_1 \\ p_1
		\end{pmatrix}
		=
		\begin{pmatrix}
		h''_{\beta\beta}(\beta,y)\beta_{1} 
			+ h''_{\beta y}(\beta,y)y_{1}
			+ K_T^\top I_{1}I_{1}^\top z_{1} \\
	  y_{2} \\
		I\adjoint _{3}(\hat\sigma u_{1}
			 -B\adjoint p_{1}) \\
		I_{1}^\top K_T\beta_{1} \\
		Ay_{1}  -BI_{3}I\adjoint _{3}u_{1}
		\end{pmatrix}
		.
		\end{equation*}
		We use the definition of the index sets and receive the equivalent expression
		\begin{equation}
			\label{eq:hat_W'_complicated}
		\hat W'
		\begin{pmatrix}
			\beta_1 \\ y_1 \\ I_3\adjoint u_1 \\ I_1^\top z_1 \\ p_1
		\end{pmatrix}
		=
		\begin{pmatrix}
			h''_{\beta\beta}(\beta,y)\beta_{1} 
			+ h''_{\beta y}(\beta,y)y_{1}
			+ K_T^\top z_{1} - K_T^\top I_{2}I_{2}^\top z_{1} \\
			h''_{y\beta}(\beta,y)\beta_1
			+ h''_{y y}(\beta,y)y_{1} + A\adjoint p_1 \\
			\hat\sigma I_{3}\adjoint (u_{1} - \hat\sigma^{-1}\chi_{\AA_3(p)}B\adjoint p_{1}) \\
			I^\top _{1}( \chi_{\AA_1(\beta,z)} K_T\beta_{1} 
			- \chi_{\AA_2(\beta,z)} z_1)\\
			Ay_{1} -Bu_{1} + BI_{4}I\adjoint _{4}u_{1}
		\end{pmatrix}
		,
		\end{equation}
		where we used 
		$I_1I_1^\top = \Id_{\R^{n+1}} - I_2I_2^\top$,
		$I_3\adjoint = I_3\adjoint \chi_{\AA_3(p)}$,
		$I_1\adjoint = I_1\adjoint \chi_{\AA_1(\beta,z)}$,
		$I^\top_1 \chi_{\AA_2(\beta,z)} = 0$,
		and
		$I_3I_3^\top = \Id_{L^2(\Omega)} - I_4I_4\adjoint$.
		Using the description of $W'(\beta,y,u,z,p)$ yields
		\begin{equation}
			\label{eq:hat_w'_intermediate}
		\hat W'
		\begin{pmatrix}
			\beta_1 \\ y_1 \\ I_3\adjoint u_1 \\ I_1^\top z_1 \\ p_1
		\end{pmatrix}
		=
		\begin{pmatrix}
		\beta_{2} - K_T^\top I_{2}I_{2}^\top z_{1} \\
	   y_{2} \\
		I\adjoint _{3}u_{2}\\
		I_{1}^\top z_{2}\\
		p_{2} + BI_{4}I\adjoint _{4}u_{1}
		\end{pmatrix}.
	\end{equation}
	Note that the claimed relations
	$I_4\adjoint u_1 = I_4\adjoint u_2$
	and
	$-I_2^\top z_1 = I_2^\top z_2$
	follow from the equations
	$\Id u_1 + \hat\sigma ^{-1} \chi_{\AA_3(p)} BG\adjoint p_1 = u_2$
	and
	$\chi_{\AA_1(\beta,z)} K_T\beta_1 - \chi_{\AA_2} z_1 = z_2$
	(which are part of \eqref{eq:W'_byuzp}).
	With these relations, we directly get \eqref{eq:hatW'_byuzp}
	from \eqref{eq:hat_w'_intermediate}.

	For the other direction, we
	first get \eqref{eq:hat_w'_intermediate} directly from \eqref{eq:hatW'_byuzp}.
	Then, a comparison with \eqref{eq:hat_W'_complicated} yields
	the equations for $\beta_2$, $y_2$, $p_2$,
	and
	$I_3\adjoint u_1 + \hat\sigma ^{-1} \chi_{\AA_3(p)} BG\adjoint p_1 = I_3\adjoint u_2$,
	$I_1^\top (\chi_{\AA_1(\beta,z)} K_T\beta_1 - \chi_{\AA_2} z_1) = I_1^\top z_2$.
	The final expression \eqref{eq:W'_byuzp} follows by utilizing
	$I_4\adjoint u_1 = I_4\adjoint u_2$
	and $-I_2^\top z_1 = I_2^\top z_2$
	again.
\end{proof}

In order to ensure the uniform invertibility
of the operators $\hat W'$,
we state an auxiliary lemma.
\begin{lemma}
	\label{lem:saddle_point_matrix}
	Let $X,Y$ be Hilbert spaces
	and $\hat A:X\to X\dualspace$, $\hat B:X\to Y\dualspace$
	be bounded linear operators.
	Let the bounded linear operator 
	$\hat D \colon X\times Y\to X\dualspace \times Y\dualspace$ be defined via
	\begin{equation*}
		\hat D=
		\begin{bmatrix}
			\hat A & \hat B\adjoint \\ \hat B & 0
		\end{bmatrix}
		.
	\end{equation*}
	Suppose that $\hat B$ is surjective
	and that $\hat A$ is coercive on $\ker \hat B$, i.e.\ there exists
	a constant $\hat \gamma > 0$
	such that
	$\dual{\hat Ax}{x}\geq\hat\gamma\norm{x}_X^2$ for all $x\in\ker \hat B$.

	Then $\hat D$ is continuously invertible.
	Moreover, the estimate
	\begin{equation*}
		\norm{\hat D^{-1}}\leq 4c^5
	\end{equation*}
	holds,
	where $c:=\max(1,\hat\gamma^{-1},\alpha,\norm{\hat A})$,
	$\alpha>0$ is a constant 
	such that $B^1_{Y\dualspace}(0)\subset \hat B(B^\alpha_X(0))$,
	and $\hat\gamma>0$ is the coercivity constant from above.
\end{lemma}
\begin{proof}
	This result follows from \cite[Proposition~II.1.3]{BrezziFortin1991}.
	Note that we have
	$\hat B(X)=Y\dualspace$ and $\ker \hat B\adjoint=\set0$.
\end{proof}

\begin{lemma}
	\label{lem:newton_derivative_invertible}
	Let
	$(\beta,y,u,z,p)\in (0,\infty)^n\times H_0^1(\Omega)\times L^2(\Omega)
  \times \R^{n+1}\times H_0^1(\Omega)$ 
	be fixed.
	Suppose that
	$I_1^\top K_T\in \R^{l_1\times n}$ is surjective,
	i.e.\ that the rows of $K_T$ which correspond to
	the index set $\AA_1(\beta,z)$ are linearly independent.
	Then, the operator $W'(\beta,y,u,z,p)$ is continuously invertible.
	Moreover, we have
	$\norm{W'(\beta,y,u,z,p)^{-1}}\leq C$
	for a constant $C>0$,
	which does not depend on $\beta,y,u,z,p$ but can depend
	on an upper bound of $\norm{y}$, on the upper and lower bounds of $\beta$,
	and on $K_T,A,B,h,\hat\sigma,\sigma_r$.
\end{lemma}
\begin{proof}
	We start with showing that
	$\hat W'$ is continuously invertible,
	which we will do using \cref{lem:saddle_point_matrix}.
	We notice that the operator $\hat W'$ has the required
	block structure if we set
	\begin{equation*}
		\begin{aligned}
			\hat A &:=
			\begin{bmatrix}
				h''_{\beta\beta}(\beta,y) & h''_{\beta y}(\beta,y) & 0
				\\
				h''_{y\beta}(\beta,y) & h''_{yy}(\beta,y) & 0
				\\
				0 & 0 & \hat\sigma I
			\end{bmatrix},
			\\
			\hat A &\colon\R^n\times H_0^1(\Omega)\times L^2(\AA_3(p))
			\to \R^n\times H^{-1}(\Omega)\times L^2(\AA_3(p)),
			\\
			\hat B &:=
			\begin{bmatrix}
				I_1^\top K_T & 0 & 0 
				\\
				0 & A & -B I_3 
			\end{bmatrix},
			\\
			\hat B &\colon\R^n\times H_0^1(\Omega)\times L^2(\AA_3(p))
			\to \R^{l_1}\times H^{-1}(\Omega).
		\end{aligned}
	\end{equation*}
	Since $A$ is invertible and $I_1^\top K_T$ is surjective
	by assumption, it follows that $\hat B$ is surjective.
	In order to show that $\hat W'$ is continuously invertible,
	it remains to show that $\hat A$ is coercive on $\ker \hat B$.

	Let $(\hat\beta,\hat y,\hat u)\in\ker \hat B$ be given.
	Then
	\begin{equation*}
		\norm{(\hat y,\hat u)}_{H_0^1(\Omega) \times L^2(\AA_3(p))}
		=
		\norm{(A^{-1}BI_3\hat u,\hat u)}_{H_0^1(\Omega) \times L^2(\AA_3(p))}
		\leq (1+\norm{A^{-1}B})\norm{\hat u}_{L^2(\AA_3(p))}
	\end{equation*}
	holds.
  Recall from \eqref{eq:h_definition} 
  that 
	$h(\beta,y)=\hat h(\beta,y)
	+\frac{\sigma_\beta}{2}\sum_{i=1}^n \left(\frac1{\beta_i}\right)^2$
  and that $\hat h$ is convex, 
  and that for
  $\frac{\sigma_\beta}{2}\sum_{i=1}^n \left(\frac1{\beta_i}\right)^2$
  we can directly calculate the second derivative, which is a diagonal matrix
  with strictly positive entries, if
	$\beta_i > 0$.
	Therefore, there exists a constant $\sigma_r>0$ for which
  \begin{equation}
    \left( h''(\beta,y)(\hat\beta,\hat y)\right)(\hat\beta,\hat y)
    \geq 
    \sigma_r \hat\beta^\top \hat\beta
  \end{equation}
  holds, where $\sigma_r$ depends on the upper bound of $\beta_i$.
	This implies
	\begin{equation*}
		\begin{aligned}
			\dual{ \hat A (\hat\beta,\hat y,\hat u)}{(\hat\beta,\hat y,\hat u)}
			&\geq \sigma_r\norm{\hat\beta}^2_{\R^n}+\hat\sigma\norm{\hat u}_{L^2(\AA_3(p))}^2
			\\ &\geq
			\sigma_r\norm{\hat\beta}^2_{\R^n} 
			+ \hat\sigma(1+\norm{A^{-1}B})^{-2} \norm{(\hat y,\hat u)}^2_{H_0^1(\Omega) \times L^2(\AA_3(p))}
			\\ &\geq
			\hat\gamma \norm{(\hat\beta,\hat y,\hat u)}^2_{\R^n \times H_0^1(\Omega) \times L^2(\AA_3(p))},
		\end{aligned}
	\end{equation*}
	where $\hat\gamma>0$ is a suitable constant.
	Thus $\hat A$ is coercive on $\ker \hat B$.
	It follows from \cref{lem:saddle_point_matrix} that
	$\hat W'$ is continuously invertible.
	Because $\hat B$ is surjective, there exists a constant
	$\alpha>0$ such that $B^1(0)\subset \hat B(B^\alpha(0))$.
	Since there are only finitely many possibilities for $I_1$
	and $I_3$ is not needed for surjectivity, the constant $\alpha$
	can be chosen such that it is independent of $I_1$ and $I_3$.
	For $\norm{\hat A}$ we note that it can be bounded
	by a constant which can depend on an upper bound on $\norm{y}_{H_0^1(\Omega)}$
	and a lower bound on $\beta_i$.

	It follows from \cref{lem:saddle_point_matrix}
	that the estimate
	$\norm{\hat W'^{-1}}\leq 4c^5$ holds for a suitable constant 
	$ c>0$
	which does not depend on $\beta,y,u,z,p$ but can depend
	on an upper bound of $\norm{y}_{H_0^1(\Omega)}$, the lower bound of $\beta_i$ and
  on $K_T,A,B,h,\hat\sigma,\sigma_r$.

	Next, we combine this result with \cref{lem:symmetrization}
	to show the invertibility of $W'(\beta,y,u,z,p)$.
	Let $(\beta_2,y_2,u_2,z_2,p_2)$ be a right-hand side 
	as in \eqref{eq:W'_byuzp}.
	Since $\hat W'$ is invertible, by \cref{lem:symmetrization} there
	exists a unique solution $(\beta_1,y_1,u_1,z_1,p_1)$
	of \eqref{eq:W'_byuzp}.
	Using the estimate $\norm{\hat W'^{-1}}\leq 4c^5$
	and \eqref{eq:hatW'_byuzp},
	one get an estimate of the form
	$\norm{(\beta_1,y_1,u_1,z_1,p_1)}\leq 
	C\norm{(\beta_2,y_2,u_2,z_2,p_2)}$,
	where $C>0$ is a suitable constant that can depend on
  $c,\hat\sigma, K_T, B$, the upper bound of $\norm{y}_{H_0^1(\Omega)}$ and the bounds of
  $\beta$. The constant $C$ however, does not depend
	on 
  $(\beta,y,u,z,p)$ or any of the embedding operators $I_1,I_2,I_3,I_4$.
  .
	Since we can estimate the norm of the unique solution in \eqref{eq:W'_byuzp}
	by the norm of the right-hand side, the claimed
	invertibility and estimate $\norm{W'(\beta,y,u,z,p)^{-1}}\leq C$
	follow.
\end{proof}

\begin{lemma}
	\label{lem:invertible_in_neighborhood}
	Let
	$(\beta,y,u,z,p)\in \Bad\times H_0^1(\Omega)\times L^2(\Omega)
\times \R^{n+1}\times H_0^1(\Omega)$ 
	be a point such that $W(\beta,y,u,z,p)=0$.
	Then the Newton derivative $W'$ is uniformly continuously
	invertible in a neighborhood of $(\beta,y,u,z,p)$.
\end{lemma}
\begin{proof}
	We want to apply \cref{lem:newton_derivative_invertible}.
	We need to verify that
	$I_1^\top K_T$ (which can depend on $\beta$ and $z$) 
	is surjective in a neighborhood.

	From the definition of $W$, we get 
	$z\geq0$,
	$K_T\beta-b_T\leq0$ and
	$z^\top (K_T\beta - b_T) =0$.
	In particular, $\beta \in T$.
	Recall that $T$ is a non-degenerate simplex.
	Thus, at most $n$ constraints in the system $K_T\beta\leq b_T$
	are active, and these active constraints are linearly independent.
	Furthermore,
	if $i\in\set{1,\ldots,n+1}$ is an index of an inactive constraint,
	we have $z_i=0$ due to the complementarity condition,
	and therefore $i\in\AA_2(\beta,z)$ and $i\not\in\AA_1(\beta,z)$.
	Thus, $\AA_1(\beta,z)$ contains at most $n$ elements.
	Therefore, the rows of $K_T$ which correspond to the index set $\AA_1(\beta,z)$
	are linearly independent, which yields that $I_1^\top K_T$ is surjective
	for this particular $\beta$, $z$.

	If $i\in \AA_2(\beta,z)$, then $i\in\AA_2(\hat\beta,\hat z)$
	holds also for $(\hat\beta,\hat z)$ that are sufficiently close to $(\beta,z)$.
	Thus, $\AA_1(\beta,z)$ cannot get larger in a neighborhood of $(\beta,z)$.
	Hence, the rows of $K_T$ that correspond to $\AA_1(\beta,z)$
	stay linearly independent in a neighborhood, i.e.\ $I_1^\top K_T$
	is surjective in a neighborhood of $(\beta,z)$.

	Now to apply \cref{lem:newton_derivative_invertible}
  we restrict the neighborhood such that $\beta > \frac 1{2a_i}$ if necessary.
  This guarantees the lower bound
  $\beta>\frac 1{2a_i}$. The upper bound of $\norm{y}_{H_0^1(\Omega)}$ is obtained from the
  coercivity of $f$ with constant $\gamma_{k,T}\mu$ (cf.
    \cref{assumptions}\ref{asmp:f_strongly_convex}.
    Hence, with \cref{lem:newton_derivative_invertible}
    there exists a constant $C>0$, such that
 $\norm{W'(\beta,y,u,z,p)^{-1}}\leq C$ 
 in the considered neighborhood of $(\beta,y,u,z,p)$.
\end{proof}

Now we are ready to give our final \lcnamecref{thm:fast_convergence},
which states that \cref{alg:semismooth_Newton_OVRP}
converges superlinearly.

\begin{theorem}
	\label{thm:fast_convergence}
	Let the function $W$ be given as in 
	\eqref{W_definition}.
	Further,
	we denote by
	$(\beta_{k,T},y_{k,T},u_{k,T})$
	the unique global solution of
	\eqref{eq:OVR_affine_penalty}
	and by $z_{k,T}$, $p_{k,T}$ the corresponding multipliers that
	satisfy \eqref{eq:KKT_OVRP_T}.
	Then there exists a neighborhood
	of the point 
	$(\beta_{k,T},y_{k,T},u_{k,T},z_{k,T},p_{k,T})$
	such that for all 
	initial values
	$(\beta_0,y_0,u_0,p_0,z_0)$
	from this neighborhood,
		the semismooth Newton method from 
	\cref{alg:semismooth_Newton_OVRP}
	either terminates in the $i$-th step with 
	$(\beta_i,y_i,u_i,z_i,p_i) = (\beta_{k,T},y_{k,T},u_{k,T},z_{k,T},p_{k,T})$
	or generates a sequence that converges $q$-superlinearly to 
	$(\beta_{k,T},y_{k,T},u_{k,T},z_{k,T},p_{k,T})$ in
  $\R^n\times H^1_0(\Omega)\times L^2(\Omega)\times \R^{n+1}\times H^1_0(\Omega)$.
\end{theorem}
\begin{proof}
	We already established that the function $W$ is semismooth 
	in the solution to 
	\eqref{eq:OVR_affine_penalty} 
	(see \cref{lem:W_semismooth}).
	We have proven in 
	\cref{lem:invertible_in_neighborhood}
	that the derivative from
	\cref{lem:W_semismooth} is invertible and the norm of the inverse is bounded
	on a neighborhood of a solution. 
	The result is now a direct application of
	\cite[Theorem 3.13]{Ulbrich2011}.
	In particular, we do not need a smoothing step,
	since the spaces in which $W$ is Newton differentiable
	coincide with the spaces in which the Newton derivative is uniformly invertible,
	see \cref{lem:W_semismooth,lem:newton_derivative_invertible}.
\end{proof}

\section{Numerical experiments}
\label{sec:numerical_experiments}

In this section we present an example for
\cref{alg:global_solution_of_IOC_penalty}
to illustrate the 
convergence behavior towards a global minimizer.
To this end, we consider the parameter identification problem
\begin{equation}
	\begin{aligned}
		 \min_{\beta,y,u} \qquad & \frac12\norm{y-y_m}^2_{L^2(\Omega)} 
		+ \frac{\sigma_u}{2}\norm{u-u_m}^2_{L^2(\Omega)} 
		+\frac{\sigma_\beta}{2}\norm{\beta - \beta_m}^2_{\R^n}
		\eqqcolon F_1\byu\\\
		\text{s.t.} \qquad & \beta\in\Bad,\\
		&(y,u)\in\Psi(\beta),
	\end{aligned}
	\label{eq:upper_level_special_problem_numerics_section}
\end{equation}
where $\Psi:\R^2\rightarrow H^1_0(\Omega) \times L^2(\Omega)$ 
denotes the solution mapping of the parameter $\beta$ to the unique
solution of the lower-level problem
\begin{equation}
\begin{aligned}
  \min_{y,u} \qquad &\frac{1}{2\beta_1}\norm{y-y_{d,1}}^2_{L^2(\Omega)} 
	+ \frac{1}{2\beta_2}\norm{y-y_{d,2}}^2_{L^2(\Omega)} 
  + \frac{\sigma_l}{2}\norm{u}^2_{L^2(\Omega)} \eqqcolon f\byu
  \\
  \text{s.t.} \qquad &  0 = -\Delta y - u\qquad \text{in}\ \Omega,\\
  &0 = y \phantom{-\Delta - u}\qquad\text{on}\ \partial\Omega,\\
  &u  \in \Uad.
\end{aligned}
\end{equation}
Let us define the data present in this bilevel optimization problem.
We use the sets $\Bad \coloneqq [0.1, 1]^2$ and $\Omega = ( -1,1 )^2$
and the two possible desired states
\begin{align*}
	y_{d,1}&:\Omega\rightarrow \R,\qquad y_{d,1}(x) = \sin(\pi x_1)\sin(\pi x_2),\\
	y_{d,2}&:\Omega\rightarrow \R,\qquad y_{d,2}(x)
	= (x_1+1)(x_1-1)(x_2+1)(x_2-1).
\end{align*}
The regularization parameter for the lower level is $\sigma_l = 0.03$.
Additionally, we introduce box constraints for the control via
\begin{align*}
	u\in\Uad&\coloneqq\{u\in L^2(\Omega)\ |\ u_a\leq u\leq u_b \text{ a.e.\ on } \Omega \},\\
u_a(x) &\coloneqq 0,
\qquad
u_b(x) \coloneqq 3.
\end{align*}
It turns out that
these constraints are active on parts of the domain
for the choice of the parameter $\beta = (0.6,0.3)^\top$. For the upper level we fix the parameters 
$\sigma_u = 0.05$ and $\sigma_\beta = 10^{-5}$.
We also choose $\beta_m\coloneqq(0.6,0.3)^\top$ and
$(y_m,u_m)\coloneqq \Psi((0.6,0.3)^\top)$,
i.e.\ the objective value of $F_1$ is zero for the solution
to the lower-level problem with $\beta=\beta_m$.
We call this setting ``fully reachable target state''. 
We mention that when this setting is 
implemented, the functions $y_m$, $u_m$ are not the 
analytical solutions,
but are calculated directly using the 
finite element solutions for the lower level.

For the setting of this \lcnamecref{sec:numerical_experiments}, 
\cref{assumptions} is valid. 
Additionally, for the
chosen functionals and parameters we can apply the semismooth Newton method from \cref{sec:semismooth_newton}
to solve the subproblems \eqref{eq:OVR_affine_penalty}.
In order to illustrate some fundamental properties of
the proposed solution algorithm, 
we consider two additional problems 
that only differ in the choice of the objective functional,
i.e.\ the functions
\begin{align*}
	F_2\byu&\coloneqq 
	\frac12\norm{y-y_m}^2_{L^2(\Omega)} 
	+ \frac{\sigma_u}{2}\norm{u-u_m}^2_{L^2(\Omega)} 
	+\frac{\sigma_\beta}{2}\norm{\beta}^2_{\R^n},\\
	F_3\byu&\coloneqq 
	\frac12\norm{y-\hat y_m}^2_{L^2(\Omega)} 
	+ \frac{\sigma_u}{2}\norm{u-\hat u_m}^2_{L^2(\Omega)} 
	+\frac{\sigma_\beta}{2}\sum_{i=1}^{2}\frac{1}{\beta_i^2}
\end{align*}
are used instead of $F_1$.
In the second objective functional $F_2$, 
the $\beta$ term is only introduced as a regularization. 
This will be called ``reachable target state''.
The functional $F_3$ is set up with desired states $\hat y_m$ and $\hat u_m$ 
that are given by
\begin{align*}
  \hat y_m&:\Omega\rightarrow \R,\qquad \hat y_m(x) = (x_1-1)(x_1+1)\sin(\pi x_2),\\
	\hat u_m&:\Omega\rightarrow \R,\qquad \hat u_m(x) = 0.
\end{align*}
This state and control have the property that
they do not arise as a solution of the lower-level problem.
This setting is named ``unreachable target state''.
We expect a noticeable difference in the convergence speed 
for the introduced settings, see \cref{rem:better_small_T}.

The refinement of the subdivision 
will be implemented by splitting the triangles 
at the midpoint of the edges. 
This refinement procedure is the application of
\cref{lem:subdivision_by_hypercube}
to the two-dimensional case. 
However, in this special case we can even guarantee that the 
diameter of the simplices is halved in each refinement.
We initialize \cref{alg:global_solution_of_IOC_penalty} with the domain 
$\Bad $ 
split into two triangles.

We use an implementation 
with the suggested improvements mentioned at the end of 
\cref{sec:algo}.
In each iteration we get a lower bound on the optimal objective value from the
element with the lowest objective value for the solution to
\eqref{eq:OVR_affine_penalty}. We obtain an upper bound from the vertex with
the lowest objective value. Hence every element whose relaxed optimal objective value
is above the upper bound can be dismissed, since
the relaxed optimal objective value is smaller than or equal to the objective value of the
original subproblem.
Further, in each iteration we refine the best $15\%$ of the active triangles
with respect to the objective value for the solution to
\eqref{eq:OVR_affine_penalty}. This is done to effectively utilize
parallelization. 
Additionally, we refine the worst $5\%$ as a measure to ``clean
up old triangles''. 
Otherwise, for some triangles that are quite far from the
actual solution but for which (by chance) the objective value comes really close, the
algorithm might take a long time to refine this element. 
Lastly, the algorithm
runs until a set amount of elements ($3\cdot 10^5$) is reached or the
difference between lower and upper bound is sufficiently mall.
For the setting of $F_1$ we chose a target bound difference of $10^{-13}$, for $F_2$ we
chose a target bound difference of $10^{-11}$. In the case of the ``unreachable target
state'' ($F_3$) the element limit was reached.
\begin{figure}[ht]
 \centering
 {
%
%
\definecolor{mycolor1}{rgb}{0.00000,0.44700,0.74100}%
\definecolor{mycolor2}{rgb}{0.85000,0.32500,0.09800}%
\begin{tikzpicture}[every node/.style={scale=0.7}]
\begin{axis}[%
width=3.2cm,
at={(0in,0in)},
scale only axis,
xmode=log,
xmin=1,
xmax=10000,
xminorticks=true,
xlabel style={font=\color{white!15!black}},
xlabel={number of subproblems},
ymode=log,
ymin=1e-35,
ymax=1,
yminorticks=true,
ylabel style={font=\color{white!15!black}},
ylabel={bound value},
axis background/.style={fill=white},
title style={font=\bfseries,font = \fontsize{14}{14}\selectfont,thick},
title={\textbf{Bounds}; $F_1$}
]
\addplot [color=mycolor1, forget plot]
  table[row sep=crcr]{%
2	0.0350894924584624\\
28	0.0350894924584624\\
28	0.00208996610212184\\
32	0.00208996610212184\\
32	0.000325247952369774\\
80	0.000325247952369774\\
176	0.000325247952369774\\
176	3.23864731455739e-05\\
272	3.23864731455739e-05\\
272	6.06504454885879e-06\\
433	6.06504454885879e-06\\
433	5.72089545423244e-06\\
776.000000000001	5.72089545423244e-06\\
776.000000000001	5.09667810764545e-07\\
1189	5.09667810764545e-07\\
1189	9.24787604052593e-08\\
1798	9.24787604052593e-08\\
1798	9.08409921223813e-08\\
2204	9.08409921223813e-08\\
2204	7.9579267702377e-09\\
2652	7.9579267702377e-09\\
2652	1.44937086525734e-09\\
2835	1.44937086525734e-09\\
2835	1.41652203546912e-09\\
3225	1.41652203546912e-09\\
3225	1.24353583959111e-10\\
3818	1.24353583959111e-10\\
3818	2.26378264376735e-11\\
4598.99999999999	2.26378264376735e-11\\
4598.99999999999	2.21387512523615e-11\\
5488.00000000001	2.21387512523615e-11\\
5488.00000000001	1.9430033021282e-12\\
6314.00000000001	1.9430033021282e-12\\
};
\addplot [color=mycolor2, forget plot]
  table[row sep=crcr]{%
2	1.96434711649957e-32\\
28	1.96434711649957e-32\\
28	3.54544339471853e-32\\
32	3.54544339471853e-32\\
80	3.54544339471853e-32\\
80	8.93829712308899e-33\\
175.999999999999	8.93829712308899e-33\\
175.999999999999	2.63496824487633e-32\\
272.000000000001	2.63496824487633e-32\\
272.000000000001	2.1097617646141e-32\\
433	2.1097617646141e-32\\
433	2.62513979573734e-32\\
775.999999999995	2.62513979573734e-32\\
775.999999999995	7.03898370887484e-32\\
1189	7.03898370887484e-32\\
1798.00000000001	7.03898370887484e-32\\
1798.00000000001	2.78690980195287e-32\\
2204.00000000001	2.78690980195287e-32\\
2204.00000000001	3.93098611810891e-32\\
2651.99999999999	3.93098611810891e-32\\
2651.99999999999	3.518574626236e-32\\
2835	3.518574626236e-32\\
2835	3.89316606455395e-32\\
3225.00000000002	3.89316606455395e-32\\
3225.00000000002	5.35362537092164e-32\\
3818.00000000002	5.35362537092164e-32\\
4598.99999999999	5.35517637994158e-32\\
4598.99999999999	3.68426409136227e-32\\
5488.00000000001	3.68426409136227e-32\\
5488.00000000001	6.53941414460493e-32\\
6313.99999999998	6.53941414460493e-32\\
};
\end{axis}
\end{tikzpicture}%
}\hspace{.5cm}%
 {
%
%
\definecolor{mycolor1}{rgb}{0.00000,0.44700,0.74100}%
\definecolor{mycolor2}{rgb}{0.85000,0.32500,0.09800}%
\begin{tikzpicture}[every node/.style={scale=0.7}]
\begin{axis}[%
width=3.2cm,
at={(0in,0in)},
scale only axis,
xmode=log,
xmin=10,
xmax=17026,
xminorticks=true,
xlabel style={font=\color{white!15!black}},
xlabel={number of subproblems},
ymode=log,
ymin=1e-07,
ymax=0.1,
yminorticks=true,
ylabel style={font=\color{white!15!black}},
ylabel={bound value},
axis background/.style={fill=white},
title style={font=\bfseries,font = \fontsize{14}{14}\selectfont,thick},
title={\textbf{Bounds}; $F_2$}
]
\addplot [color=mycolor1, forget plot]
  table[row sep=crcr]{%
24	0.0350921924584624\\
93.9999999999999	0.0350921924584624\\
164	0.0350921924584624\\
164	0.0318581335142266\\
234	0.0318581335142266\\
418	0.0318581335142266\\
418	0.0296740877059077\\
694.000000000001	0.0296740877059077\\
694.000000000001	0.0133977631707967\\
1070	0.0133977631707967\\
1070	0.0123339981773896\\
1542	0.0123339981773896\\
1542	0.00209199110212184\\
2084	0.00209199110212184\\
2368	0.00209199110212184\\
2744	0.00209199110212184\\
3216	0.00209199110212184\\
3216	0.000405850697168083\\
3688	0.000405850697168083\\
4063	0.000405850697168083\\
4479	0.000405850697168083\\
4893	0.000405850697168083\\
4893	7.93251618582745e-05\\
5456	7.93251618582745e-05\\
5915.00000000001	7.93251618582745e-05\\
5915.00000000001	5.29090391305726e-05\\
6364	5.29090391305726e-05\\
6364	4.22995180814405e-05\\
7015	4.22995180814405e-05\\
7015	3.46645981455739e-05\\
7654	3.46645981455739e-05\\
7654	8.25879454885878e-06\\
8293.00000000001	8.25879454885878e-06\\
8834	8.25879454885878e-06\\
8834	6.00514223392541e-06\\
9407.00000000001	6.00514223392541e-06\\
9407.00000000001	2.99369267359992e-06\\
9889.00000000001	2.99369267359992e-06\\
9889.00000000001	2.75615218576455e-06\\
10256	2.75615218576455e-06\\
10742	2.75615218576455e-06\\
10742	2.33468864837238e-06\\
11279	2.33468864837238e-06\\
11653	2.33468864837238e-06\\
11929	2.33468864837238e-06\\
11929	2.30515249837449e-06\\
12155	2.30515249837449e-06\\
12155	2.26581223364283e-06\\
12266	2.26581223364283e-06\\
12337	2.26581223364283e-06\\
12397	2.26581223364283e-06\\
12647	2.26581223364283e-06\\
12870	2.26581223364283e-06\\
13144	2.26581223364283e-06\\
13322	2.26581223364283e-06\\
13322	2.25769899241855e-06\\
13547	2.25769899241855e-06\\
13777	2.25769899241855e-06\\
13928	2.25769899241855e-06\\
14054	2.25769899241855e-06\\
14186	2.25759122356223e-06\\
14288	2.25759122356223e-06\\
14426	2.25759122356223e-06\\
14560	2.25759122356223e-06\\
14658	2.25759122356223e-06\\
14770	2.25759122356223e-06\\
14878	2.25759122356223e-06\\
15018	2.25759122356223e-06\\
15116	2.25759122356223e-06\\
15240	2.25759122356223e-06\\
15366	2.25759122356223e-06\\
15550	2.25759122356223e-06\\
15704	2.25759122356223e-06\\
15852	2.25759122356223e-06\\
16010	2.25759122356223e-06\\
16200	2.25759122356223e-06\\
16364	2.25759122356223e-06\\
16540	2.25759122356223e-06\\
16674	2.25759122356223e-06\\
16836	2.25759122356223e-06\\
17026	2.25759122356223e-06\\
};
\addplot [color=mycolor2, forget plot]
  table[row sep=crcr]{%
24	1e-07\\
93.9999999999999	1e-07\\
164	1e-07\\
234	1e-07\\
418	1e-07\\
694.000000000001	1e-07\\
694.000000000001	2.7578125e-07\\
1070	2.7578125e-07\\
1542	2.7578125e-07\\
1542	4.83203125e-07\\
2084	4.83203125e-07\\
2084	5.78125e-07\\
2368	5.78125e-07\\
2368	6.501953125e-07\\
2744	6.501953125e-07\\
2744	1.05625e-06\\
3216	1.05625e-06\\
3688	1.05625e-06\\
3688	1.087890625e-06\\
4063	1.087890625e-06\\
4063	1.4447265625e-06\\
4479	1.4447265625e-06\\
4479	1.580078125e-06\\
4893	1.580078125e-06\\
4893	1.7470703125e-06\\
5456	1.7470703125e-06\\
5456	1.80244140625e-06\\
5915.00000000001	1.80244140625e-06\\
5915.00000000001	1.91241455078125e-06\\
6364	1.91241455078125e-06\\
6364	2.040625e-06\\
7015	2.040625e-06\\
7654	2.040625e-06\\
7654	2.10835202510231e-06\\
8293.00000000001	2.10835202510231e-06\\
8293.00000000001	2.11181640625e-06\\
8834	2.11181640625e-06\\
9407.00000000001	2.11181640625e-06\\
9407.00000000001	2.17457389720292e-06\\
9889.00000000001	2.17457389720292e-06\\
9889.00000000001	2.19410400390625e-06\\
10256	2.19410400390625e-06\\
10256	2.22078537773083e-06\\
10742	2.22078537773083e-06\\
10742	2.22662195730486e-06\\
11279	2.22662195730486e-06\\
11279	2.23077496798576e-06\\
11653	2.23077496798576e-06\\
11653	2.23871727136106e-06\\
11929	2.23871727136106e-06\\
12266	2.24253665208817e-06\\
12337	2.24253665208817e-06\\
12397	2.24253665208817e-06\\
12647	2.24253665208817e-06\\
12870	2.24253665208817e-06\\
13144	2.24253665208817e-06\\
15240	2.24319403469563e-06\\
15366	2.24319403469563e-06\\
15550	2.24319403469563e-06\\
15852	2.24327624146827e-06\\
16010	2.24327624146827e-06\\
16200	2.24327624146827e-06\\
16836	2.24335845578462e-06\\
17026	2.24335845578462e-06\\
};
\end{axis}
\end{tikzpicture}%
}\hspace{.5cm}%
 {
%
%
\definecolor{mycolor1}{rgb}{0.00000,0.44700,0.74100}%
\definecolor{mycolor2}{rgb}{0.85000,0.32500,0.09800}%
\begin{tikzpicture}[every node/.style={scale=0.7}]

\begin{axis}[%
width=3.2cm,
at={(0in,0in)},
scale only axis,
xmode=log,
xmin=15.1999402459192,
xmax=64695.5096898424,
xminorticks=true,
xlabel style={font=\color{white!15!black}},
xlabel={number of subproblems},
ymode=log,
ymin=0.47792585846082,
ymax=0.71534762051392,
ytick={0.5,0.55,0.6,0.65,0.7},
yticklabels={{0.5},{0.52},{0.54},{0.56},{0.58},{0.6},{0.62},{0.64},{0.66},{0.68},{0.7}},
yminorticks=true,
ylabel style={font=\color{white!15!black}},
ylabel={bound value},
axis background/.style={fill=white},
title style={font=\bfseries,font = \fontsize{14}{14}\selectfont,thick},
title={\textbf{Bounds}; $F_3$}
]
\addplot [color=mycolor1, forget plot]
  table[row sep=crcr]{%
24	0.705001674798699\\
116	0.705001674798699\\
208	0.705001674798699\\
278	0.705001674798699\\
462.000000000001	0.705001674798699\\
646	0.705001674798699\\
933.000000000001	0.705001674798699\\
1122	0.705001674798699\\
1410	0.705001674798699\\
1564	0.705001674798699\\
1762	0.705001674798699\\
2044	0.705001674798699\\
2350	0.705001674798699\\
2773	0.705001674798699\\
2773	0.704686085584673\\
3295	0.704686085584673\\
3955	0.704686085584673\\
4684	0.704686085584673\\
5439	0.704686085584673\\
6257	0.704686085584673\\
7077	0.704676999316602\\
7868	0.704676999316602\\
8544	0.704669122482258\\
9300.00000000001	0.704669122482258\\
10029	0.704669122482258\\
10747	0.704669122482258\\
11508	0.704669122482258\\
12319	0.704669122482258\\
13146	0.704669122482258\\
14062	0.704668815592625\\
15092	0.704668815592625\\
16304	0.704668815592625\\
17730	0.704668776340371\\
19284	0.704668776340371\\
21062	0.704668776340371\\
23150	0.704668747955993\\
25434	0.704668747955993\\
27946	0.704668747955993\\
30852	0.704668747955993\\
34072	0.704668747955993\\
37414	0.704668747955993\\
41516	0.704668746037356\\
45790	0.704668746037356\\
50136.0000000001	0.704668746037356\\
55218	0.704668746037356\\
60614.0000000001	0.704668746037356\\
65784	0.704668746037356\\
};
\addplot [color=mycolor2, forget plot]
  table[row sep=crcr]{%
24	0.482783331878682\\
116	0.482783331878682\\
278	0.482806389729922\\
462.000000000001	0.482806389729922\\
1122	0.482868006434895\\
1410	0.482868006434895\\
1410	0.496820299197651\\
1564	0.496820299197651\\
1564	0.542411100835221\\
1762	0.542411100835221\\
1762	0.559175154763825\\
2044	0.559175154763825\\
2044	0.587393069943672\\
2350	0.587393069943672\\
2350	0.613168601292877\\
2773	0.613168601292877\\
2773	0.63586850283422\\
3295	0.63586850283422\\
3295	0.648072628740488\\
3955	0.648072628740488\\
3955	0.659290004736015\\
4684	0.659290004736015\\
4684	0.670166394657857\\
5439	0.670166394657857\\
5439	0.67881498941947\\
6257	0.67881498941947\\
6257	0.683354964769943\\
7077	0.683354964769943\\
7077	0.688601243870168\\
7868	0.688601243870168\\
7868	0.693495676012781\\
8544	0.693495676012781\\
8544	0.695769835057999\\
9300.00000000001	0.695769835057999\\
9300.00000000001	0.697393106408869\\
10029	0.697393106408869\\
10029	0.698727293471714\\
10747	0.698727293471714\\
10747	0.699830516286559\\
11508	0.699830516286559\\
11508	0.700515153549152\\
12319	0.700515153549152\\
12319	0.701211096012878\\
13146	0.701211096012878\\
13146	0.701903463939289\\
14062	0.701903463939289\\
14062	0.702399313872278\\
15092	0.702399313872278\\
15092	0.702784358406693\\
16304	0.702784358406693\\
16304	0.703156622344231\\
17730	0.703156622344231\\
17730	0.703552201598273\\
19284	0.703552201598273\\
19284	0.703723896701676\\
21062	0.703723896701676\\
21062	0.703940836197266\\
23150	0.703940836197266\\
23150	0.704138384204676\\
25434	0.704138384204676\\
25434	0.704226003242586\\
27946	0.704226003242586\\
27946	0.704325795885088\\
30852	0.704325795885088\\
30852	0.704407292411842\\
34072	0.704407292411842\\
50136.0000000001	0.704593324328571\\
55218	0.704593324328571\\
65784	0.704617463695351\\
};
\end{axis}
\end{tikzpicture}%
}%
 \caption{Upper bound (blue) and lower bound (red) for the setting of $F_1$, $F_2$ and
$F_3$ w.r.t. the number of solved subproblems.}
  \label{fig:lower_upper_bounds}
\end{figure}
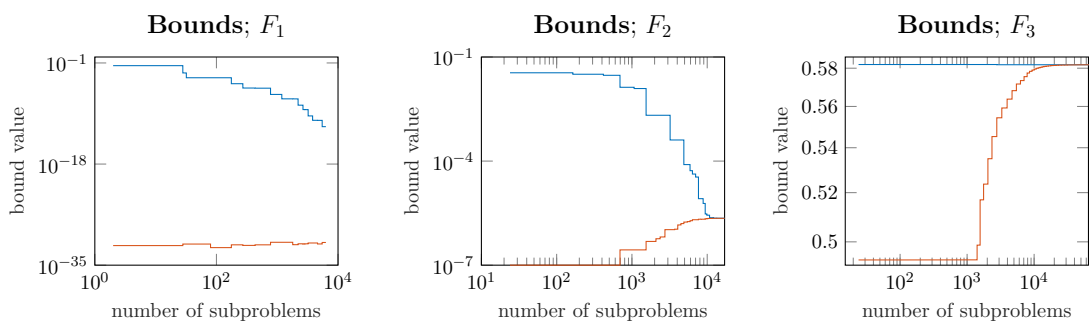
\begin{figure}[ht]
 \centering
 {
%
%
\definecolor{mycolor1}{rgb}{0.00000,0.44700,0.74100}%
\definecolor{mycolor2}{rgb}{0.85000,0.32500,0.09800}%
\begin{tikzpicture}[every node/.style={scale=0.7}]
\begin{axis}[%
width=3.2cm,
at={(0in,0in)},
scale only axis,
xmode=log,
xmin=1,
xmax=100000,
xminorticks=true,
xlabel style={font=\color{white!15!black}},
xlabel={number of subproblems},
ymode=log,
ymin=1e-06,
ymax=1,
yminorticks=true,
ylabel style={font=\color{white!15!black}},
ylabel={distance to $\bar\beta$},
axis background/.style={fill=white},
title style={font=\bfseries,font = \fontsize{14}{14}\selectfont,thick},
title={\textbf{Convergence to}  $\bar\beta$; $F_1$}
]
\addplot [color=mycolor1, forget plot]
  table[row sep=crcr]{%
2	0.254950975679639\\
9.00000000000001	0.254950975679639\\
9.00000000000001	0.0559016994374947\\
17	0.0559016994374947\\
17	0.0673145600891814\\
23	0.0673145600891814\\
35	0.0673145600891814\\
35	0.00698771242968691\\
52	0.00698771242968691\\
52	0.00841432001114755\\
124	0.00841432001114755\\
124	0.00737029775941929\\
305	0.00737029775941929\\
305	0.000873464053710731\\
425	0.000873464053710731\\
425	0.00105179000139357\\
572.000000000001	0.00105179000139357\\
572.000000000001	0.000921287219927242\\
862	0.000921287219927242\\
862	0.000109183006713956\\
1165	0.000109183006713956\\
1165	0.000131473750174063\\
1829	0.000131473750174063\\
1829	0.000115160902491068\\
3353	0.000115160902491068\\
3353	1.36478758392073e-05\\
4736	1.36478758392073e-05\\
4736	1.64342187719047e-05\\
5971	1.64342187719047e-05\\
5971	1.43951128112798e-05\\
7515.00000000001	1.43951128112798e-05\\
7515.00000000001	1.70598447995366e-06\\
8965.00000000001	1.70598447995366e-06\\
10772	1.70598447995366e-06\\
12808	1.70598447995366e-06\\
14604	1.70598447995366e-06\\
15999	1.70598447995366e-06\\
17065	1.70598447995366e-06\\
18008	1.70598447995366e-06\\
18822	1.70598447995366e-06\\
19426	1.70598447995366e-06\\
19897	1.70598447995366e-06\\
20242	1.70598447995366e-06\\
20440	1.70598447995366e-06\\
20565	1.70598447995366e-06\\
20707	1.70598447995366e-06\\
20843	1.70598447995366e-06\\
20913	1.70598447995366e-06\\
20984	1.70598447995366e-06\\
21112	1.70598447995366e-06\\
21249	1.70598447995366e-06\\
21384	1.70598447995366e-06\\
21452	1.70598447995366e-06\\
};
\addplot [color=mycolor2, forget plot]
  table[row sep=crcr]{%
2	0.860232526704263\\
9.00000000000001	0.860232526704263\\
17	0.860232526704263\\
23	0.860232526704263\\
35	0.860232526704263\\
52	0.860232526704263\\
124	0.860232526704263\\
305	0.860232526704263\\
425	0.860232526704263\\
425	0.559016994374947\\
572.000000000001	0.559016994374947\\
862	0.559016994374947\\
862	0.254950975679639\\
1165	0.254950975679639\\
1165	0.107529065838033\\
1829	0.107529065838033\\
1829	0.0332192056647958\\
3353	0.0332192056647958\\
3353	0.0261922728347503\\
4736	0.0261922728347503\\
4736	0.00737029775941929\\
5971	0.00737029775941929\\
5971	0.00305712418798804\\
7515.00000000001	0.00305712418798804\\
8965.00000000001	0.00305712418798804\\
8965.00000000001	0.000360469312994121\\
10772	0.000360469312994121\\
12808	0.000360469312994121\\
14604	0.000360469312994121\\
14604	0.000337905530656076\\
15999	0.000337905530656076\\
17065	0.000337905530656076\\
18008	0.000337905530656076\\
18008	0.000323198298214853\\
18822	0.000323198298214853\\
19426	0.000323198298214853\\
19426	0.000286930158045663\\
19897	0.000286930158045663\\
20242	0.000286930158045663\\
20440	0.000286930158045663\\
20565	0.000286930158045663\\
20707	0.000286930158045663\\
20843	0.000286930158045663\\
20913	0.000286930158045663\\
20984	0.000286930158045663\\
21112	0.000286930158045663\\
21112	0.000177107154889358\\
21249	0.000177107154889358\\
21249	0.000154075995622228\\
21384	0.000154075995622228\\
21452	0.000154075995622228\\
};
\end{axis}
\end{tikzpicture}%
}\hspace{.2cm}%
 {
%
%
\definecolor{mycolor1}{rgb}{0.00000,0.44700,0.74100}%
\definecolor{mycolor2}{rgb}{0.85000,0.32500,0.09800}%
\begin{tikzpicture}[every node/.style={scale=0.7}]

\begin{axis}[%
width=3.2cm,
at={(0in,0in)},
scale only axis,
xmode=log,
xmin=10,
xmax=100000,
xminorticks=true,
xlabel style={font=\color{white!15!black}},
xlabel={number of subproblems},
ymode=log,
ymin=1e-06,
ymax=1,
yminorticks=true,
ylabel style={font=\color{white!15!black}},
ylabel={distance to $\bar\beta$},
axis background/.style={fill=white},
title style={font=\bfseries,font = \fontsize{14}{14}\selectfont,thick},
title={\textbf{Convergence to} $\bar\beta$; $F_2$}
]
\addplot [color=mycolor1, forget plot]
  table[row sep=crcr]{%
22	0.25496115287285\\
27	0.25496115287285\\
63.0000000000001	0.25496115287285\\
63.0000000000001	0.436030422828491\\
77.0000000000001	0.436030422828491\\
118	0.436030422828491\\
118	0.423868514993799\\
130	0.423868514993799\\
130	0.18452163708661\\
281	0.18452163708661\\
470.000000000001	0.18452163708661\\
629	0.18452163708661\\
814.000000000001	0.18452163708661\\
1081	0.18452163708661\\
1081	0.0558785066532805\\
1375	0.0558785066532805\\
1620	0.0558785066532805\\
1620	0.0589238874624517\\
1969	0.0589238874624517\\
2322	0.0589238874624517\\
2322	0.022060405750059\\
2561	0.022060405750059\\
2897	0.022060405750059\\
3059	0.022060405750059\\
3059	0.024079080223414\\
3426	0.024079080223414\\
3426	0.0179474167081606\\
3851	0.0179474167081606\\
3851	0.00701098202609818\\
4357	0.00701098202609818\\
4357	0.00837578528639266\\
4880	0.00837578528639266\\
4880	0.00740879339345492\\
5220	0.00740879339345492\\
5220	0.000850825277129078\\
5623	0.000850825277129078\\
5889.00000000001	0.000850825277129078\\
6028.00000000001	0.000850825277129078\\
6325.00000000001	0.000850825277129078\\
6325.00000000001	0.00169242624943998\\
6464	0.00169242624943998\\
6464	0.000882800105580021\\
6691	0.000882800105580021\\
6866	0.000882800105580021\\
6866	0.000787806223671348\\
7072.00000000001	0.000787806223671348\\
7244	0.000787806223671348\\
7244	0.000308824619745108\\
7392	0.000308824619745108\\
7500	0.000308824619745108\\
7653	0.000308824619745108\\
7653	0.000135950724276567\\
7802.00000000001	0.000135950724276567\\
7802.00000000001	9.29938725291881e-05\\
7937.00000000001	9.29938725291881e-05\\
7997	9.29938725291881e-05\\
8080	9.29938725291881e-05\\
8175	9.29938725291881e-05\\
8274	9.29938725291881e-05\\
8405	9.29938725291881e-05\\
8493	9.29938725291881e-05\\
8621	9.29938725291881e-05\\
8621	4.53220661328392e-05\\
8732.00000000001	4.53220661328392e-05\\
8732.00000000001	6.77982574836372e-05\\
8843.00000000001	6.77982574836372e-05\\
8843.00000000001	3.28054762183901e-05\\
8976	3.28054762183901e-05\\
9093.00000000001	3.28054762183901e-05\\
9093.00000000001	6.51366268062722e-06\\
9233	6.51366268062722e-06\\
9382	6.51366268062722e-06\\
9547.00000000001	6.51366268062722e-06\\
9547.00000000001	1.26264029951876e-05\\
9724	1.26264029951876e-05\\
9892.00000000001	1.26264029951876e-05\\
9892.00000000001	5.80484505029229e-06\\
10096	5.80484505029229e-06\\
10096	3.31577828821729e-06\\
10338	3.31577828821729e-06\\
10338	5.12591089002427e-06\\
10606	5.12591089002427e-06\\
10886	5.12591089002427e-06\\
11170	5.12591089002427e-06\\
11170	3.53879680394758e-06\\
11526	3.53879680394758e-06\\
12061	3.53879680394758e-06\\
12641	3.53879680394758e-06\\
13269	3.53879680394758e-06\\
13269	2.73409594650632e-06\\
14111	2.73409594650632e-06\\
15252	2.73409594650632e-06\\
15252	2.57210579414765e-06\\
16557	2.57210579414765e-06\\
17947	2.57210579414765e-06\\
19116	2.57210579414765e-06\\
20073	2.57210579414765e-06\\
20889	2.57210579414765e-06\\
21554	2.57210579414765e-06\\
22078	2.57210579414765e-06\\
22502	2.57210579414765e-06\\
22841	2.57210579414765e-06\\
23118	2.57210579414765e-06\\
23333	2.57210579414765e-06\\
23525	2.57210579414765e-06\\
23667	2.57210579414765e-06\\
23775	2.57210579414765e-06\\
23875	2.57210579414765e-06\\
23941	2.57210579414765e-06\\
24005	2.57210579414765e-06\\
24081	2.57210579414765e-06\\
24153	2.57210579414765e-06\\
24177	2.57210579414765e-06\\
};
\addplot [color=mycolor2, forget plot]
  table[row sep=crcr]{%
22	0.860226496649195\\
27	0.860226496649195\\
63.0000000000001	0.860226496649195\\
77.0000000000001	0.860226496649195\\
118	0.860226496649195\\
130	0.860226496649195\\
281	0.860226496649195\\
470.000000000001	0.860226496649195\\
629	0.860226496649195\\
814.000000000001	0.860226496649195\\
1081	0.860226496649195\\
1375	0.860226496649195\\
1620	0.860226496649195\\
1969	0.860226496649195\\
2322	0.860226496649195\\
2561	0.860226496649195\\
2897	0.860226496649195\\
3059	0.860226496649195\\
3426	0.860226496649195\\
3851	0.860226496649195\\
4357	0.860226496649195\\
4880	0.860226496649195\\
5220	0.860226496649195\\
5220	0.252865131075219\\
5623	0.252865131075219\\
5623	0.0545473018072054\\
5889.00000000001	0.0545473018072054\\
5889.00000000001	0.0425427385947336\\
6028.00000000001	0.0425427385947336\\
6325.00000000001	0.0425427385947336\\
6325.00000000001	0.036111315959981\\
6464	0.036111315959981\\
6464	0.0231086872313216\\
6691	0.0231086872313216\\
6866	0.0231086872313216\\
6866	0.0172429405218277\\
7072.00000000001	0.0172429405218277\\
7072.00000000001	0.0138818476485538\\
7244	0.0138818476485538\\
7392	0.0138818476485538\\
7392	0.00740879339345492\\
7500	0.00740879339345492\\
7653	0.00740879339345492\\
7653	0.00616958914208737\\
7802.00000000001	0.00616958914208737\\
7802.00000000001	0.00253747851449872\\
7937.00000000001	0.00253747851449872\\
7997	0.00253747851449872\\
8080	0.00253747851449872\\
8175	0.00253747851449872\\
8274	0.00253747851449872\\
8405	0.00253747851449872\\
8493	0.00253747851449872\\
8621	0.00253747851449872\\
8621	0.00111546045980278\\
8732.00000000001	0.00111546045980278\\
8843.00000000001	0.00111546045980278\\
8976	0.00111546045980278\\
9093.00000000001	0.00111546045980278\\
9233	0.00111546045980278\\
9382	0.00111546045980278\\
9547.00000000001	0.00111546045980278\\
9724	0.00111546045980278\\
9724	0.00054557853978408\\
9892.00000000001	0.00054557853978408\\
9892.00000000001	0.000510661833043947\\
10096	0.000510661833043947\\
10338	0.000510661833043947\\
10338	0.000135950724276567\\
10606	0.000135950724276567\\
10886	0.000135950724276567\\
11170	0.000135950724276567\\
11526	0.000135950724276567\\
12061	0.000135950724276567\\
12641	0.000135950724276567\\
12641	8.34952647173841e-05\\
13269	8.34952647173841e-05\\
14111	8.34952647173841e-05\\
14111	3.4129472014382e-05\\
15252	3.4129472014382e-05\\
15252	2.89114868205822e-05\\
16557	2.89114868205822e-05\\
17947	2.89114868205822e-05\\
19116	2.89114868205822e-05\\
19116	2.8006429210213e-05\\
20073	2.8006429210213e-05\\
20889	2.8006429210213e-05\\
21554	2.8006429210213e-05\\
22078	2.8006429210213e-05\\
22078	2.77027682564991e-05\\
22502	2.77027682564991e-05\\
22841	2.77027682564991e-05\\
22841	2.60059952086524e-05\\
23118	2.60059952086524e-05\\
23118	1.8624545900078e-05\\
23333	1.8624545900078e-05\\
23525	1.8624545900078e-05\\
23667	1.8624545900078e-05\\
23775	1.8624545900078e-05\\
23875	1.8624545900078e-05\\
23941	1.8624545900078e-05\\
23941	1.74013930051918e-05\\
24005	1.74013930051918e-05\\
24081	1.72636490362484e-05\\
24153	1.72636490362484e-05\\
24153	1.28565967605299e-05\\
24177	1.28565967605299e-05\\
};
\end{axis}
\end{tikzpicture}%
}\hspace{.2cm}%
 {
%
%
\definecolor{mycolor1}{rgb}{0.00000,0.44700,0.74100}%
\definecolor{mycolor2}{rgb}{0.85000,0.32500,0.09800}%
\begin{tikzpicture}[every node/.style={scale=0.7}]
\begin{axis}[%
width=3.2cm,
at={(0in,0in)},
scale only axis,
xmode=log,
xmin=10,
xmax=110028,
xminorticks=true,
xlabel style={font=\color{white!15!black}},
xlabel={number of subproblems},
ymode=log,
ymin=1e-06,
ymax=10,
yminorticks=true,
ylabel style={font=\color{white!15!black}},
ylabel={distance to $\bar\beta$},
axis background/.style={fill=white},
title style={font=\bfseries,font = \fontsize{14}{14}\selectfont,thick},
title={\textbf{Convergence to} $\bar\beta$; $F_3$}
]
\addplot [color=mycolor1, forget plot]
  table[row sep=crcr]{%
28	0.0645095840096478\\
32	0.0645095840096478\\
36	0.0645095840096478\\
40	0.0645095840096478\\
62.0000000000001	0.0645095840096478\\
127	0.0645095840096478\\
242	0.0645095840096478\\
500	0.0645095840096478\\
652.000000000001	0.0645095840096478\\
944	0.0645095840096478\\
1212	0.0645095840096478\\
1464	0.0645095840096478\\
1776	0.0645095840096478\\
2192	0.0645095840096478\\
2192	0.00825958400964788\\
2647	0.00825958400964788\\
3290	0.00825958400964788\\
4008	0.00825958400964788\\
4738	0.00825958400964788\\
5479	0.00825958400964788\\
5479	0.00580291599035209\\
6225	0.00580291599035209\\
7010.00000000001	0.00580291599035209\\
7010.00000000001	0.00122833400964795\\
7756	0.00122833400964795\\
8522.00000000001	0.00122833400964795\\
9330.00000000001	0.00122833400964795\\
10188	0.00122833400964795\\
10992	0.00122833400964795\\
11838	0.00122833400964795\\
12729	0.00122833400964795\\
12729	0.000529478490352008\\
13735	0.000529478490352008\\
14819	0.000529478490352008\\
16149	0.000529478490352008\\
16149	0.000349427759647969\\
17747	0.000349427759647969\\
19518	0.000349427759647969\\
21676	0.000349427759647969\\
21676	9.00253653519645e-05\\
24170	9.00253653519645e-05\\
26998	9.00253653519645e-05\\
30444	9.00253653519645e-05\\
34088	9.00253653519645e-05\\
38084	9.00253653519645e-05\\
38084	1.98379158979911e-05\\
43002	1.98379158979911e-05\\
48256	1.98379158979911e-05\\
53792	1.98379158979911e-05\\
60006	1.98379158979911e-05\\
66718	1.98379158979911e-05\\
73878	1.98379158979911e-05\\
73878	7.62790441455331e-06\\
81619.9999999999	7.62790441455331e-06\\
90043.9999999999	7.62790441455331e-06\\
99390.0000000001	7.62790441455331e-06\\
99390.0000000001	6.1050057417189e-06\\
110028	6.1050057417189e-06\\
};
\addplot [color=mycolor2, forget plot]
  table[row sep=crcr]{%
28	1.03668708491896\\
32	1.03668708491896\\
36	1.03668708491896\\
40	1.03668708491896\\
62.0000000000001	1.03668708491896\\
127	1.03668708491896\\
242	1.03668708491896\\
500	1.03668708491896\\
652.000000000001	1.03668708491896\\
944	1.03668708491896\\
1212	1.03668708491896\\
1464	1.03668708491896\\
1776	1.03668708491896\\
2192	1.03668708491896\\
2192	0.988248032903573\\
2647	0.988248032903573\\
2647	0.960164521754838\\
3290	0.960164521754838\\
4008	0.960164521754838\\
4738	0.960164521754838\\
5479	0.960164521754838\\
5479	0.787033370361282\\
6225	0.787033370361282\\
7010.00000000001	0.787033370361282\\
7010.00000000001	0.604186154125279\\
7756	0.604186154125279\\
8522.00000000001	0.604186154125279\\
9330.00000000001	0.604186154125279\\
9330.00000000001	0.550890077658675\\
10188	0.550890077658675\\
10992	0.550890077658675\\
11838	0.550890077658675\\
12729	0.550890077658675\\
13735	0.550890077658675\\
13735	0.374940925727064\\
14819	0.374940925727064\\
14819	0.32532380271851\\
16149	0.32532380271851\\
16149	0.217196136774853\\
17747	0.217196136774853\\
19518	0.217196136774853\\
19518	0.139573883537674\\
21676	0.139573883537674\\
24170	0.139573883537674\\
26998	0.139573883537674\\
26998	0.131306459009648\\
30444	0.131306459009648\\
30444	0.101423646509648\\
34088	0.101423646509648\\
38084	0.101423646509648\\
43002	0.101423646509648\\
43002	0.0785720840096479\\
48256	0.0785720840096479\\
53792	0.0785720840096479\\
60006	0.0785720840096479\\
60006	0.0623123183846479\\
66718	0.0623123183846479\\
73878	0.0623123183846479\\
81619.9999999999	0.0623123183846479\\
81619.9999999999	0.0475906386971479\\
90043.9999999999	0.0475906386971479\\
99390.0000000001	0.0475906386971479\\
110028	0.0475906386971479\\
};
\end{axis}
\end{tikzpicture}%
}%
 \caption{Distance between the calculated solution 
	 $\bar\beta$ and the best known vertex (blue) and
the furthest active vertex (red) respectively for each iteration.} 
  \label{fig:beta_conv}
\end{figure}
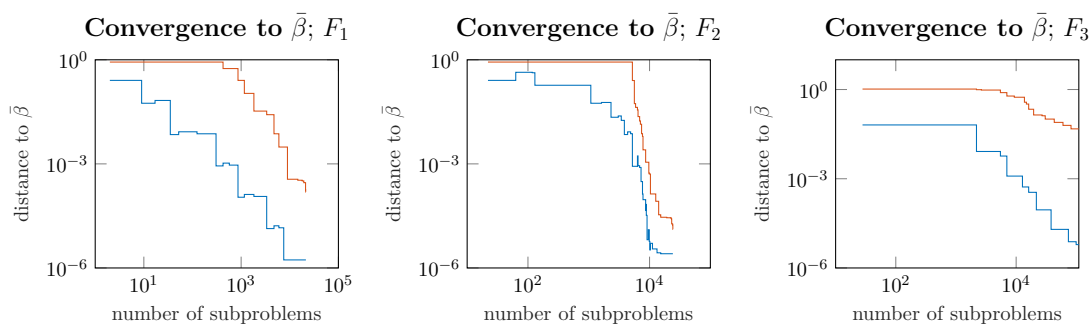
\begin{figure}[ht]
	\centering
  {\input{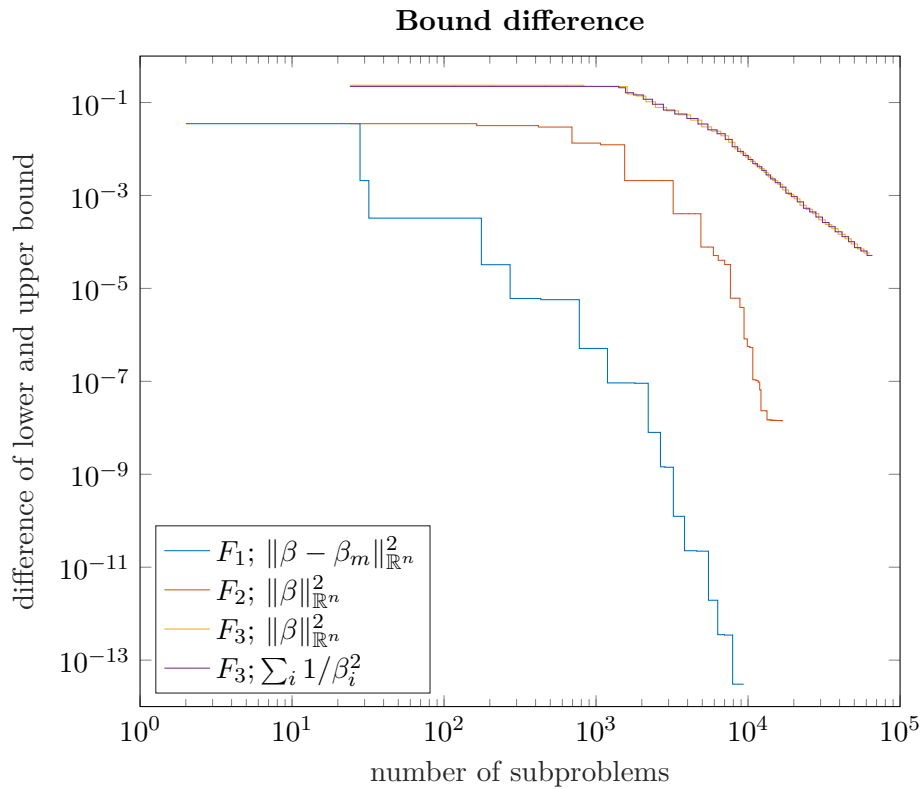}}%
    \caption{Difference of upper
    and lower bound for the settings of $F_1$, $F_2$ and $F_3$ w.r.t. the
  number of solved subproblems. For the setting of $F_3$ the results for two different
regularization terms are displayed}
\label{fig:bound_diffs}
\end{figure}

\begin{figure}
  \centering
  \includegraphics[scale=.18]{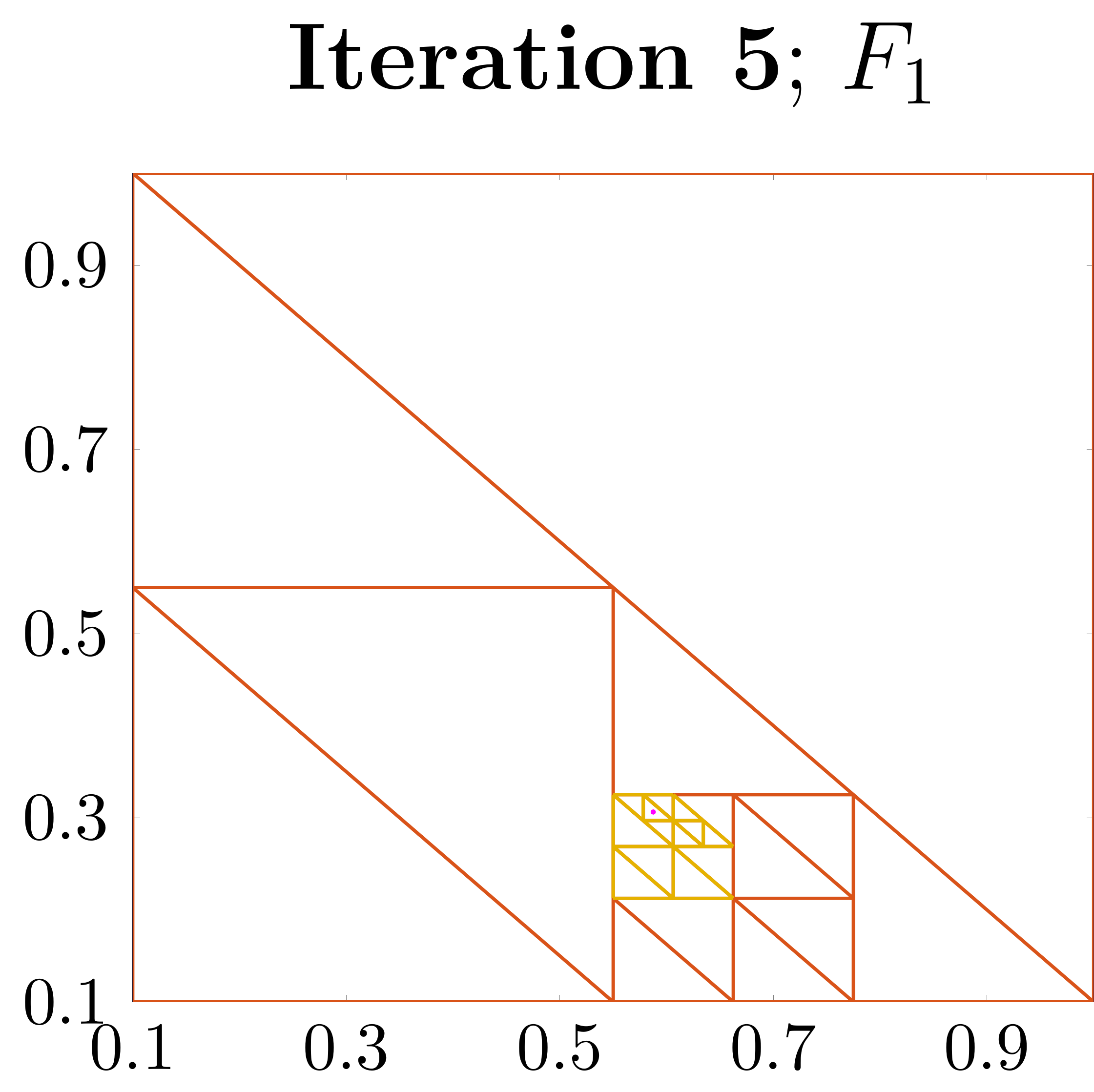}\hspace{.1cm}%
  \includegraphics[scale=.18]{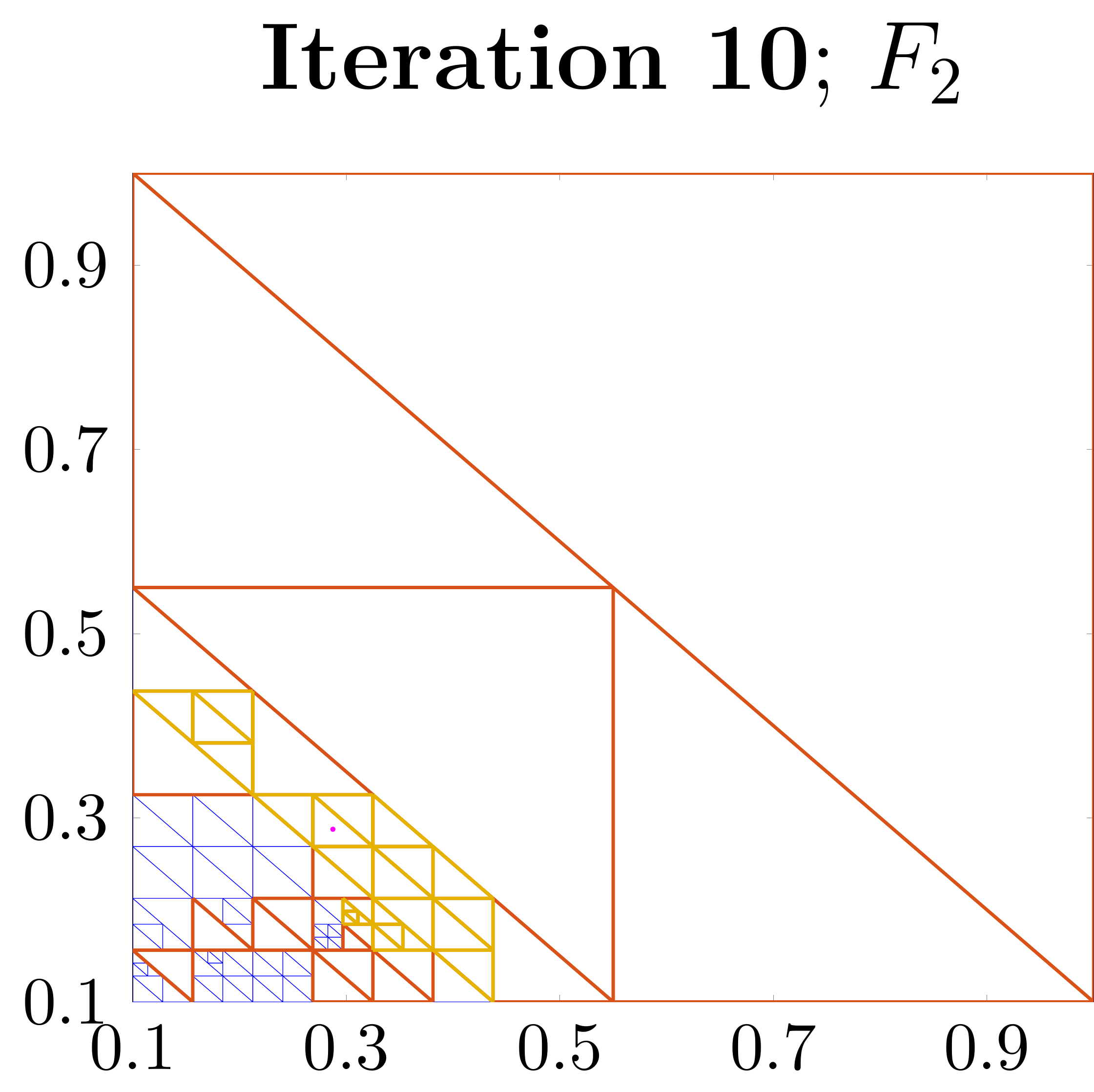}\hspace{.1cm}%
  \includegraphics[scale=.18]{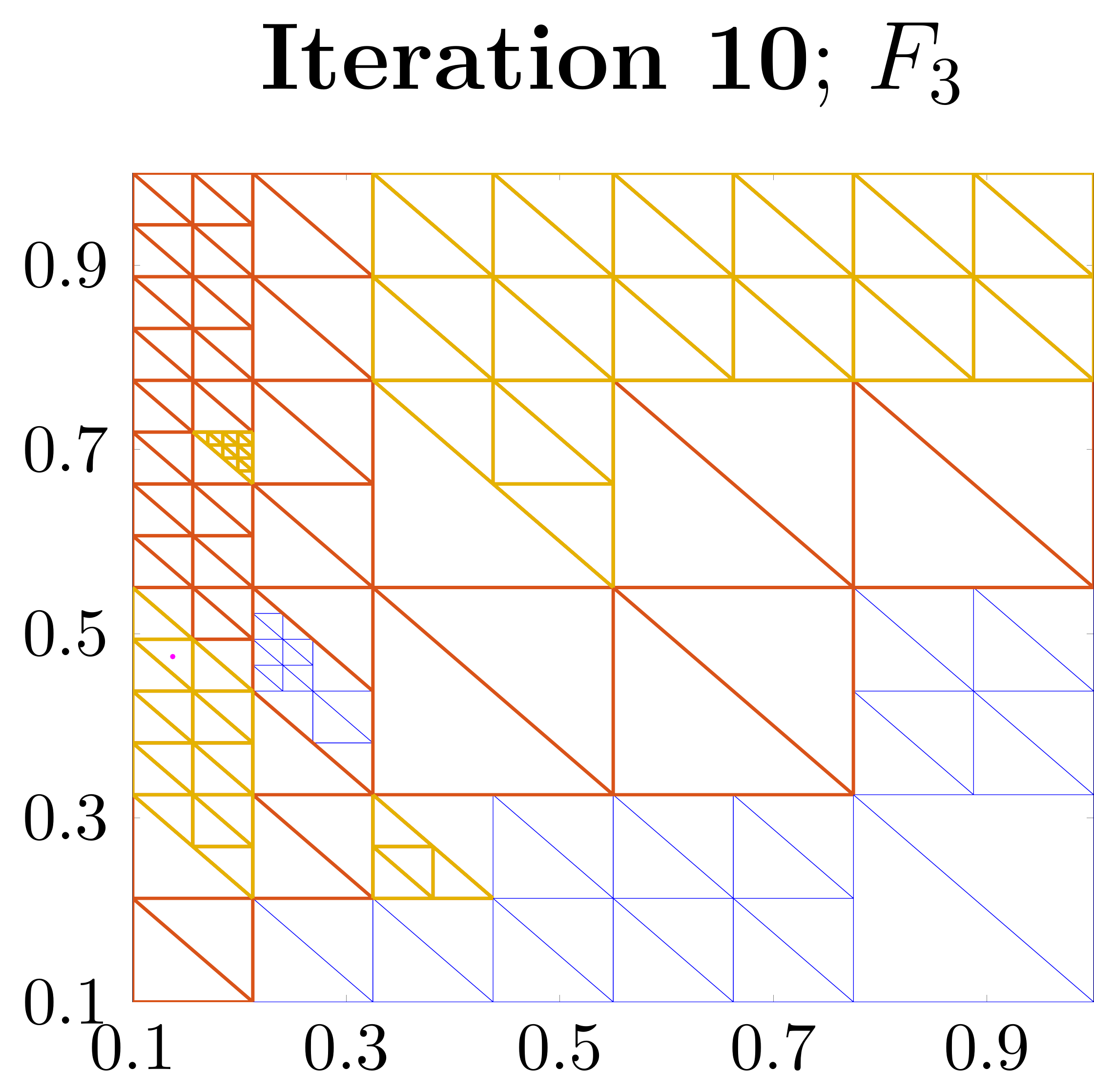}%
  \\%
  \includegraphics[scale=.18]{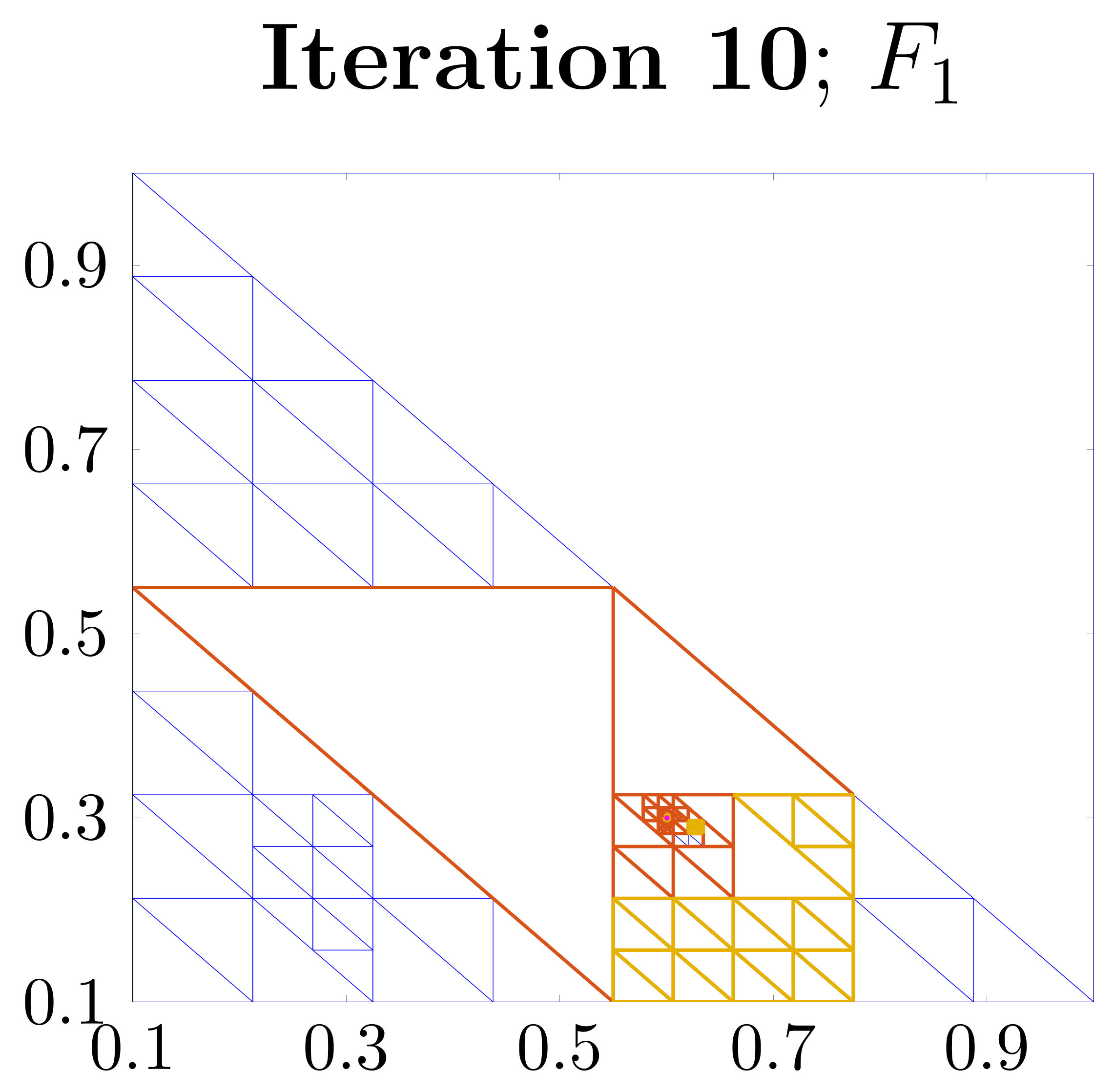}\hspace{.1cm}%
  \includegraphics[scale=.18]{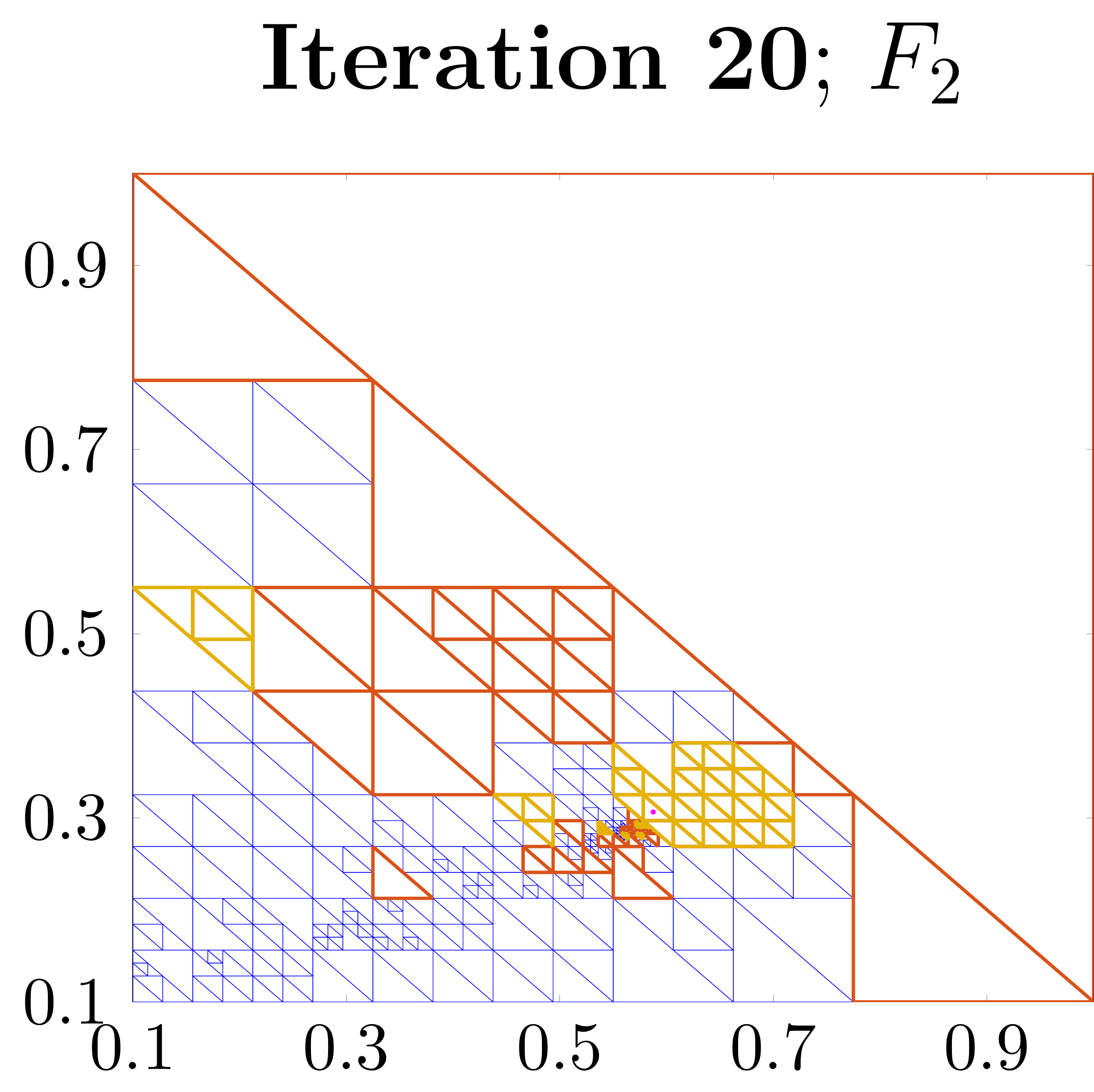}\hspace{.1cm}%
  \includegraphics[scale=.18]{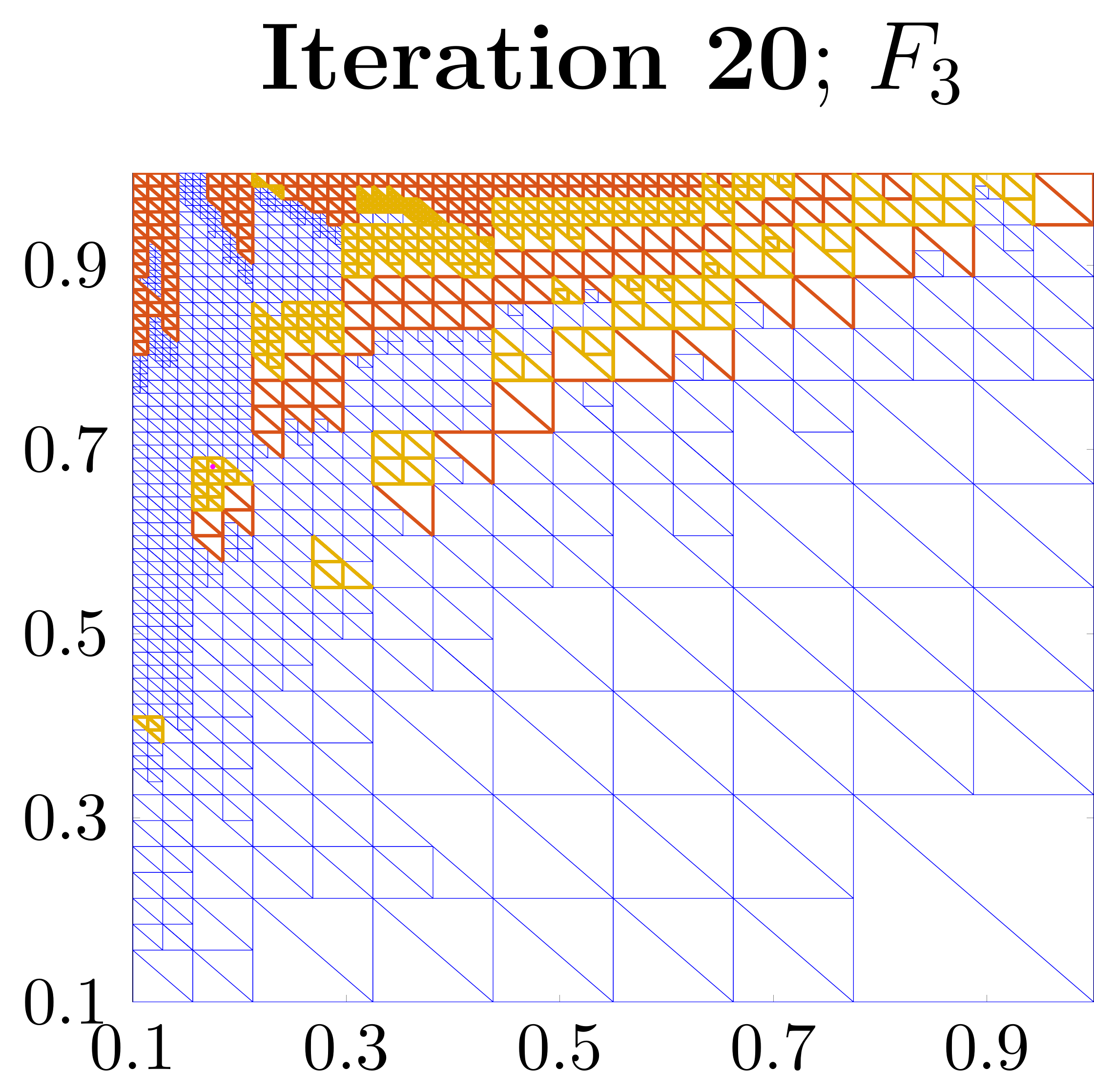}%
  \\%
  \includegraphics[scale=.18]{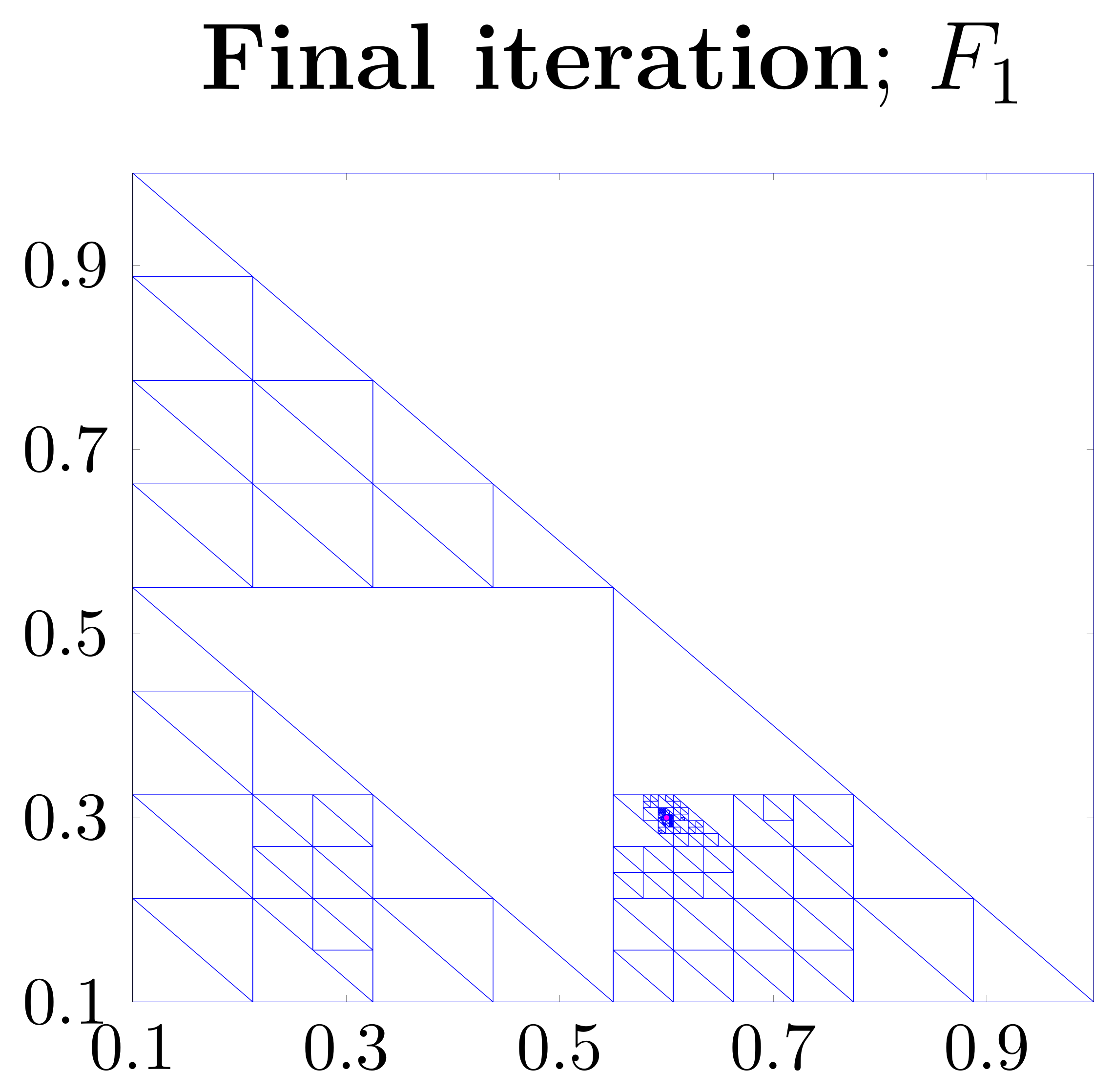}\hspace{.1cm}%
  \includegraphics[scale=.18]{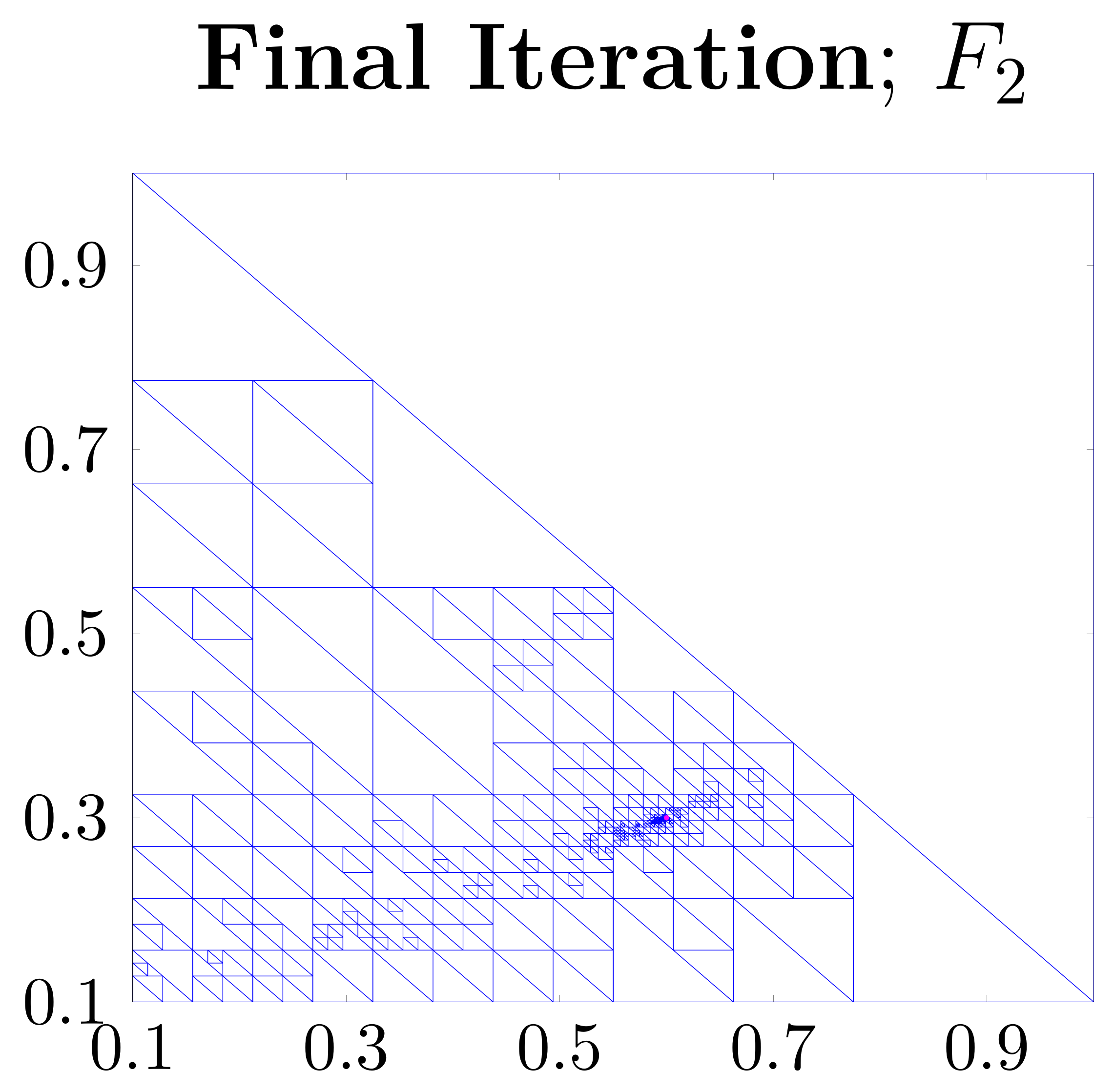}\hspace{.1cm}%
  \includegraphics[scale=.18]{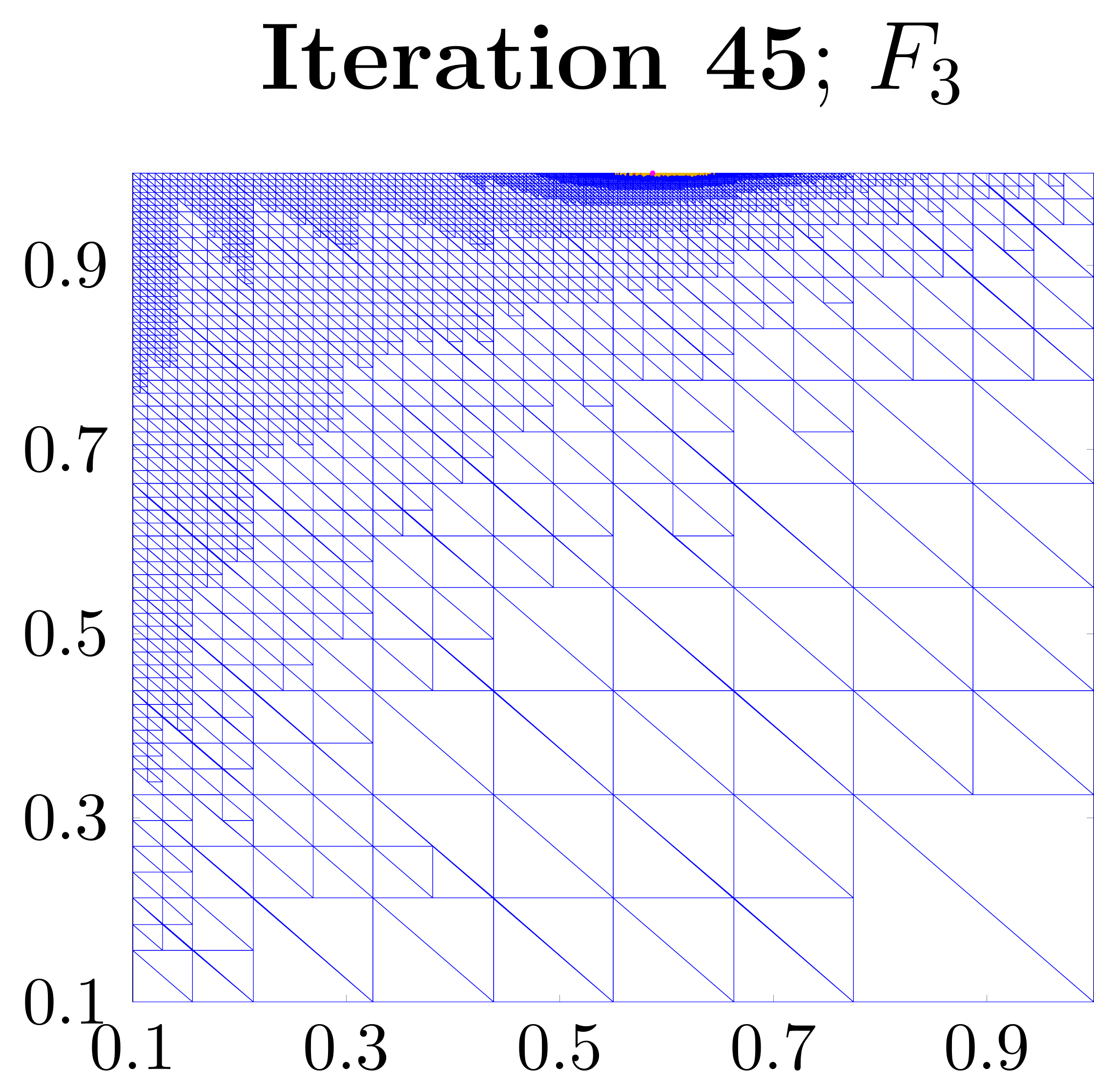}%
  \\%
  \includegraphics[scale=.18]{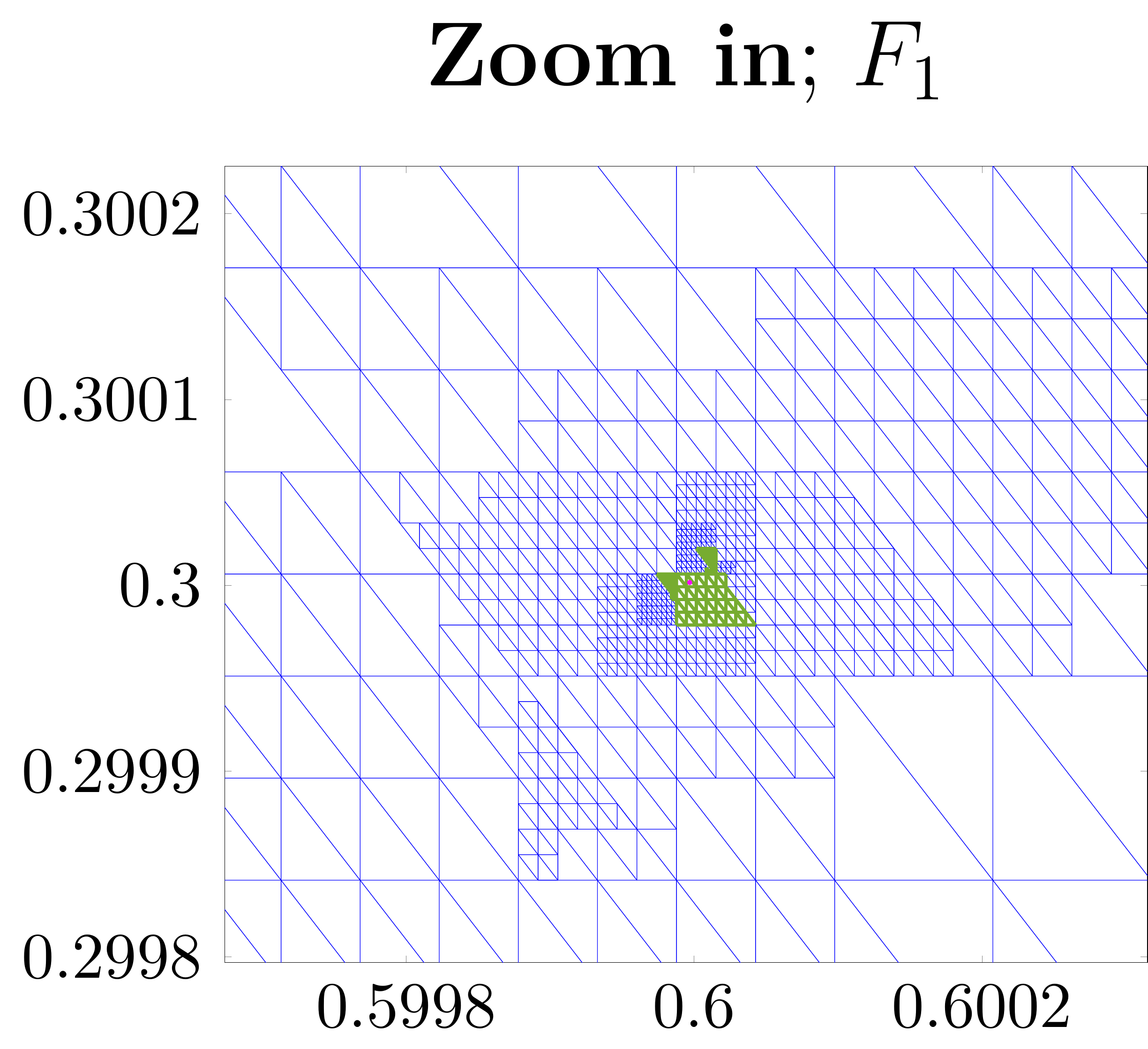}\hspace{.1cm}%
  \includegraphics[scale=.18]{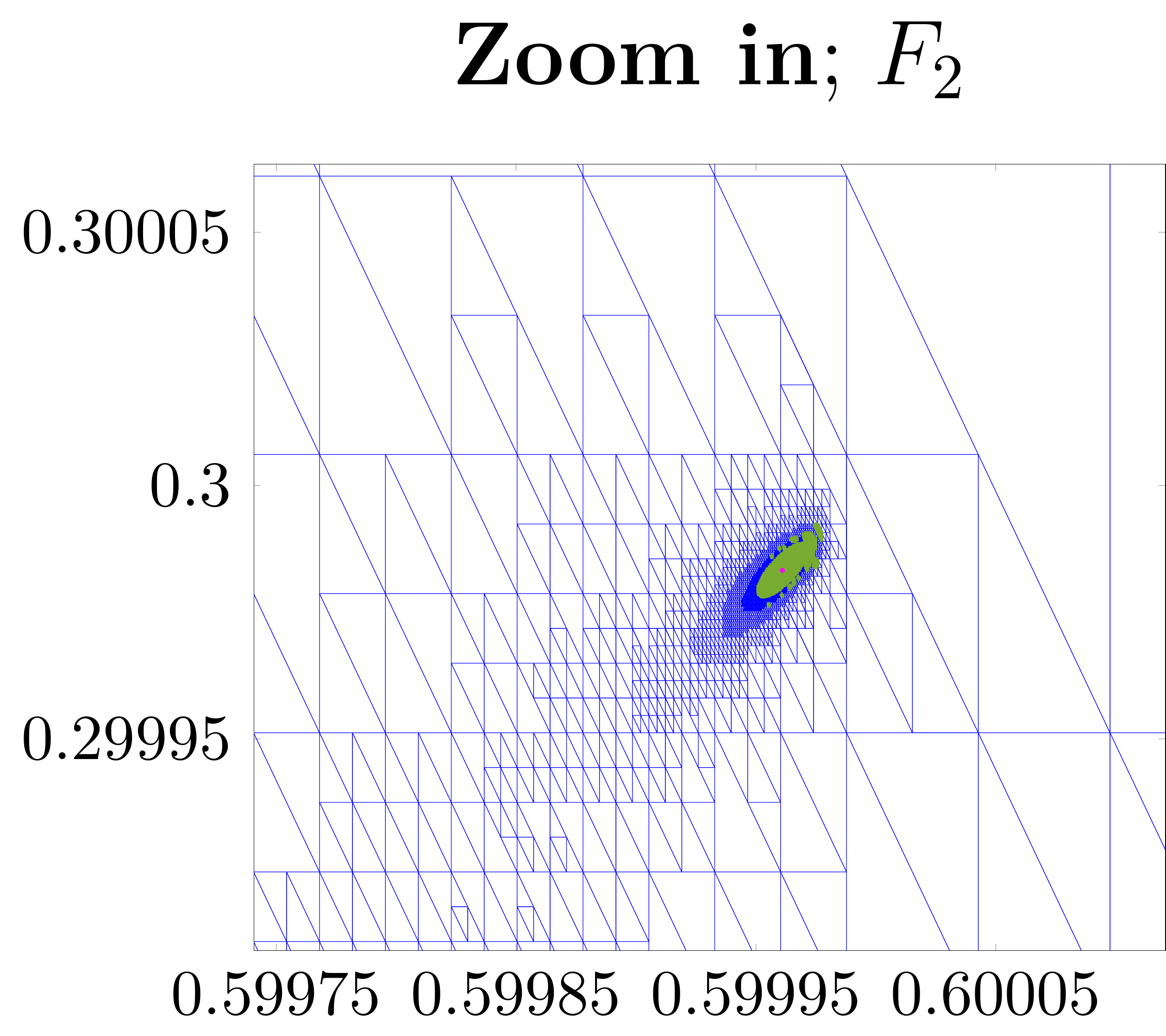}\hspace{.1cm}%
  \includegraphics[scale=.18]{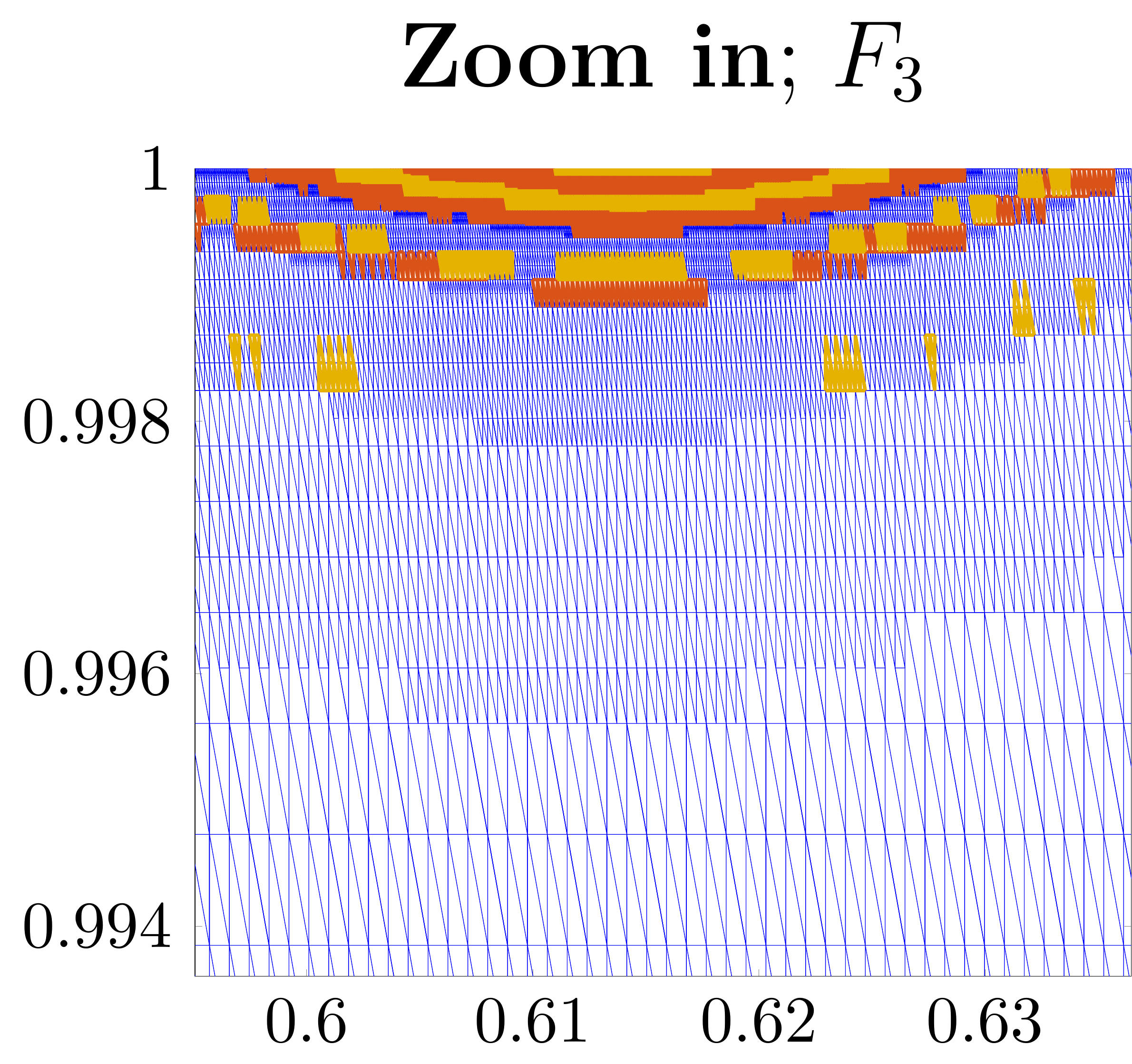}%
 \caption{From left to right: Progression of the splitting of the domain
    $\Bad$ for 
    \hyperref[eq:OVR_affine_penalty]{(OVRP($\xi_\TT$))}
		for the settings of $F_1$, $F_2$ and $F_3$.
  Simplices are differentiate by the color of their outline:
Dismissed (blue), 
relevant(red), 
split in the last iteration (yellow),
difference of lower and upper bound 
for the element is within $10^{-9}$ (green).
The element with the current best objective value is marked with a pink dot.}
  \label{fig:conv_graphic}
\end{figure}

We now visualize the convergence of \cref{alg:global_solution_of_IOC_penalty}
in \cref{fig:lower_upper_bounds}--\ref{fig:conv_graphic}.
These graphics indicate the convergence $\beta_k\rightarrow \bar\beta$ as predicted in 
\cref{lem:penalty_convergence}, see in particular \cref{fig:beta_conv}.
In
\cref{fig:lower_upper_bounds} we show the difference of lower and upper bound
compared for
all mentioned settings.
Note that these bounds are theoretically
nondecreasing, but in the setting of $F_1$ the lower bound in
\cref{fig:lower_upper_bounds} is close to zero
with repsect to machine accuracy, which explains the slightly perturbed behaviour.

We have a stark difference of
convergence speed for the different settings introduced in this section.
Additionally there is a noticeable difference between looking at the vertex that
provides the upper bound and the furthest active vertex. Note that only for the latter
the distance to $\bar\beta$ is guaranteed to be nonincreasing, while the vertex
providing the upper bound might be more interesting from a heuristic
point of view if one considers a depth-search.
The splitting of the domain can be seen in \cref{fig:conv_graphic}. 
For the
purpose of better visualization in the setting of $F_1$ and $F_2$,
the algorithm
was continued for \cref{fig:conv_graphic} until every element either had a vertex for which the corresponding upper level objective was
close ($10^{-9}$) to the upper bound or was dismissed.
We show the difference 
of lower and upper bound 
for all the cases discussed in
\cref{fig:bound_diffs}.

Finally, we give some explanation for the difference in convergence speed. As
discussed in
\cref{convergence_speed_issues} and
\cref{rem:better_small_T}, a growth condition for the upper-level
objective functional for a solution w.r.t.\ $\beta$ allows for an estimate of
convergence speed. This is exactly what we have for the setting of $F_1$. Thus,
we get the estimate from
\cref{convergence_speed_issues} and the number of active subproblems does not
substantially increase between iterations. 
For the case of $F_2$, we have the
second case from \cref{rem:better_small_T}, where the derivative of $F_2$ is close to
zero in the solution. This is, because the term $\norm{\beta}_{\R^n}^2$ only comes up as
a regularization with a small parameter for the upper-level objective
functional. The solution of the parameter estimation problem is still close to
$(y_m,u_m)$. For the case of $F_3$, we no longer have a setting for which we
obtain a nice bound on the number of required subproblems to reach a certain accuracy.
Especially, the number of of active subproblems might heavily increase during
the runtime of \cref{alg:global_solution_of_IOC_penalty}. This can be seen well
in
\cref{fig:conv_graphic}. Finally \cref{fig:bound_diffs} indicates, that the
important property in the setting of $F_3$ is that the solution is no longer close to 
$(\hat y_m,\hat u_m)$, i.e.\ that the target state is ``unreachable'' and that the choice of regularization term
$\frac{\sigma_\beta}{2}\norm{\beta}_{\R^n}^2$ or
$\frac{\sigma_\beta}{2}\sum_{i=1}^2\frac{1}{\beta_i^2}$ is of minor importance
regarding convergence speed for this case.

\subsection*{Data availability}
The datasets generated during and/or analysed during the current study are
available from the corresponding author on reasonable request.

\newrefcontext[sorting=nyt]
\printbibliography

@ARTICLE{AlbrechtLeiboldUlbrich2012,
	author = {Albrecht, Sebastian and Leibold, Marion and Ulbrich, Michael},
	title = {A bilevel optimization approach to obtain optimal cost
		functions for human arm movements},
	journal = {Numerical Algebra, Control and Optimization},
	volume = {2},
	year = {2012},
	number = {1},
	pages = {105--127},
	issn = {2155-3289},
	doi = {10.3934/naco.2012.2.105},
}

@ARTICLE{AlbrechtUlbrich2017,
	author = {Albrecht, Sebastian and Ulbrich, Michael},
	title = {Mathematical programs with complementarity constraints in the
		context of inverse optimal control for locomotion},
	journal = {Optimization Methods \& Software},
	volume = {32},
	year = {2017},
	number = {4},
	pages = {670--698},
	issn = {1055-6788},
	doi = {10.1080/10556788.2016.1225212},
}

@book{Bard1998,
	author="Bard, J. F.",
	title="Practical {B}ilevel {O}ptimization: {A}lgorithms and {A}pplications",
	year=1998,
	publisher="Kluwer Academic",
	address="Dordrecht"
}

@article{BenitaMehlitz2016,
author = {Benita, F. and Mehlitz, P.},
title = {Bilevel {O}ptimal {C}ontrol {W}ith {F}inal-{S}tate-{D}ependent {F}inite-{D}imensional {L}ower {L}evel},
journal = {SIAM Journal on Optimization},
volume = {26},
number = {1},
pages = {718-752},
year = {2016},
doi = {10.1137/15M1015984}
}

@BOOK{BonnansShapiro2000,
	title                    = {Perturbation Analysis of Optimization Problems},
	author                   = {J. Fr{\'{e}}d{\'{e}}ric Bonnans and Alexander Shapiro},
	publisher                = {Springer},
	year                     = {2000},
	doi = {10.1007/978-1-4612-1394-9},
	address                  = {Berlin}
}

@BOOK{BrennerScott2008,
	doi = {10.1007/978-0-387-75934-0},
	year = 2008,
	publisher = {Springer New York},
	author = {Susanne C. Brenner and L. Ridgway Scott},
	title = {The Mathematical Theory of Finite Element Methods}
}

@BOOK{BrezziFortin1991,
	doi = {10.1007/978-1-4612-3172-1},
	year = 1991,
	publisher = {Springer New York},
	editor = {Franco Brezzi and Michel Fortin},
	title = {Mixed and Hybrid Finite Element Methods}
}

@article{DacorognaMarechal2008,
	title = {The role of perspective functions in convexity,  polyconvexity, rank-one convexity and separate convexity},
	author = {Dacorogna, Bernard and Maréchal, Pierre},
	journal = {Journal of Convex Analysis},
	number = {2},
	volume = {15},
	pages = {271-284},
	year = {2008},
}

@BOOK{Dempe2002,
	title                    = {Foundations of bilevel programming},
	author                   = {Stephan Dempe},
	publisher                = {Kluwer Academic Publishers},
	year                     = {2002},
	address                  = {Dordrecht}
}

@ARTICLE{DempeHarderMehlitzWachsmuth2018:1,
	doi = {10.1007/s10898-019-00758-1},
	year = 2019,
	publisher = {Springer Nature},
	author = {Stephan Dempe and Felix Harder and Patrick Mehlitz and Gerd Wachsmuth},
	title = {Solving inverse optimal control problems via value functions to global optimality},
	journal = {Journal of Global Optimization},
	volume = {74},
	number = {2},
	pages = {297--325},
}

@BOOK{DempeKalashnikovPerezValdesKalashnykova2015,
	author = {Dempe, S. and Kalashnikov, V. and P{\'{e}}rez-Vald{\'{e}}z, G. and Kalashnykova, N.},
	title = {Bilevel Programming Problems - Theory, Algorithms and Applications to Energy Networks},
	year=2015,
	publisher = {Springer},
	address = {Berlin}
}

@article{FischLenzHolzapfelSachs2012,
title	="On the {S}olution of {B}ilevel {O}ptimal {C}ontrol {P}roblems  to {I}ncrease the {F}airness in {A}ir {R}aces",
author	="Fisch, F. and Lenz, J. and Holzapfel, F. and Sachs, G.",
year	=2012,
journal	="Journal of Guidance, Control, and Dynamics",
volume	=35,
number	=4,
pages	="1292-1298",
doi		="10.2514/1.54407"
}

@PHDTHESIS{Harder2021,
  doi = {10.26127/BTUOPEN-5375},
  author = {Harder, Felix},
  title = {On bilevel optimization problems in infinite-dimensional spaces},
  school = {BTU  Cottbus-Senftenberg},
  year = {2021},
}

@ARTICLE{HarderWachsmuth2017:2,
	title                    = {Comparison of Optimality Systems for the Optimal Control of the Obstacle Problem},
	author                   = {Felix Harder and Gerd Wachsmuth},
	year                     = {2018},
	journal = {GAMM-Mitteilungen},
	volume = {40},
	number = {4},
	pages = {312--338},
	doi = {10.1002/gamm.201740004},
}

@ARTICLE{HarderWachsmuth2018:1,
	doi = {10.1080/02331934.2018.1495205},
	year = 2018,
	publisher = {Informa {UK} Limited},
	volume = {68},
	number = {2-3},
	pages = {615--643},
	author = {Felix Harder and Gerd Wachsmuth},
	title = {Optimality conditions for a class of inverse optimal control problems with partial differential equations},
	journal = {Optimization}
}

@phdthesis{Hatz2014,
	title		="Efficient {N}umerical {M}ethods for {H}ierarchical {D}ynamic {O}ptimization with {A}pplication to {C}erebral {P}alsy {G}ait {M}odeling",
	author		="Hatz, K.",
	school		="University of Heidelberg, Germany",
	year		=2014
}

@ARTICLE{HatzSchloederBock2012,
	author = {Hatz, Kathrin and Schl\"oder, Johannes P. and Bock, Hans Georg},
	title = {Estimating parameters in optimal control problems},
	journal = {SIAM Journal on Scientific Computing},
	volume = {34},
	year = {2012},
	number = {3},
	pages = {A1707--A1728},
	issn = {1064-8275},
	doi = {10.1137/110823390},
}

@ARTICLE{HintermuellerItoKunisch2002,
	doi = {10.1137/s1052623401383558},
	year = 2002,
	publisher = {Society for Industrial {\&} Applied Mathematics ({SIAM})},
	volume = {13},
	number = {3},
	pages = {865--888},
	author = {Michael Hintermüller and Kazufumi Ito and Karl Kunisch},
	title = {The Primal-Dual Active Set Strategy as a Semismooth Newton Method},
	journal = {{SIAM} Journal on Optimization}
}

@BOOK{HinzePinnauUlbrichUlbrich2009,
	doi = {10.1007/978-1-4020-8839-1},
	year = {2009},
	author = {Michael Hinze and Rene Pinnau and Michael Ulbrich and Stefan Ulbrich},
	publisher = {Springer Netherlands},
	title = {Optimization with {PDE} Constraints}
}

@ARTICLE{HollerKunischBarnard2018,
	doi = {10.1088/1361-6420/aade77},
	year = 2018,
	publisher = {{IOP} Publishing},
	volume = {34},
	number = {11},
	pages = {115012},
	author = {Gernot Holler and Karl Kunisch and Richard C Barnard},
	title = {A bilevel approach for parameter learning in inverse problems},
	journal = {Inverse Problems}
}

@article{KalashnikovBenitaMehlitz2015,
	author	="Kalashnikov, V. and Benita, F. and Mehlitz, P.",
	title	="The natural gas cash-out problem: {A} bilevel optimal control approach",
	journal	="Math. Probl. Eng.",
	year	=2015,
	pages 	="1-17",
	doi		="10.1155/2015/286083"
}

@Incollection{KnauerBueskens2010,
author="Knauer, M. and B{\"u}skens, C.",
editor="Diehl, M. and Glineur, F. and Jarlebring, E. and Michiels, W.",
title="Hybrid {S}olution {M}ethods for {B}ilevel {O}ptimal {C}ontrol {P}roblems with {T}ime {D}ependent {C}oupling",
bookTitle="Recent Advances in Optimization and its Applications in Engineering: The 14th Belgian-French-German Conference on Optimization",
year="2010",
publisher="Springer",
address="Berlin",
pages="237--246",
doi="10.1007/978-3-642-12598-0_20"
}

@book{LewisVrabieSyrmos2012,
	author		="Lewis, F. L. and Vrabie, D. and Syrmos, V. L.",
	title		="Optimal Control",
	year		=2012,
	publisher	="John Wiley \& Sons",
	address		="Hoboken"
}

@PHDTHESIS{Mehlitz2017,
	title                    = {Contributions to complementarity and bilevel programming in {B}anach spaces},
	author                   = {Mehlitz, Patrick},
	school                   = {Technische Universität Bergakademie Freiberg},
	year                     = {2017},
	eprinttype = {urn},
	eprint = {urn:nbn:de:bsz:105-qucosa-227091}
}

@ARTICLE{MehlitzWachsmuth2015:1,
	title                    = {Weak and strong stationarity in generalized bilevel programming and bilevel optimal control},
	author                   = {Patrick Mehlitz and Gerd Wachsmuth},
	journal                  = {Optimization},
	year                     = {2016},
	number                   = {5},
	pages                    = {907--935},
	volume                   = {65},
	doi                      = {10.1080/02331934.2015.1122007}
}

@INBOOK{MehlitzWachsmuth2019:1,
	author = {Patrick Mehlitz and Gerd Wachsmuth},
	title = {Bilevel optimal control: existence results and stationarity conditions},
	editor = {Dempe, Stephan and Zemkoho, Alain},
	booktitle = {Bilevel Optimization: Advances and Next Challenges},
	year = {2020},
	publisher = {Springer International Publishing},
	address = {Cham},
	pages = {451--484},
	doi = {10.1007/978-3-030-52119-6_16},
}

@ARTICLE{Outrata1990,
	doi = {10.1007/bf01416737},
	year = {1990},
	month = {jul},
	publisher = {Springer Science and Business Media {LLC}},
	volume = {34},
	number = {4},
	pages = {255--277},
	author = {Jiří V. Outrata},
	title = {On the numerical solution of a class of Stackelberg problems},
	journal = {{ZOR} Zeitschrift für Operations Research Methods and Models of Operations Research},
}

@ARTICLE{Robinson1976:1,
	title                    = {Stability theory for systems of inequalities. {II}. {D}ifferentiable
		nonlinear systems},
	author                   = {Robinson, Stephen M.},
	journal                  = {SIAM Journal on Numerical Analysis},
	year                     = {1976},
	number                   = {4},
	pages                    = {497--513},
	volume                   = {13},
	doi = {10.1137/0713043},
	issn                     = {0036-1429}
}

@book{ShimizuIshizukaBard1997,
	author="Shimizu, K. and Ishizuka, Y. and Bard, J. F.",
	title="Nondifferentiable and two-level mathematical programming",
	year=1997,
	publisher="Kluwer Academic",
	address="Dordrecht"
}

@book{Troeltzsch2009,
	author		="Tr{\"o}ltzsch, F.",
	title		="Optimale Steuerung partieller Differentialgleichungen",
	year		=2009,
	publisher	="Vieweg",
	address		="Wiesbaden"
}

@book{Troutman1996,
	author	="Troutman, J. L.",
	title	="Variational Calculus and Optimal Control",
	year	=1996,
	publisher	="Springer",
	address	="New York"
}

@BOOK{Ulbrich2011,
	author = {Ulbrich, Michael},
	title                    = {Semismooth Newton Methods for Variational Inequalities and Constrained
		Optimization Problems in Function Spaces},
	series = {MOS-SIAM Series on Optimization},
	volume = {11},
	publisher = {Society for Industrial and Applied Mathematics (SIAM),
		Philadelphia, PA; Mathematical Optimization Society,
		Philadelphia, PA},
	year = {2011},
	pages = {xiv+308},
	isbn = {978-1-611970-68-5},
	doi = {10.1137/1.9781611970692},
}

@article{Ye1995,
author = {Ye, J. J.},
title = {Necessary {C}onditions for {B}ilevel {D}ynamic {O}ptimization {P}roblems},
journal = {SIAM Journal on Control and Optimization},
volume = {33},
number = {4},
pages = {1208-1223},
year = {1995},
doi = {10.1137/S0363012993249717}
}

@article{Ye1997,
author = {Ye, J. J.},
title = {Optimal {S}trategies {F}or {B}ilevel {D}ynamic {P}roblems},
journal = {SIAM Journal on Control and Optimization},
volume = {35},
number = {2},
pages = {512-531},
year = {1997},
doi = {10.1137/S0363012993256150}
}

@ARTICLE{ZoweKurcyusz1979,
	title                    = {Regularity and stability for the mathematical programming problem
		in {B}anach spaces},
	author                   = {Zowe, Jochem and Kurcyusz, Stanisław},
	journal                  = {Applied Mathematics and Optimization},
	year                     = {1979},
	number                   = {1},
	pages                    = {49--62},
	doi = {10.1007/BF01442543},
	volume                   = {5}
}

\end{document}